\documentclass[11pt]{article}
\title{The largest eigenvalue of small rank perturbations of Hermitian random matrices }
\author{S. P\'ech\'e \\ \\Institut Fourier, Universit\'e Joseph Fourier,\\
 BP 74, 38402 St MARTIN D'HERES Cedex, France,\\
sandrine.peche@ujf-grenoble.fr}
\usepackage{epsfig}
\usepackage{amsfonts}
\usepackage{amssymb}
\usepackage{amsthm}
\usepackage{amstext}
\usepackage{amscd}
\usepackage{amsmath}
\usepackage{delarray}
\begin{document}
\maketitle
\newtheorem{theo}{Theorem}[section]
\newtheorem{prop}{Proposition}[section]
\newtheorem{lemme}{Lemma}[section]
\newtheorem{conjecture}{Conjecture}[section]
\newtheorem{definition}{Definition}[section]
\newtheorem{fact}{Fact}[section]
\newtheorem{hyp}{Assumption}[section]
\theoremstyle{remark}
\newtheorem*{rem}{Remark}
\newtheorem{remark}{Remark}[section]
\newtheorem{Remark}{Remark}[section]
\newtheorem{Notationnal remark}{Remark}[section]
\newcommand{\bremnot}{\begin{Notationnal remark}}
\newcommand{\eremnot}{\end{Notationnal remark}}
\newcommand{\brem}{\begin{remark}}
\newcommand{\erem}{\end{remark}}
\newcommand{\bconj}{\begin{conjecture}}
\newcommand{\econj}{\end{conjecture}}
\newcommand{\bdefi}{\begin{definition}}
\newcommand{\edefi}{\end{definition}}
\newcommand{\bt}{\begin{theo}}
\newcommand{\bfa}{\begin{fact}}
\newcommand{\efa}{\end{fact}}
\newcommand{\Si}{\Sigma}
\newcommand{\et}{\end{theo}}
\newcommand{\bp}{\begin{prop}}
\newcommand{\ep}{\end{prop}}
\newcommand{\bl}{\begin{lemme}}
\newcommand{\el}{\end{lemme}}
\newcommand{\be}{\begin{equation}}
\newcommand{\ee}{\end{equation}}
\newcolumntype{L}{>{$}l<{$}}
\newenvironment{Cases}{\begin{array}\{{lL.}}{\end{array}}
\begin{abstract}
We compute the limiting eigenvalue statistics at the edge of the spectrum of large Hermitian random matrices perturbed by the addition of small rank deterministic matrices. We consider random Hermitian matrices with independent Gaussian entries $M_{ij}, i\leq j$ with various expectations.  We prove that the largest eigenvalue of such random matrices exhibits, in the large $N$ limit, various limiting distributions depending on both the eigenvalues of the matrix $\left(\mathbb{E}M_{ij}\right)_{i,j=1}^N$ and its rank. This rank is also allowed to increase with $N$ in some restricted way.
\end{abstract}
\section{Introduction and results}
The aim of this paper is to investigate how a small rank perturbation of a standard $N\times N$ random matrix can affect significatively the limiting properties of the spectrum, as the size $N$ of the matrix goes to infinity. The statistics of extreme eigenvalues is here of interest. Note that it is not clear what is meant by ``a small rank perturbation of a random matrix'' and we shall define it formally in the sequel. Actually, a first study of eigenvalue statistics for such perturbed random matrices has been achieved in \cite{BaikGBAPeche}. Therein the authors consider non homogeneous Wishart random matrices $R_N=1/N XX^*$, where $X$ is a $p\times N$ random matrix with independent complex Gaussian entries with a \emph{spiked} covariance matrix $\Si.$ That is, $\Si-Id$ ($Id$ is the identity matrix) is a fixed rank (independent of $N$) diagonal matrix, while both $p$ and $N$ go to infinity. \\In this paper, we consider Hermitian random matrices. Let $\mu$ (resp. $ \mu'$) be a probability distribution on $\mathbb{C}$ (resp. $\mathbb{R}$). A $N\times N$ random Hermitian matrix $M(\mu,\mu')$ is then a Hermitian matrix with entries being mutually independent random variables of distribution $\mu$ (resp $\mu'$) strictly above the diagonal (resp. on the diagonal).
Define then $\displaystyle{M_N(\mu,\mu')=\frac{1}{\sqrt N}M(\mu,\mu')}$.
Let also $\lambda_1\geq  \lambda_2\geq \cdots \geq \lambda_N $ be the ordered eigenvalues of $M_N$ and $\mu_N=\dfrac{1}{N}\sum_{i=1}^N \delta_{ \lambda_i}$ its spectral measure. A famous result of Wigner (\cite{Wigner}) asserts that $\mu_N$ admits a non-random limit as $N$ goes to infinity.
\bp \cite{Wigner} Assume that $\int x d\mu(x)=\int x d\mu'(x)=0,$ and that
$\int |x|^2 d\mu(x)=\sigma^2,$ $\int |x|^2d\mu'(x)<\infty.$ Then, almost surely, $\displaystyle{\lim_{N \rightarrow \infty}\mu_N=\tilde\rho_{\sigma}, }$ where $\tilde \rho_{\sigma} $ is the semi-circular law with parameter $\sigma^2$, defined by the density with respect to Lebesgue measure \be\label{semicerclesigma^2}\rho_{\sigma}(x)=\frac{2}{\pi \sigma^2}\sqrt{4\sigma^2-x^2}1_{[-2\sigma,2\sigma]}(x).\ee
\ep
%This fundamental result is a so-called "universality result", since the limiting semi-circular law does not depend on the details of the distributions $\mu, \mu'$, provided they are centered and of finite variance.
 %In view of this global behavior of the spectrum, it is then natural to study the behavior of extreme eigenvalues.
 %Note that Wigner's Theorem does not give much insight about the limiting behavior of $\lambda_1$.
 Let $\lambda^*=2\sigma$ be the top edge of the support of $\tilde\rho_{\sigma}.$
 It is then a fundamental result of \cite{Geman}
  that, for the archetypical of Hermitian ensemble, the so-called GUE, $\displaystyle{\lim_{N\rightarrow\infty }\lambda_1=\lambda^*.}$
\bdefi The $N\times N$ GUE with parameter $\sigma^2$ is the distribution of a $N\times N$ random matrix $ M(\mu,\mu')$,
if $\mu$ (resp. $\mu'$) is the centered complex (resp. real) Gaussian distribution of variance $\sigma^2$.
\edefi
The result obtained in \cite{Geman} has later been precised in \cite{TracyWidom}.
Consider the Airy function defined by
$\displaystyle{Ai(u)=\frac{1}{2\pi}\int_{\infty e^{i5\pi/6}}^{\infty e^{i\pi/6}} \exp{\{iua+\frac{1}{3}a^3\}}da,}$
and define the Airy kernel
\be \label{noyauAiry}Ai(u,v)=\int_{0}^{\infty} Ai(y+u)Ai(y+v)dy.\ee
\bdefi The Tracy-Widom distribution is defined by the distribution function \\$ F^{TW}_2(x):=\det(I-A_x),$ where $A_x$ is the trace class operator acting on $L^2(x, \infty)$ with kernel $Ai(u,v).$
 \edefi

\bp \cite{TracyWidom}% Let $\mu$ be any symmetric distribution satisfying the following conditions: $\displaystyle{\int|x|^2d\mu=\sigma^2}$ and there exists $C>0$ such that $$\forall k\geq 0, \: \int |x|^{2k}d\mu'<(Ck)^k, \:\int |x|^{2k}d\mu<(Ck)^k.$$
Let $\lambda_1$ be the largest eigenvalue of $V_N=\frac{1}{\sqrt N}V$, where $V$ is drawn from the GUE with parameter $\sigma^2$. Then,
$\displaystyle{\lim_{N \rightarrow \infty}P\left ( \sigma^{-1}N^{2/3}\left (\lambda_{1}- \lambda^*\right )\leq x\right)=  F^{TW}_2(x). }$
  \ep
\brem It is shown in \cite{Soshnikovedge} that the above result actually holds for a  wide class of random matrices $M_N(\mu,\mu')$ with centered distributions $\mu,$ $\mu'$.
\erem
The scope of this paper is to define a suitable "small" rank perturbation of a random matrix $V_N$ drawn from the GUE, so that the largest eigenvalue separates from "the bulk", $[-\lambda^*,\lambda^*],$ and study in this case, how it interacts with the ``bulk" of eigenvalues in $[-\lambda^*,\lambda^*].$
Due to the rotational invariance of the Gaussian distribution, it is enough to consider  diagonal perturbations.
%The scope of this paper will then be to determine a critical scale, if any, for both the rank and the entries of $W_N$, for which $\lambda_1$ separates from the "bulk", in the large $N$ limit.
\subsection{The model}
The model studied here is known in random matrix litterature as the deformed Wigner ensemble. The first study of such an ensemble goes back to \cite{BH1} and \cite{Johansson}.
\bdefi Given $k \in \mathbb{N}$, $r\in \mathbb{N}$ and ordered real numbers $\pi_1>\pi_2\geq\cdots\geq \pi_{r+1}$, a \emph{deformed Wigner matrix} is a $N \times N$ random matrix $M_N=W_N+\frac{1}{\sqrt N} V$ where $V$ is of the $N\times N$ GUE with parameter $1$ and $W_N$ is the diagonal matrix $W_N=\text{ diag }(\pi_1,\ldots,\pi_1,\pi_2,\ldots, \pi_{r+1},0,\ldots, 0),$ with rank $k+r,$ and
 %with  % $\pi_1> \pi_2 \geq\cdots\geq \pi_{r+1},$ and
 where the largest eigenvalue $\pi_1$ has multiplicity $k$.
\edefi
\brem We assume that $\pi_i=0,$ $\forall i\geq 2$ if $r=0.$ The $\pi_i, i=1,\ldots, r+1$ can be negative but lie in a compact set independent of $N$.
\erem

In this paper, we consider matrices $W_N$ with rank $k+r$ such that \be \label{condrang}\lim_{N \rightarrow \infty}\dfrac{k+r}{N}=0.\ee
In particular, $k$ and $r$ may depend on $N$ . %$\displaystyle{\lim_{N \rightarrow\infty}\dfrac{1}{N}\sum_{i=1}^{N}\delta _{W_N_{ii}} =\delta_{0}.}$
Noting $\lambda_1\geq \lambda_2\geq \cdots \geq \lambda_N $ the ordered eigenvalues of $M_N$ and $\mu_N$ its spectral measure, condition (\ref{condrang}) ensures that
 $\displaystyle{\lim_{N \rightarrow \infty}\mu_N=\tilde \rho_1}$, where $\tilde \rho_1$ is the semi-circle law defined in (\ref{semicerclesigma^2}), with parameter $\sigma^2=1$.

\subsection{Results}First, we fix the rank of $W_N$ independently of $N$ and identify the critical scale $\pi_1=\pi_1^c$ for which $\lambda_1$ separates from the bulk. Results in this part are similar to those in \cite{BaikGBAPeche}. Then, and this is the main result of the paper, we study the limiting properties of largest eigenvalues when the rank of $W_N$ is allowed to increase with $N$, focusing on the case where $\lambda_1$ is separated from the bulk.
\subsubsection{A fixed rank perturbation}
We consider matrices $W_N$ with fixed rank $k+r$, independent of $N$.
\begin{hyp}\label{Hyp1} $W_N=\text{diag}(\pi_1,\ldots, \pi_1,\pi_2,\ldots,\pi_{r},0,\ldots,0),$ with $\pi_1$ of multiplicity $k,$ such that \begin{itemize}
 \item $k$ and $r$ are given integers \emph{independent of $N,$}
 \item $\pi_1$ is a given real number independent of $N$,
 \item $\pi_i, i=2,\ldots, r+1$ lie in a compact set of $(-\infty,\pi_1)$ \emph{independent of $N$}.
\end{itemize}
\end{hyp}
%The results obtained here are similar to those in \cite{BaikGBAPeche}.
Before stating the results, we need a few definitions.
%Let us define the generalized Airy kernel for $k\geq 0$.
Given an integer $m \geq 1,$ and a contour ${\cal C}$ going from $\infty e^{5i\pi/6}$ to $\infty e^{i\pi/6}$, with $0$ lying above ${\cal C}$, we set
\begin{equation}
t^{(m)}(v)=\frac{1}{2\pi}\int_{{\cal C}}
\exp{\{iua+\frac{1}{3}a^3i\}}(-ia)^{m-1}da,
\quad
s^{(m)}(u)=\frac{1}{2\pi}\int_{{\cal C}}
\exp{\{iua+\frac{1}{3}a^3i\}}\frac{1}{(ia)^m}da.\label{smBaik}
\end{equation}
 Given $x\in\mathbb{R}$, let also $A_x$ be the operator acting on $L^2(x, \infty)$ with kernel $Ai(u,v)$ defined
in (\ref{noyauAiry}), and $<\:,\:>$ denote the standard scalar product of operators on $L^2(x, \infty).$
 \bdefi Given an integer $k\geq 0$, $F^{TW}_{k+2}$ is the distribution function defined by \be
\label{airygal}F^{TW}_{k+2}(x)=\det(1-A_x)\det\left (\delta_{m,n}-<\frac{1}{1-A_x}s^{(m)}, t^{(n)}>\right )_{1\leq
m,n\leq k}, \: x\in \mathbb{R}. \ee
\edefi
\brem $F^{TW}_{k+2}$ was proved to be distribution function in
\cite{BaikGBAPeche}.
\erem

\paragraph{} The first theorem gives a necessary condition to have
$\displaystyle{\lim_{N \rightarrow \infty }\lambda_1=\lambda^*=2.}$ Still, we prove that the limiting distribution of $\lambda_1$ depends on both the value and the multiplicity of $\pi_1.$

 \bt \label{Th: UnigalWig}
 Assume Assumption \ref{Hyp1} holds.  \begin{itemize}\item If $\pi_1<1,$ then,
 $\displaystyle{\lim_{N \rightarrow \infty}P \left(N^{2/3}\left ( \lambda_{1}-2\right ) \leq x\right )= F^{TW}_2(x).}$
\item If $\pi_1=1$, then,
$\displaystyle{\lim_{N \rightarrow \infty}P \left(N^{2/3}\left ( \lambda_{1}-2\right ) \leq x\right )= F^{TW}_{k+2}(x).}$
\end{itemize}
\et
In the next theorem, we prove that, as soon as $\pi_1>1$, with probability one, the largest eigenvalue $\lambda_1$ exits the support of the semi-circular law.

\bdefi  Given $k\geq 0,$ define the probability distribution
$$F_{GUE,\sigma^2}^{k}(x)=\frac{1}{Z_k}\int_{-\infty}^x \cdots \int_{-\infty}^x \prod_{1\leq i<j\leq k}|u_i-u_j|^2\prod_{i=1}^k \exp{\{-\frac{u_i^2}{2\sigma^2}\}} du_1\cdots du_k,$$
where $Z_k$ is the normalizing constant $Z_k=\displaystyle{\int_{\mathbb{R}^k}\prod_{1\leq i<j\leq k}|u_i-u_j|^2\prod_{i=1}^k \exp{\{-\frac{u_i^2}{2\sigma^2}\}} du_1\cdots du_k.}$
\edefi
\brem It can be shown (see e.g. \cite{Mehta}, Chapter 5) that $F_{GUE,\sigma^2}^{k}$ is the probability distribution of the largest eigenvalue of the $k\times k $ GUE with parameter $\sigma^2.$
\erem

\bt \label{Th: GUEWig}
Assume Assumption \ref{Hyp1} holds with $\pi_1>1.$ Then,
$$\lim_{N \rightarrow \infty}P \left (\sigma^2(\pi_1) N^{1/2}\left ( \lambda_{1}-C(\pi_1)\right )\leq x\right )= F_{GUE, \sigma^2(\pi_1)}^{k}(x),\text{ where}$$
 \begin{equation}C(\pi_1)=\pi_1+\dfrac{1}{\pi_1}\quad\text{and}\quad
 \label{sigma^2}\sigma^2(\pi_1)=\dfrac{\pi_1^2}{\pi_1^2-1}.\end{equation}
\et
\brem
This result should be compared with the result of \cite{KomlosFuredi}. Therein, the authors 
consider Hermitian random matrices $M_N(\mu,\mu')$, where $\mu$, $\mu'$ are distributions with compact support such that
$\displaystyle{\int x d\mu = \int x d\mu'=m\not=0,}$ $\displaystyle{\int|x|^2d\mu=\int|x|^2d\mu'=\sigma^2+m^2.}$ Then, for $C(\cdot)$ defined as in (\ref{sigma^2}), it is proved that
% $$\sqrt N \lambda_1= \sum_{i,j=1}^N \frac{M_{ij}}{N}+\frac{\sigma^2}{2m}+O(\frac{1}{\sqrt N}).$$
%which implies that
$\sqrt N \left (\lambda_1-C(\sqrt N m)\right )$ has asymptotically Gaussian fluctuations ${\cal N}(0, \sigma^2)$. %In our setting, this would correspond to $W_N=\text{ diag }(m\sqrt N, 0, \ldots, 0)$ (which does not satisfy Assumption \ref{Hyp1}). Then, $\sigma^2=1\simeq \sigma^2(\sqrt N m)$, for large $N$, where $\sigma^2(\cdot)$ is defined in (\ref{sigma^2}).
Here, we obtain that the scale at which $\lambda_1$ actually separates from the bulk (when $\mu$, $\mu'$ are Gaussian distributions) is $m=m_N =\dfrac{1}{\sqrt N}$. \erem
Theorem \ref{Th: GUEWig} gives the intuition that a "bulk" of $k$ eigenvalues exits the support of the semi-circular law, provided $\pi_1>1.$ Furthermore, these $k$ eigenvalues seem to behave as those of a typical $k\times k$ random matrix. We now show that this still holds if $k$ goes to infinity in some restricted way.

\subsubsection{A large rank perturbation }We investigate the case where the rank of $W_N$ is increasing with $N$. To our knowledge, the kind of perturbation that we now define, is new.
%has not been studied yet.
Let $k_N, r_N$ be given sequences of integers such that
\begin{equation}\lim_{N \rightarrow \infty}k_N=\infty, \text{  } \lim_{N\rightarrow\infty}\frac{k_N}{N} = 0,\text{ and }
 \lim_{N\rightarrow\infty}\frac{r_N}{N} = 0.\label{condH2suites}\end{equation}
We first consider the case where $\pi_1>1$, so that the largest eigenvalue separates from the \emph{bulk}.
\begin{hyp}\label{Hyp2} $W_N=\text{ diag }(\pi_1, \ldots, \pi_1,\pi_2,\ldots, \pi_{r_N+1},0,\ldots, 0),$ with $\pi_1$ of multiplicity $k_N$ and  \begin{itemize}
\item $(k_N)_{ N \in \mathbb{N}}$ and $ (r_N)_{ N \in \mathbb{N}}$ satisfy (\ref{condH2suites}),                                                                \item $\pi_1>1$ is given, independent of $N$,
\item $\pi_i, i=2,\ldots, r_N+1$ lie in a compact set of $(-\infty,\pi_1),$ independent of $N$.
\end{itemize}
\end{hyp}
%From the preceding results, there should exist a nucleus of $k_N$ eigenvalues around $C(\pi_1)$, separated from the bulk $[-2,2]$.
%We actually show that, suitably rescaled, these $k_N$ largest eigenvalues behave as the eigenvalues of a typical $k_N\times k_N$ random Hermitian matrix,  in the "bulk" as well as at the edge.\\
We first deal with local eigenvalue statistics in the ''bulk'' of the $k_N$ largest eigenvalues and consider the so-called spacing function between nearest neighbor eigenvalues. Let $\rho=\rho_{\sigma^2}$ be the density of the semi-circular law (\ref{semicerclesigma^2}) with parameter $\sigma^2(\pi_1),$ defined in (\ref{sigma^2}).
Define \be \label{parametrealpha}\alpha_N=\frac{\sqrt{k_N}}{\sqrt N} \text{ and } \beta_N=\frac{r_N}{N}.\ee
Let $t_N$ be a sequence such that $\displaystyle{\lim_{N \rightarrow \infty}t_N =\infty, \lim_{N \rightarrow \infty}\frac{t_N}{ k_N}= 0}.$
\bdefi \label{def: spacing}Given $|\alpha|<2\sigma(\pi_1)$, and for $\displaystyle{u=C(\pi_1)+\alpha_N\frac{\alpha }{ \sigma^2(\pi_1)}+\frac{\beta_N}{\pi_1}-\frac{1}{N}\sum_{i=1}^{N \beta_N}\dfrac{1}{\pi_1-\pi_{i+1}}}$, the "spacing function", $S_N(\alpha ,s,\lambda),$ is the symmetric function which, if $\lambda_1\leq \lambda_2\leq \cdots \leq\lambda_N$, equals 
$$\displaystyle{S_N(\alpha, s,\lambda)=\frac{1}{2t_N}\sharp \{1\leq j\leq N-1; \quad \lambda_j-\lambda_{j+1} \leq\frac{\alpha_N s}{\sigma^2k_N\rho (\alpha)}, \quad|\lambda_j -u|\leq \frac{\alpha_N t_N}{k_N \rho(\alpha)\sigma^2}\}}.$$
\edefi
Here $\alpha$ is to be seen as a point in the "bulk" of (\ref{semicerclesigma^2}). Then, we obtain the following result.

\bt \label{theo: bulkaubord} Assume Assumption \ref{Hyp2} holds. Then,
 $\displaystyle{
\lim_{N \rightarrow \infty}\left |ES_N(\alpha ,s,\lambda)-\int_0^s H''(u)du\right |=0,}$
where $\displaystyle{H(s)=\sum_{m=0}^{\infty}\frac{(-1)^m}{m!}\int_{[0,s]^m}\text{det}\Bigl (\frac{\sin \pi(x_i-x_j)}{\pi(x_i-x_j)}\Bigr )_{i,j=1}^m \prod_{i=1}^mdx_i.}$
\et
\brem
 The above theorem states that the archetypical repulsion of eigenvalues of Hermitian random matrices is exhibited amongst the $k_N$ largest eigenvalues, in the large $N$ limit.\erem

\brem The case $\pi_1\leq 1$ has already been studied in \cite{Johansson} and \cite{GBAPeche} (Appendix A), showing a similar repulsion of eigenvalues (up to changes in the rescalings).%In this case, given a point $|u|<2$ and a sequence $t_N$ such that $\displaystyle{\lim_{N\rightarrow \infty}\frac{t_N}{N}=0}$, one considers the spacing function
  %  $$S_N(u, s,\lambda)=\frac{1}{2t_N}\sharp \{1\leq j\leq N-1; \: \lambda_{j}-\lambda_{j+1} \leq\frac{s}{N\rho_1 (u)}, |\lambda_j -u|\leq \frac{t_N}{N \rho_1(u)}\}.$$ It is then proved that (\ref{spacingequ}) also holds in this case if $S_N(\alpha ,s,\lambda)$ is replaced with $S_N(u, s,\lambda)$.
                   \erem

\paragraph{}We then turn to local eigenvalue statistics at the edge.
%We need some preliminary definitions to state the results.
Let $\alpha_N$, $\beta_N$ be given as in (\ref{parametrealpha}), and $\log$ be the principal branch of the logarithm. Set, for $w\in \mathbb{C}\setminus (-\infty,\pi_1],$
\be
F_u(w):=w^2/2-uw+(1-\alpha_N^2-\beta_N)\log w +\alpha_N^2 \log(w-\pi_1)+\frac{1}{N}\sum_{i=1}^{N\beta_N}\log(w-\pi_{i+1}), \label{defindefuengen}\ee
so that
\begin{eqnarray} 
&&F_u'(w)=w-u+\frac{1-\alpha_N^2-\beta_N}{w}+\frac{\alpha_N^2}{w-\pi_1}+\frac{1}{N}\sum_{i=1}^{N\beta_N}\frac{1}{w-\pi_{i+1}},\label{defdeuo} \\
& &F_u''(w)=1-\frac{1-\alpha_N^2-\beta_N}{w^2}-\frac{\alpha_N^2}{(w-\pi_1)^2}-\frac{1}{N}\sum_{i=1}^{N\beta_N}\frac{1}{(w-\pi_{i+1})^2}.\label{defw_o}
\end{eqnarray}
Note that $F_{u}''$ does not depend on $u$. We then define $w_o$ as follows.
\be w_o \text{ is the largest solution of the equation }F_u''(w)=0.\label{defwophrase}\ee
In particular, it can be shown that $w_o>\pi_1$.
Finally define $u_o$ and $t_r$ by \begin{equation}\label{defdeu_o}
F'_{u_o}(w_o)=0, \quad t_r=\dfrac{w_o-\pi_1}{\alpha_N},
\end{equation}
where $F_{u}'$ and $w_o$ are respectively given by (\ref{defdeuo}) and (\ref{defwophrase}).
 \bt \label{theo: rhoN<<3/7} Assume Assumption \ref{Hyp2} holds and let $u_o$ and $t_r$ be given by (\ref{defdeu_o}). Then,
$$\lim_{N \rightarrow \infty}P\left(t_r \frac{k_N^{2/3}}{\alpha_N}\left( \lambda_{1}-u_o\right )\leq x\right )= F^{TW}_2(x).$$
%where $F^{TW}_2$ is the Tracy-Widom distribution.
\et

\brem The above theorem states that, as long as $\displaystyle{\alpha_N \rightarrow 0}$, the suitably scaled largest eigenvalue of the deformed Wigner ensemble also behaves as the largest eigenvalue of a $k_N \times k_N$ GUE.
The rescaling is such that $\displaystyle{t_r \frac{k_N^{2/3}}{\alpha_N}}=N^{2/3}\left (\frac{F_{u_o}^{(3)}(w_o)}{2}\right )^{-1/3}(1+o(1))$ and, if $r_N=N\beta_N=0$,
 $\displaystyle{u_o=C(\pi_1)+\alpha_N\frac{2}{\sigma(\pi_1)}+O(\alpha_N^2)}.$ \erem
 \brem The case $\displaystyle{\lim_{N \rightarrow \infty}\frac{k_N}{N}=\alpha \in (0,1)}$ will be
  the object of a subsequent paper, and is not examined here. In this context,
  the limiting statistics of extreme eigenvalues are determined in \cite{BleherKuij}, when
 $W_N=\text{diag }(a,\ldots, a,-a,\ldots,-a)$, where numbers of $a$ and $-a$ are both approximately $N/2.$ \erem

The proof of Theorem \ref{theo: rhoN<<3/7} is based on an extension of the method developed in
\cite{BaikGBAPeche} and may bring some new tools for the study of such deformed models. In
particular, we can also consider the case where $\pi_1\leq 1$.
%, for which $\lambda_1$ may not separate from the bulk. 
Then , we obtain the following result.
 \begin{hyp}\label{Hyp2'} $W_N=\text{ diag }(\pi_1, \ldots, \pi_1,\pi_2,\ldots, \pi_{r_N+1},0,\ldots, 0),$ with $\pi_1$ of multiplicity $k_N$ and  \begin{itemize}
  \item $(k_N)_{ N \in \mathbb{N}},$ and $ (r_N)_{ N \in \mathbb{N}},$ satisfy (\ref{condH2suites}), 
  \item $\pi_i, i=2,\ldots, r+1$ lie in a compact set of $(-\infty,\pi_1),$ independent of $N,$                              
   \item $\pi_1\leq 1$ is given, independent of $N$.                                                                                                \end{itemize}

 \end{hyp}
 Define $F_{u_o}$ as in (\ref{defdeuo}). Let then $w_o$ (greater than 1 here) be given as in (\ref{defwophrase}) and $u_o$ as in (\ref{defdeu_o}).
\bt \label{theopi_1<1} Assume Assumption \ref{Hyp2'} holds.
Then,
$$\lim_{N \rightarrow \infty}P\left ( N^{2/3}\left (\frac{F_{u_o}^{(3)}(w_o)}{2}\right )^{-1/3}(\lambda_1-u_o)\leq x\right )=F^{TW}_2(x).$$
\et
\brem If $\pi_1<1$ and $r_N=0,$ Theorem \ref{theopi_1<1} proves in particular that $\lambda_1$ exhibits the archetypical behavior of the largest eigenvalue of a $N \times N$ GUE with parameter $1$, as long as $k_N <<N^{1/3}.$ Otherwise $\lambda_1$ is slightly translated.
\erem
\brem The proofs of Theorem \ref{theo: rhoN<<3/7} and Theorem \ref{theopi_1<1} are very similar and the second one will only be sketched.
\erem
\subsection{Sketch of the proof}
 Basically, the idea is to deduce the large $N$ limit of local eigenvalue statistics of the deformed Wigner ensemble from the asymptotics of the so-called ``$m$ point correlation functions", defined as follows.
 %correlation functions.
%We here study the local eigenvalue statistics of the deformed Wigner ensemble, by computing
Let $P_N$ be the joint eigenvalue distribution on $(\mathbb{R}^N, {\cal B}(\mathbb{R}^N))$ induced by the deformed Wigner ensemble. It is known that $P_N$ admits a density with respect to Lebesgue measure. We denote this density $g$.
Then, given an integer $m\leq N$,
 the $m$-point correlation function, $R_N^m(\cdot),$ induced by $P_N$ is defined by
$\displaystyle{R_N^m(x_1,\ldots,x_m)=\frac{N!}{(N-m)!}\int_{\mathbb{R}^m} g(x_1,\ldots,
x_N)\prod_{i=m+1}^N dx_i.}$
We refer to \cite{Johansson}, Section 4) for the use of correlation functions in the study of local eigenvalue statistics.
%The reason why we work with correlation functions is that (see e.g. \cite{Johansson}, Section 4) the distribution of the largest eigenvalue, as well as the spacing functions, can be conveniently expressed in terms of the correlation functions.\\

It happens that, for the deformed Wigner ensemble, the computation of the asymptotics of correlation functions is quite simple. This follows from
 beautiful results of \cite{Johansson}, \cite{BH1},\cite{IZ}, \cite{HC}. %based on the famous Itzykson-Zuber-Harisch-Chandra integral formula (see \cite{IZ}, \cite{HC}).
%Indeed, in \cite{BH1} and \cite{Johansson}, it is proved that the deformed Wigner ensemble induces a determinantal random point field, that is, the $m$-point correlation functions of the induced joint eigenvalue density are given by a determinant involving some correlation kernel $K_N$, as we now recall.
\bp \label{Prop: corrJoh}
 \cite{Johansson}The $m$-point correlation function of the deformed Wigner ensemble is given by
$\displaystyle{R_N^m(x_1,\ldots, x_m)=\det \Bigl ( K_N(x_i,x_j)\Bigr )_{i,j=1}^m,}$ with the correlation kernel $K_N$ defined by
\be \label{noyauJohansson}
K_N(u,v)=\frac{N}{(2i\pi)^2}\int_{\Gamma} dz \int_{\gamma}dw e^{N\{\frac{w^2}{2}-vw-\frac{z^2}{2}+uz\}}\left (\frac{w}{z}\right)^{N-k-r}\left (\frac{w-\pi_1}{z-\pi_1}\right)^k\:\prod_{i=2}^{r+1} \frac{w-\pi_i}{z-\pi_i}\frac{1}{w-z},
\ee
where $\Gamma $ encircles $0$ and $\pi_i, i=1,\ldots, r+1,$ and is oriented counterclockwise, and $\gamma =A+i\mathbb{R}$, with $A$ large enough to ensure that $\Gamma\cap \gamma= \emptyset,$ is oriented from bottom to top.
\ep
\brem Actually, the integral representation (\ref{noyauJohansson}) has been established in the case where $W_N$ has pairwise distinct eigenvalues $W_{ii}, i=1, \ldots, N$. Yet, by a straightforward use of l'Hopital's rule, one can see this formula also holds in the case where $W_{jj}=W_{kk},$ for some $j\not=k$.
\erem
Thanks to the above expression (\ref{noyauJohansson}), the asymptotic expansion of $K_N$ can be computed through a saddle point analysis. We then deduce the asymptotic expansion of correlation functions $R_N^m(\cdot)$ and of local eigenvalue statistics. 
Let us develop some of the ideas used for computing the limiting distribution of the largest eigenvalue.
By an inclusion-exclusion formula, one can show that
$\displaystyle{P(\lambda_1\leq s)=\det (I-K_N)_{L^2(s,\infty)},}$
where $\displaystyle{\det (I-K_N)_{L^2(s,\infty)}}$ is the Fredholm determinant of the trace class operator acting on $L^2(s, \infty)$ with kernel $K_N.$
First, we prove that the correlation kernel can be written as
\be\label{formeproduit}K_N(x,y)=\int_0^{\infty}H_N(x+t)J_N(y+t)dt,\ee for some kernels $H_N$, $J_N.$
Using a saddle point analysis, we then prove that there exist kernels $H_{\infty}$, $J_{\infty}$ such that
$$\lim_{N\rightarrow \infty}\int_0^{\infty}|H_N(x+u)-H_{\infty}(x+u)|^2du = 0,\: \:
\lim_{N\rightarrow \infty}\int_0^{\infty}|J_N(x+u)-J_{\infty}(x+u)|^2du = 0,$$
for all $x$ in a compact interval. This ensures that
$\displaystyle{\lim_{N \rightarrow \infty}\det (I-K_N)_{L^2(s, \infty)}=\det (I-K_{\infty})_{L^2(s, \infty)},}$
where $K_{\infty}(x,y)=\int_0^{\infty}H_{\infty}(x+t)J_{\infty}(y+t)dt.$ This eventually gives the convergence in distribution of the largest eigenvalue of the deformed Wigner ensemble.

\section{Proof of Theorem \ref{Th: UnigalWig} \label{Sec: TW}}

\begin{figure}[htbp]
 \begin{center}
  \epsfig{figure=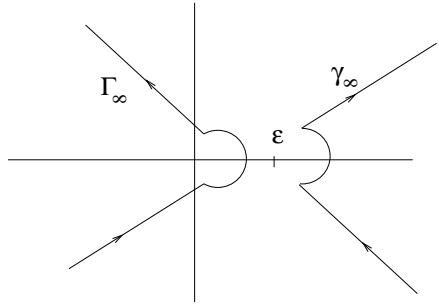,height=4cm, angle=0}
 \caption{Contours $\Gamma_{\infty}$ and $\gamma_{\infty}$.
 \label{fig:contourgammasinfinis}}
 \end{center}
\end{figure}

In this section, we assume that all the eigenvalues of $W_N$ are smaller than,
 or equal to one. We further assume (in this section only) that $\pi_1=1$
 and has multiplicity $k$. It is assumed that $k=0$ if all the eigenvalues
 of $W_N$ are strictly smaller than $1.$ In all cases, we assume that
 $\pi_i<1-\eta,$ for $i=2,\ldots, r+1$ where $\eta >0$ is fixed; $k$ and $r$ are here given
integers independent of $N$.\\
Let then some $\epsilon >0$, that will be fixed later, be given and set
\begin{equation} \label{lesrescalingsTh1.1}
u=2+\frac{x}{N^{2/3}},\:\: v=2+\frac{y}{N^{2/3}},\:\: w_c=1,\:\: \tilde
w_c=w_c+\dfrac{\epsilon}{N^{1/3}}, \:\:K_N'(x,y)=\dfrac{e^{N(u-v)\tilde w_c}}{N^{2/3}}K_N(u,v).
\end{equation}
Note that the rescaled correlation kernel
$K_N'(x,y)$ defines the same correlation functions as \\$\dfrac{1}{N^{2/3}}K_N(u,v).$ Define also
\be  g(w)=\prod_{i=2}^{r+1}\frac{w-\pi_i}{w}\frac{1}{w^k},\quad w\in \mathbb{C}^*,
\text{ and }F(z)=z^2/2-2z+\log z, \quad z\in \mathbb{C}\setminus \mathbb{R}_-. \label{defdeF}
\ee
 Here we use the principal branch of the logarithm and $\exp{\{N \log w\}}$ stands for $w^N.$
Then, from (\ref{noyauJohansson}), we readily obtain that $K_N'(x,y)$  can be cast to the form
(\ref{formeproduit}). Let $\Gamma$ and $\gamma$ be as in (\ref{noyauJohansson}). \bp \label{Prop : removalofsing}
$\displaystyle{K_N'(x,y)=-\int_0^{\infty}H_N(x+t)J_N(y+t)dt,}$ where \be
\label{HN}H_N(x)=\frac{N^{1/3}}{2\pi}
\int_{\Gamma}\frac{1}{g(z)(z-w_c)^k}\exp{\{-NF(z)\}}\exp{\{N^{1/3}(x+t)(z-\tilde w_c)\}}dz, \ee \be
\label{JN}J_N(y)=\frac{N^{1/3}}{2\pi}\int_{\gamma}g(w)(w-w_c)^k\exp{\{NF(w)\}}\exp{\{-N^{1/3}(y+t)(w-\tilde
w_c)\}}dw. \ee \ep
\paragraph{Proof: }We use the fact that $\displaystyle{\frac{1}{w-z}=N^{1/3}\int_0^{\infty}\exp{\{-N^{1/3}t(w-z)\}}dt}$.$\blacksquare$\\
We now indicate the idea of the proof, which is very similar to that in \cite{BaikGBAPeche}. We will perform a saddle point analysis of the kernels $H_N$ and $J_N$. The critical points
for $F$ satisfy  $ F'(w)=w+\dfrac{1}{w}-2 =0.$ Such an equation admits a single critical point,
$w_c=1=\pi_1$, and \be\label{deriveesen1}F''(w_c)=1-\frac{1}{w_c^2}=0,\: F^{(3)}(w_c)=2. \ee 
Intuitively, the leading terms of the asymptotic expansions of (\ref{HN}), (\ref{JN}) are
obtained by performing the corresponding integrals on a neighborhood of width $N^{-1/3}$ of $w_c$.
The steepest descent (resp. ascent ) curve for $F$ comes to $w_c$ with an
angle of $\pm \pi/3$ (resp $2\pi/3$) with respect to the real axis. Yet, as the integrand has a
pole at $w_c$, one needs to deform these path so
that $\Gamma$ encircles $w_c$ but does not cross $\gamma.$ 
Essentially, we will have to show that the ascent and descent contours, deformed in such a way,
still satisfy the saddle point analysis requirements.
We now define the expected limiting kernels. 
Let $\Gamma_{\infty}$ be the contour going from $\infty e^{-2i\pi/3}$ to $\infty e^{2i\pi/3}$, crossing the real axis on the right of the origin, oriented counterclockwise. Let $\gamma_{\infty}$ be the contour going from $\infty e^{-i\pi/3}$ to $\infty e^{i\pi/3}$, oriented from bottom to top and crossing the real axis on the right handside of $\Gamma_{\infty}.$ A plot of these contours is given on Figure \ref{fig:contourgammasinfinis}.

We then set \begin{eqnarray}&&\label{Hinfini}
H_{\infty}(x)=\frac{\exp{\{-\epsilon x\}}}{2\pi}\int_{\Gamma_{\infty}}\exp{\{xa-\frac{a^3}{3}\}}\frac{1}{a^k}da,\\
&& \label{Jinfini}J_{\infty}(y)=\frac{\exp{\{\epsilon y\}}}{2\pi}\int_{\gamma_{\infty}}\exp{\{-yb+\frac{b^3}{3}\}}b^kdb.
\end{eqnarray}
The end of this section is devoted to the proof of the following result.

\bp \label{Prop: HinfiniJinfini}
Fix $\epsilon >0$ and let $Z_N=g(w_c)\exp{\{N F(w_c)\}}N^{-k/3}.$

 For any fixed $y_o \in \mathbb{R},$ there exists $C>0,\:c>0,$ an integer $N_o>0$ such that
\be \label{Hx}
\Big|Z_N H_N(x)-H_{\infty}(x)\Big|\leq \frac{C \exp{\{-cx\}}}{N^{1/3}}, \text{ for any }x\geq y_o, \:N\geq N_o.
\ee
\be \label{Jy}
\Big|\frac{1}{Z_N} J_N(y)-J_{\infty}(y)\Big|\leq \frac{C \exp{\{-cy\}}}{N^{1/3}}, \text{ for any }y\geq y_o, \:N\geq N_o.
\ee

\ep
\brem The fact that Proposition \ref{Prop: HinfiniJinfini} implies Theorem \ref{Th: UnigalWig} is proved in \cite{BaikGBAPeche}, Section 3.3. It follows in particular from the fact that $J_{\infty}(y)=it^{(k+1)}(y)e^{\{\epsilon y\}}$ and $H_{\infty}(x)=is^{(k)}(x)e^{\{-\epsilon x\}},$ where $t^{(k+1)}$, $s^{(k)}$ are defined in (\ref{smBaik}).\erem
\brem Before beginning the proof of Proposition \ref{Prop: HinfiniJinfini}, it is convenient to note that the exponential term $F$, given in (\ref{defdeF}), satisfies $F(z)=\overline{F(\overline{z})}$. Thus, we may only consider the parts of the contours $\Gamma$ or $\gamma$ lying in the upper half plane $\{z\in \mathbb{C}, \:Im(z)>0\}$. Estimates for the remaining contours are obtained by conjugation when needed. This is valid for the whole paper.
\erem
\subsection{Estimate for $Z_NH_N$.}
This subsection is devoted to the proof of Formula (\ref{Hx}). We first define an ascent curve
$\Gamma$ for $F$. We then deduce the asymptotic expansion of $H_N$.

\subsubsection{Contour for the saddle point analysis}
In this part, we give an ascent curve for $F$ and also prove that the third order Taylor expansion of $F$ (as heuristically explained in the preamble) can be made in some disk around $w_c$.\\
Let $\Gamma$ be the contour defined in the following way.
\begin{eqnarray}
& \Gamma_o= w_c+\frac{\epsilon e^{i\theta}}{2 N^{1/3}},\quad \theta \in [0, 2\pi/3],
&\Gamma_1= w_c+te^{i2\pi/3},\quad \frac{\epsilon}{2 N^{1/3}}\leq t \leq 2,\cr
&\Gamma_2=\sqrt{3}i-t,\quad 0\leq t \leq R_o,
&\Gamma_3=i(\sqrt{3}-t)-R_o,\quad 0\leq t\leq \sqrt{3}.
\end{eqnarray}
Here $R_o$ is chosen large enough so that $\Gamma$ encircles all the eigenvalues $\pi_i, i=1,\ldots, r+1,$ and will be fixed later.
Finally define $\displaystyle{\Gamma =\cup_{i=0}^3 \Gamma_i \cup \overline{\cup_{i=0}^3 \Gamma_i},}$
oriented counterclockwise, as on Figure \ref{fig:contourGamma} below.
\begin{figure}[h]
\begin{center}
\begin{tabular}{c}
\epsfig{figure=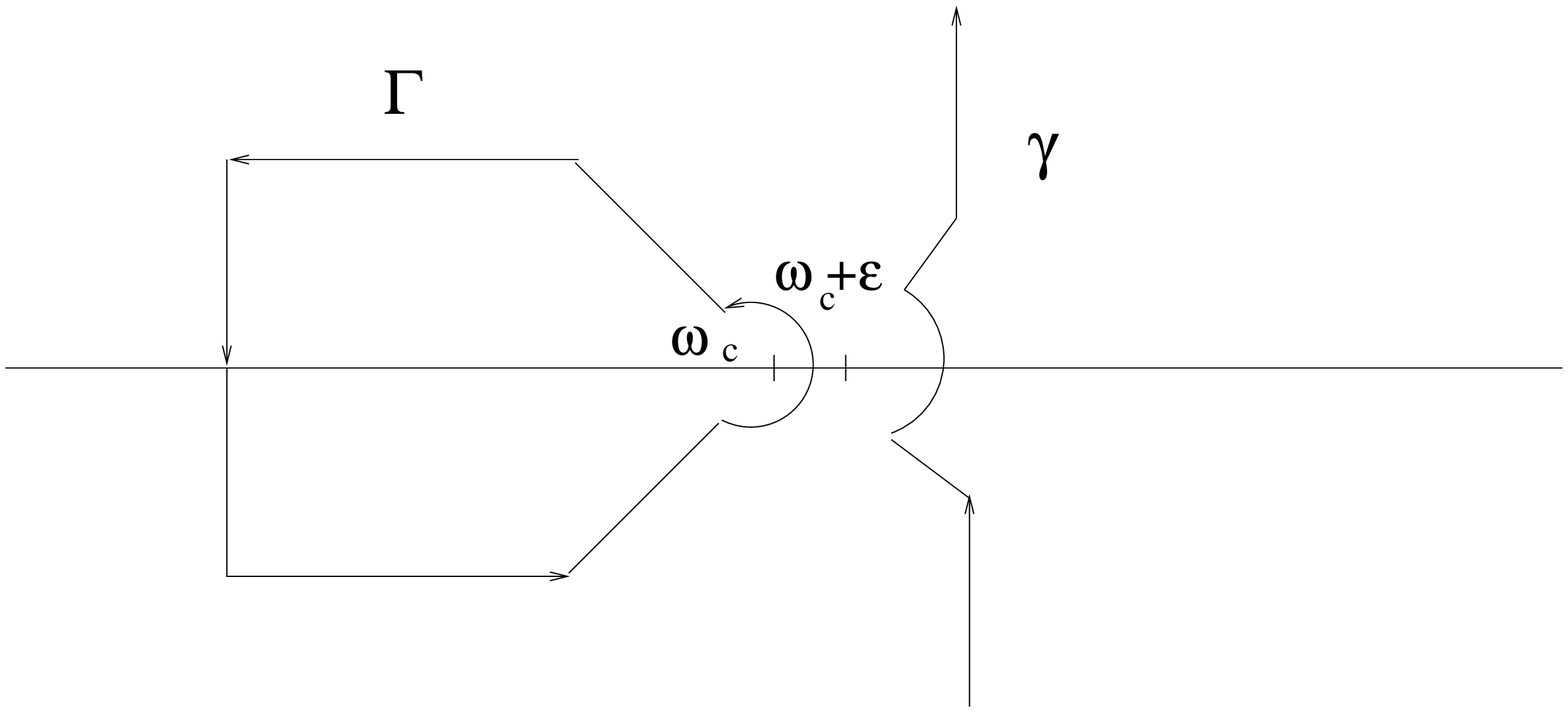,height=3.8cm}
\end{tabular}
\caption{Contours $\Gamma$ and $\gamma.$}
\label{fig:contourGamma}
\end{center}
\end{figure}
\bl
\label{Lem: steepascent}
$Re (F)$ increases as $z$ along $\Gamma_1 \cup\Gamma_2$ and if $z^*=\Gamma_1\cap \Gamma_2,$
$\displaystyle{\min_{z \in \Gamma_2\cup\Gamma_3} Re \left (F(z)\right )=Re \left (F(z^*)\right ).}$
\el
\paragraph{Proof of Lemma \ref{Lem: steepascent} :}
For $z\in \Gamma_1$, we have that
%$$Re \left (F(z)\right )=\frac{(1-t-t^2/4)}{2}-2(1-\frac{t}{2}-\tilde w_c)+\frac{1}{2} \log (1-t+t^2),$$
%so that 
$\displaystyle{\frac{d}{dt}Re\left ( F(w_c+t e^{2i\pi/3})\right )=\frac{1}{2}\frac{t^2(2-t)}{1-t+t^2}\geq 0,}$ for $t\leq 2.$ 
%Thus $Re (F)$ is a strictly increasing function on $[0,2]$ which achieves its maximum at $t=2$.
Then, along $\Gamma_2,$
$\displaystyle{\frac{d}{dt}Re \left (F(\sqrt{3}i-t)\right )=t+2+\frac{t}{|t+\sqrt{3}+t|^2}>2+t,}$
so that $Re (F)$ achieves its minimum on $\Gamma_2$ at $z^*.$
Finally we choose $R_o$ such that  
\be \label{defdeR_o}Re \left (F(-R_o+it)\right )=\frac{R_o^2}{2}-\frac{t^2}{2}+2R_o-\frac{1}{2}\log{|R_o+it|^2}>Re \left ( F(z^*)\right ),\: \forall \: |t|\leq \sqrt 3.\ee
This can be achieved if $R_o$ is chosen large enough.$\blacksquare$

\paragraph{}We now determine some disk around $w_c$ where the third order Taylor expansion of $F$ holds.
Let now $\delta$ be chosen so that
\be 0<\delta<1/2 \quad \text{and }\frac{\delta}{4(1-\delta)^4}\leq 1/6 .\label{delta} \ee

\bl \label{Lem: taylor2} In the disk $\{|z-w_c| \leq \delta\}$, 
$\displaystyle{\Big|F(z)-F(w_c)-\dfrac{F^{(3)}(w_c)}{3!} (z-w_c)^3\Big|\leq \frac{F^{(3)}(w_c)}{3!}\frac{|z-w_c|^3}{2}.}$
\el
\paragraph{Proof of the Lemma \ref{Lem: taylor2}: }
This follows from the Taylor expansion
%Such a choice for $\delta$ is explained by the following. For $z\in \Gamma'$, one has
\begin{equation}
\Big|F(z)-F(w_c)-\dfrac{F^{(3)}(w_c)}{3!} (z-w_c)^3\Big|\leq \max_{\Gamma'}\frac{|F^{(4)}(z)|}{4!}|z-w_c|^4
\leq \frac{\delta }{4!(1-\delta )^4}|z-w_c|^3 \leq \frac{|z-w_c|^3}{6}. \label{145}\blacksquare
\end{equation}

\brem \label{rem: deltazn}The above Lemmas imply in particular that  $\forall t\leq \delta$,
 $ \displaystyle{Re \left (F(w_c+te^{2i\pi/3})\right )\geq F(w_c)+t^3/6,}$
 and
 $\displaystyle{ \min_{\Gamma''}Re \left (F(z)\right )\geq F(w_c)+\delta^3/6.}$
 Here we have used that $\delta \leq 1/2$ and (\ref{deriveesen1}).
 \erem
The latter remark suggests that the leading term in the asymptotic expansion of $H_N$ is given by
the integral performed on the disk $\{|z-w_c| \leq \delta\}$. We thus split the contour $\Gamma=\Gamma'\cup\Gamma"$ where $\Gamma'=\Gamma \cap
\{|z-w_c| \leq \delta\}$ and $\Gamma''=\Gamma\setminus \Gamma'.$ Let also $\Gamma'_{\infty}$ be
the image of $\Gamma'$ under the map $z\mapsto N^{1/3}(z-w_c) $ and
$\Gamma''_{\infty}=\Gamma_{\infty}\setminus \Gamma'_{\infty}.$ 
We split accordingly the kernels $H_N$ and $H_{\infty}$, which we write $H_N(x)=H_N'(x)+H_N''(x)$ and $H_{\infty}(x)=H_{\infty}'(x)+H_{\infty}''(x),$ where
\begin{equation}H_N'(x)=\frac{N^{1/3}}{2\pi}
\int_{\Gamma'}\frac{e^{\{-NF(z)\}}}{g(z)(z-w_c)^k}e^{\{N^{1/3}x(z-\tilde w_c)\}}dz\text{ and }
H_{\infty}'(x)=\frac{e^{\{-\epsilon x\}}}{2\pi}\int_{\Gamma'_{\infty}}e^{\{xa-\frac{a^3}{3}\}}\frac{1}{a^k}da.\nonumber
\end{equation}
We now turn to the end of the proof of Formula (\ref{Hx}).
 \subsubsection{The case $x$ is bounded.}
Formula (\ref{Hx}) follows in this case from the following Lemma.
 \bl \label{Lem: TWsecond}
 Let $y_o>0$ be given. Then, there exists constants $C(y_o)>0$, $N_o>0$ such that,  for any $|x|\leq y_o,$ and $N\geq N_o$,
 $$|Z_N H_N(x)-H_{\infty}(x)| \leq \frac{C(y_o)}{ N^{1/3}}.$$
 \el
 \paragraph{Proof of Lemma \ref{Lem: TWsecond}: }We consider the contributions of $\Gamma'$ and $\Gamma''$ separately.
We first prove
 \begin{eqnarray} \label{148}
 &&|Z_N H''(x)|\leq \exp{\left\{-N \frac{\delta^3}{12}\right \}}, \quad
 |H_{\infty}''(x)|\leq\exp{\{-N\delta^3/6\}}. %\label{152}
 \end{eqnarray}
 Let then $L_{\Gamma''}$ be the length of $\Gamma''$, $C(R_o)=R_o+2$, and $C_g$ be a constant such that
  \be \label{C_g}\frac{1}{C_g}\leq \min_{\Gamma''}|g|\leq C_g,\ee
  which is well defined since, by Assumption \ref{Hyp1}, $\pi_i,i=2,\ldots, r+1,$ lie in a compact interval of $(-\infty,1).$
 Then, using Remark \ref{rem: deltazn}, we have that
\begin{eqnarray}&|Z_N H''(x)|&\leq \frac{|g(w_c)|}{2\pi N^{(k-1)/3}} \int_{\Gamma''}\Big|\exp{\{N F(w_c)-NF(z)\}}\exp{\{N^{1/3}x(z-\tilde w_c)\}}\Big | \frac{|dz|}{|g(z)(z-z_c)^k|}\cr
&&\label{majoZN}\\
&&\leq\frac{|g(w_c)|}{2 \pi N^{(k-1)/3}}L_{\Gamma''}\frac{C_g}{\delta^k}\exp{\left \{N^{1/3}y_o C(R_o)\right \}}\exp{\left\{-N \frac{\delta^3}{6}\right \}} \leq \exp{\left\{-N \frac{\delta^3}{12}\right \}},\nonumber \end{eqnarray}
 for $N$ large enough. This yields the first part of (\ref{148}).
 The second inequality is straightforward from \cite{BaikGBAPeche}, formula (152), for instance.\\ %to obtain that
%$$\displaystyle{|H_{\infty}''(x)|\leq \frac{e^{\epsilon y_o}}{\pi} \int_{\delta N^{1/3}}^{\infty}\frac{1}{t^k}\exp{\{y_ot-\frac{t^3}{3}\}} dt\leq\exp{\{-N\frac{\delta^3}{6}\}} ,}$$
% for $N$ large enough. This gives (\ref{152}).\\

 We then turn to the contour $\Gamma'=\Gamma'_o\cup\Gamma'_1,$ where $\Gamma'_o:=\Gamma_o\cap \{|w-w_c|\leq \delta\}=\Gamma_o$ and
 $\Gamma'_1= \Gamma'\setminus\Gamma'_o.$  Here we assume that $\epsilon$ is chosen so that 
  $\epsilon\leq \delta$ and thus $\Gamma'_o=\Gamma_o$.
 Here we prove that
 \be|Z_N H_N'(x)-H_{\infty}'(x)|\leq \frac{C}{N^{1/3}}. \label{164}\ee
 One has
 \be \label{153}
 |Z_N H'_N(x)-H'_{\infty}(x)|\leq \frac{N^{1/3}}{2\pi}\int_{\Gamma'}\frac{e^{N^{1/3}xRe(z-\tilde w_c)}}{(N^{1/3}|z-w_c|)^k} \Big|e^{N(-F(z)+F(w_c))}\frac{g(w_c)}{g(z)}-e^{-N(z-w_c)^3/3}\Big||dz|.
 \ee
 We now skip the details (given in \cite{BaikGBAPeche}, page 26).
 %We first need to show that the fact $\Gamma$ does not pass exactly through $w_c$ does not affect the behavior of the leading term in the asymptotic expansion.
 Then for $z\in \Gamma'_o$, using (\ref{145}),
\begin{eqnarray}
&& \Big|\exp{\{N(F(w_c)-F(z))\}}-\exp{\{-N\frac{(w-w_c)^3}{3}\}}\Big|\cr &&\leq \max \left
(\Big|e^{N(F(w_c)-F(z))}\Big|, \Big|e^{-N\frac{(z-w_c)^3}{3}}\Big|\right
)N\Big|F(z)-F(w_c)-\frac{(z-w_c)^3}{3}\Big| \cr &&\leq  N C_o|z-w_c|^4\exp{\{N Re \left (
\frac{(z-w_c)^3}{16}\right )\}}\leq \frac{C_o\epsilon ^4}{16 N^{1/3}}
\exp{\{\frac{\epsilon^3}{16}\} },\nonumber
\end{eqnarray}
 where $C_o=1/(1-\delta)^4$ is well defined since $\delta<1/2.$
 Similarly
 $\displaystyle{\big|\frac{g(w_c)}{g(z)}-1\big |\leq \frac{C_g C'_g}{2N^{1/3}}},$
 where $C_g$ is given by (\ref{C_g}) and
 $C_g'=\max\{|g'(s)|, s\in \Gamma'_o\cup \Gamma'_1\},$ which is well defined since the $\pi_i, i\geq 2$ in a compact set of $(-\infty,1)$. Thus, $\forall z\in \Gamma'_o$
\begin{equation}
 \Big|\exp{\{N(F(w_c)-F(z))\}}\frac{g(w_c)}{g(z)}-\exp{\{-N\frac{(w-w_c)^3}{3}\}}\Big|\leq \left (\exp{\{\frac{\epsilon^3}{16}\} }\frac{C_o\epsilon^4  }{16}+C_gC_g'/2\right )\frac{1}{N^{1/3}}.\label{158}
\end{equation}
  Using now that the length of $\Gamma'_o$ is $2\pi \epsilon N^{-1/3}/3,$ we obtain from (\ref{158})
  that there exists $C_1>0$ such that
   \be \frac{N^{1/3}}{2\pi}\int_{\Gamma'_o}\frac{e^{N^{1/3}xRe(z-\tilde w_c)}}{(N^{1/3}|z-w_c|)^k}
   \Big|e^{N(-F(z)+F(w_c))}\frac{g(w_c)}{g(z)}-e^{-N(z-w_c)^3/3}\Big||dz| \leq \frac{C_1}{N^{1/3}}. \label{estG'o}
  \ee
\paragraph{}Similarly for $z\in \Gamma'_1,$ one has that $\displaystyle{\exp{\{N^{1/3}xRe(z-\tilde w_c)\}}\leq \exp{\{N^{1/3}y_ot/2+\epsilon y_o\}}}$ and
  \be\label{166}\Big|e^{N(F(w_c)-F(z))}\frac{g(w_c)}{g(z)}-e^{-N(z-w_c)^3/3}\Big|\leq  (C_o+C_gC_g')(Nt^4+t)\exp{\{-N\frac{t^3}{6}\}}.\ee
 Now (see \cite{BaikGBAPeche}), we obtain from (\ref{166}) that there exists some $C_2>0$ such
 that, for $N$ large enough,
 \begin{eqnarray}
 &&\frac{N^{1/3}}{2\pi}\int_{\Gamma'_1}\frac{e^{N^{1/3}xRe(z-\tilde w_c)}}{(N^{1/3}|z-w_c|)^k}
   \Big|e^{N(-F(z)+F(w_c))}\frac{g(w_c)}{g(z)}-e^{-N(z-w_c)^3/3}\Big||dz|\cr
 &&\leq\frac{N^{1/3}}{\pi}(C_o^3+C_gC'_g)\int_{\frac{\epsilon}{2N^{1/3}}}^\delta
\!\left (\frac{1}{N^{1/3}t}\right )^k (Nt^4+t)\exp{\left \{\epsilon
y_o+\frac{y_otN^{1/3}}{2}-\frac{Nt^3}{6}\right \}} \leq \frac{C_2}{N^{1/3}}.\label{estG'1}
\end{eqnarray}
  Finally, combining (\ref{153}), (\ref{estG'o}), and (\ref{estG'1}), we obtain (\ref{164}).
Using now (\ref{148}), (\ref{164}), we then obtain that
  $\displaystyle{|Z_N H_N(x)-H_{\infty}(x)| \leq \frac{C(y_o)}{ N^{1/3}},}$
   for $N$ large enough. $\blacksquare$

   \subsubsection{The case $x$ positive}
Fromula (\ref{Hx}) follows in this case from the following Lemma.
   \bl \label{LemTWJN}Assume $x>0$, then there exist $C>0$, $N_o>0$ such that for $N\geq N_o$, $$|Z_NH_N(x)-H_{\infty}(x) |\leq \frac{C \exp{\{-\epsilon x/2}\}}{N^{1/3}}.$$
   \el

\paragraph{Proof of Lemma \ref{LemTWJN}:} The thing that makes it all here is that the whole contour $\Gamma $ lies on the half plane $Re (z-\tilde w_c)<0,$ where $\tilde w_c$ has been defined in (\ref{lesrescalingsTh1.1}). This gives that, for large positive $x$, the kernel $Z_N H_N$ decays exponentially, as we now explain.\\
 For $z \in \Gamma''$, one has that
 $\displaystyle{Re (z-\tilde w_c)\leq -\frac{\epsilon}{N^{1/3}}-\frac{\delta}{2},}$
 yielding from (\ref{majoZN}) that
 $$|Z_NH''(x)|\leq \exp{\left\{-\epsilon x-\frac{\delta N^{1/3}x}{2}-N\frac{\delta^3}{12}\right\}}\text{ for $N$ large enough}.
$$
It is also easy to check that $\displaystyle{|H''_{\infty}(x)|\leq \exp{\{-\epsilon x-\delta \frac{N^{1/3}x}{2}-\frac{N\delta^3}{6}\}},}$
   for $N$ large enough.\\
We now consider the part of $Z_NH_N$ (resp. $H'_{\infty}$) corresponding to the integral performed over $\Gamma'$ (resp. $\Gamma'_{\infty}$), along which one has that 
 $\:\displaystyle{\exp{\{N^{1/3}xRe(z-\tilde w_c)\}}\leq \exp{\{-\epsilon x/2\}}.}\:$
 Inserting the latter in (\ref{153}), and performing the same computations as for the case where $x$ lies in a compact set, we obtain that
 $$|Z_NH_N'(x)-H'_{\infty}(x) |\leq \frac{C_2 \exp{\{-\epsilon x/2}\}}{N^{1/3}}.\blacksquare$$
\brem \label{rem: crucialpoints}There are two crucial steps in the above proof. The first one is the definition of an ascent curve $\Gamma$, which coincides with the steepest ascent curve for $F$ in an annulus $\{\epsilon \leq |z-w_c|\leq \delta\}.$ The second step is to determine a $\delta >0$ such that Lemma \ref{Lem: taylor2} holds. This second step also ensures that we can find $\epsilon$ small enough so that $\Gamma$ encircles $w_c$ but remains on the left handside of $\tilde w_c$.  Once these two points obtained, one only needs a good enough control of the perturbative term $g$ along $\Gamma$, so that the end of the proof follows. This remark will be the basis for the proof of Theorem \ref{theo: rhoN<<3/7}.
\erem
\subsection{Estimate for $\frac{1}{Z_N}J_N(y)$ }
This subsection is devoted to the proof of Formula (\ref{Jy}). We first define a descent curve
$\gamma$ for $F$. Then, we obtain the asymptotic expansion of $J_N$.
%, examining separately the cases
%where $y$ is bounded or positive.
\subsubsection{Contour for the saddle point analysis }
We now give a descent curve for $F$. Define
\begin{equation}
\gamma_0= w_c+\frac{3\epsilon e^{i\theta}}{N^{\frac{1}{3}}} ,\:\: 0\leq \theta \leq \frac{\pi}{3}; \:\:
\gamma_1=w_c+te^{i\frac{\pi}{3}}, \:\: \frac{3\epsilon }{N^{\frac{1}{3}}}\leq t\leq t_o; \:\: \gamma_2=w_c+t_oe^{i\frac{\pi}{3}}+it, \:\:  t\geq 0.\label{gamma2}
\end{equation}
Actually, we choose $t_o>\delta,$ where $\delta $ is given by (\ref{delta}).
Finally let $\gamma$ be the contour $\gamma=\cup_{i=0}^2 \gamma_i \cup \overline{\cup_{i=0}^2 \gamma_i}$ oriented from bottom to top, as on Figure \ref{fig:contourGamma}.

%\begin{figure}[htbp]
% \begin{center}
%\begin{tabular}{c}
% \epsfig{figure=gammabord.eps,height=4cm, angle=0}
%\end{tabular}
% \caption{Contour $\gamma$.
%\label{fig:contourgamma}}
%\end{center}
%\end{figure}
\bl \label{steepdescent}
$Re (F)$ is decreasing on $\gamma_1\cup\gamma_2$ as $Im (w)$ increases. And $\exists \:C_o>0$ such that, if $w^*=w_c+t_oe^{i\pi/3}$,
$\displaystyle{Re \left (F(w^*+it)\right )\leq Re \left (F(w^*)\right )-\frac{C_ot^2}{2}, \:t\geq 0 .}$
\el
\paragraph{Proof of Lemma \ref{steepdescent}:}
One can check that 
$\dfrac{d}{dt}Re \left (F(w_c+t e^{i\pi/3})\right)=\dfrac{-t^2(2+t)}{2(1+t+t^2)}<0, \: \forall t>0.$
%Instead of having a finite length contour as for $\Gamma$, we have to bound the real part of points on $\gamma$. Otherwise the kernel could be hardly be majorized in the case $y<0$. This is the reason why we have imposed the contour $\gamma_1$ to stop at $t_o$.
Then, along $\gamma_2$, and for $C_o=C_o=1-1/|w^*|^2 >0,$ one has
$\displaystyle{\frac{d}{dt}Re \left (F(w^*+it)\right )
%=-Im(w^*+t)(1-\frac{1}{|w|^2})
\leq -C_o(t+\frac{\sqrt{3}}{2}t_o)}.\blacksquare$\\

Let now $\delta $ be given as in (\ref{delta}), so that (\ref{145}) still holds. We split as before the contour $\gamma$. Set $\gamma'=\gamma\cap\{|w-w_c|\leq \delta\}$ and  $\gamma''=\gamma\setminus \gamma'.$
Let also $\gamma'_{\infty}$ be the image of $\gamma'$ under the map $w\mapsto N^{1/3}(w-w_c),$ and $\gamma''_{\infty}=\gamma_{\infty}\setminus \gamma'_{\infty}.$ Define then $J_1''=J_N-J'_N-J''_2, \text{ and } J''_{\infty}=J_{\infty}-J'_{\infty}$ where
\begin{eqnarray}&&J_N'(y)=\frac{N^{1/3}}{2\pi}\int_{\gamma'}g(w)(w-w_c)^k\exp{\{NF(w)\}}\exp{\{-N^{1/3}y(w-\tilde
w_c)\}}dw,\cr 
&&J''_2(y)=\frac{N^{1/3}}{2\pi}\int_{\gamma_2}g(w)(w-w_c)^k\exp{\{NF(w)\}}
\exp{\{-N^{1/3}y(w-\tilde
w_c)\}}dw,\cr
&&J'_{\infty}(y)=\frac{\exp{\{\epsilon y\}}}{2\pi}\int_{\gamma'_{\infty}}\exp{\{-yb+\frac{b^3}{3}\}}b^kdb
.\nonumber
\end{eqnarray}
We now prove Formula (\ref{Jy}).
\subsubsection{The case $y$ in a compact interval}
Formula (\ref{Jy}) readily follows in this case from the following Lemma.
\bl \label{TWJNVraiment} Let $y_o>0$ be given. Then, there exists $C=C(y_o)>0$, $N_o$ such that for any $|y|\leq y_o$,
$$|\frac{1}{Z_N}J_N(y)-J_{\infty}(y)| \leq \frac{C}{N^{1/3}},\forall \: N\geq N_o.$$ 
\el

\paragraph{Proof of Lemma \ref{TWJNVraiment}: }Let us first consider the kernel $J_N''=J_1''+J_2''$.
We now show that there exists $C>0$ such that \be|\frac{1}{Z_N}J_N''(y)|\leq C\exp{\{-N\frac{\delta^3}{12}\}}\label{estp19}
\ee
for $N$ large enough.
The only difference from the preceding subsection is that $\gamma''$ is not of finite length. We first consider the integral performed on $\gamma_2$. \be |\frac{1}{Z_N}J''_2(y)|\leq
\frac{(N^{1/3})^{k+1}}{2\pi |g(w_c)|}\int_{\gamma_2}e^{N^{1/3}y_o Re(w-\tilde w_c)}e^{N\left (Re\left (
F(w)-F(w_c)\right)\right )}|w-w_c|^k||g(w)| |dw|.\label{encoreunexpli1}\ee Now, using Lemma \ref{steepdescent},
and the fact that $Re \left (F(w^*)\right) \leq Re \left (F(w_c)\right)-\delta^3/6$ (which follows from the fact that $t_o\geq
\delta$ and Remark \ref{rem: deltazn}), one has

$$
|\frac{1}{Z_N}J''_2(y)|\leq \frac{(N^{1/3})^{k+1}}{2\pi
|g(w_c)|}e^{N^{1/3}\frac{y_ot_o}{2}-N\frac{\delta^3}{6}}\int_{i\mathbb{R}_+} e^{-N
\frac{C_ot^2}{2}}\prod_{i=2}^{r+1}\sqrt{\frac{(1+t_o/2+A)^2+(\sqrt3t_o/2+t)^2}{(1+t_o/2)^2+(t_o\sqrt
3 /2+t)^2 }}dt,
$$
where $A$ is chosen such that $|\pi_i|<A, i=2,\ldots, r+1.$ Now, under Assumption \ref{Hyp1}, $A$ can be chosen independently of $N$. Thus, for $N$ large enough,
\be\label{estimeejnexplique}|\frac{1}{Z_N}J_2''(y)|\leq \exp{\{-N\frac{\delta^3}{12}\}}.\ee
The remaining contour $\gamma''_1=\gamma''\setminus \gamma_2$ has a finite length
$L_{\gamma''_1}$, independent of $N$. Define also $\tilde C_g =\max_{w \in \gamma''_1} , |g(w)|$,
which is uniformly bounded. Now, using that, for $N$ large enough, $ Re \left (F(w)\right )\leq
F(w_c)-\delta^3/6 \text{ and }Re (w-\tilde w_c)\leq t_o, \: \forall w \in \gamma''_1$, we obtain
that
\be\label{estimeejnexplique2}|\frac{1}{Z_N}J_1''(y)|\leq \frac{(N^{1/3})^{k+1}}{2\pi |g(w_c)|}\delta^k \tilde C_g L_{\gamma''_1}\exp{\left\{N^{1/3}y_ot_o-N\frac{\delta^3}{6}\right\}}\leq \exp{\{-\frac{N\delta^3}{12}\}}, \ee for $N$ large enough. Combining  (\ref{estimeejnexplique}) and (\ref{estimeejnexplique2}) yields (\ref{estp19}).
\\And inserting $b=te^{i\pi/3}$, with $ t\geq \delta N^{1/3}$, in (\ref{Jinfini}) yields that (see formula (188) in \cite{BaikGBAPeche})\\
$\displaystyle{|J''_{\infty}(y)|\leq \exp{\{-\frac{N\delta^3}{6}\}}, \text{  for $N$ large enough.}}$ 
Finally, mimicking the proof of (\ref{166}), (\ref{158}), and using the same arguments as for $H_N$ (see Remark \ref{rem: crucialpoints}), it is easy to show that $\exists C_3(y_o)>0,$ such that
$\displaystyle{|\frac{1}{Z_N}J_N'(y)-J'_{\infty}(y)| \leq \frac{C_3(y_o)}{N^{1/3}}.}\blacksquare$
\subsubsection{The case $y>0$}
Formula (\ref{Jy}) follows in this case from the following Lemma.
\bl \label{Lem ypposirajout}There exist $C>0$, $N_o>0$ such that, $\forall y>0$, and $\forall N\geq N_o,$
$$|\frac{1}{Z_N}J_N(y)-J_{\infty}(y)| \leq \frac{C}{N^{1/3}}\exp{\{-\frac{\epsilon y}{2}\}}.$$
\el
\paragraph{Proof of Lemma \ref{Lem ypposirajout}: }We only give the mains ideas. One has
\begin{eqnarray}&\forall w \in \gamma_2,\:\:
Re(w-\tilde w_c)=Re(w^*-\tilde w_c)=\dfrac{t_o}{2}-\dfrac{\epsilon}{N^{1/3}};
&\forall w\in \gamma''_1, \:\:
Re(w-\tilde w_c)\geq \dfrac{\delta}{2}-\dfrac{\epsilon}{ N^{1/3}};\cr
&\forall w\in\gamma'_1,\:\: Re(w-\tilde w_c)=\dfrac{ |w-w_c|}{2}-\dfrac{\epsilon}{N^{1/3}}\geq \dfrac{\epsilon}{2 N^{1/3}};\:\:
&\forall w \in \gamma_o,\:\: Re(w-\tilde w_c)\geq \dfrac{\epsilon}{2 N^{1/3}}.\label{43}
\end{eqnarray}
Note that (\ref{43}) explains why we choose a circle of ray $3\epsilon $ for $\gamma_o$.
Here we assume that $\epsilon$ is small enough so that $\epsilon-\delta/2<-\epsilon/2.$ This gives that the whole contour $\gamma$ lies on the right of $\tilde w_c.$
We then insert the above estimates in e.g. (\ref{encoreunexpli1}) and copy the proof of (\ref{estimeejnexplique}). The same is done for the integral performed on $\gamma_1''$. Then, one readily obtains that, for $N$ large enough,
$$|\frac{1}{Z_N}J_N''(y)|\leq \exp{\left \{\epsilon y-\frac{\delta}{2} N^{1/3}y-N\frac{\delta^3}{12}\right \}},\text{ 
while }|J_{\infty}''(y)|\leq \exp{\left \{\epsilon y-\frac{\delta N^{1/3}y}{2}-N\frac{\delta^3}{6}\right \}}.$$
Finally, using (\ref{43}) and mimicking the estimates of the preceding subsection, we obtain that
$$|\frac{1}{Z_N}J'_N(y)-J'_{\infty}(y)| \leq \frac{C}{N^{1/3}}\exp{\{-\frac{\epsilon y}{2}\}}.\blacksquare$$
%This finishes the proof of Lemma \ref{Lem ypposirajout} and gives (\ref{Jy}) in this case. $\blacksquare$
\section{Proof of Theorem \ref{Th: GUEWig}}
In this section, we assume that $\pi_1$ lies in a compact interval of $(1,\infty)$ and that Assumption \ref{Hyp1} holds.
Let now $\epsilon>0 $ be fixed and set
\be \tilde \pi_1=\pi_1+\frac{\epsilon}{\sqrt{N}}, \:
u=C(\pi_1)+\frac{x}{\sigma^2\sqrt N}, \:v=C(\pi_1)+\frac{y}{\sigma^2\sqrt N},\label{rescGGUE} \ee and let 
$\displaystyle{
K_N'(x,y)=\frac{1}{\sigma^2 \sqrt N} K_N\left (C(\pi_1)+\frac{x}{\sigma^2\sqrt N},C(\pi_1)+\frac{y}{\sigma^2\sqrt N}\right )\exp{\left \{\frac{(y-x)}{\sigma^2}\tilde \pi_1\right \}}}
$ be the associated rescaled correlation kernel.
Define $F_{C(\pi_1)}(w)=w^2/2-C(\pi_1)(w-\tilde \pi_1)+\log w,$ where we use the principal branch of the logarithm
($e^{\pm N \log w}=w^{\pm N}$). We now bring $K_N'(x,y)$ to the form (\ref{formeproduit}).
Let $\Gamma$ and $\gamma$ be contours as in Proposition \ref{Prop: corrJoh}.
\bp %Let $K_N'$ be defined by (\ref{refKN'GUEfixe}), then
 $\displaystyle{K_N'(x,y)=-\int_{0}^{\infty}H_N(\frac{x+t}{\sigma^2})J_N(\frac{y+t}{\sigma^2})dt,}$
with
\begin{eqnarray} &&H_N(\frac{x}{\sigma^2})=\frac{\sqrt N }{2\pi\sigma^2}\int_{\Gamma}\frac{1}{(z-\pi_1)^k g(z)}\exp{\left\{\sqrt N \frac{x}{\sigma^2}(z-\tilde \pi_1)\right \}}\exp{\left\{-NF_{C(\pi_1)}(z)\right\}}dz ,\label{HNGGUE}\cr
&&J_N(\frac{y}{\sigma^2})=\frac{\sqrt N}{2\pi\sigma^2}\int_{\gamma}(w-\pi_1)^k g(w)\exp{\left\{-\sqrt N \frac{y}{\sigma^2}(w-\tilde \pi_1)\right\}}\exp{\left\{NF_{C(\pi_1)}(w)\right\}}dw .\label{JNGGUE}
\end{eqnarray}
\ep

\paragraph{} We briefly indicate the idea of the proof of Theorem \ref{Th: GUEWig}. Here, the critical points to be considered
satisfy
$F_{C(\pi_1)}'(w)=w+1/w-C(\pi_1)=0.$
They are given by $\pi_1$ and $1/\pi_1$ and are non degenerate. 
%Thus a second order Taylor expansion of the exponential term should give the asymptotics of the kernels $H_N$ and $J_N$.
One can check that
\be\label{deriveesecondeenpi}F''(\pi_1)=1-\frac{1}{\pi_1^2}=\frac{1}{\sigma^2(\pi_1)}>0.\ee
 We will show that, as $\Gamma$ has to encircle the pole $\pi_1$, the contribution of the sole pole $\pi_1$ will give the leading term in the asymptotic expansion.
In the subsequent, we note $\sigma^2=\sigma^2(\pi_1)$ and define now
 the expected limiting kernels.\\ Let $\gamma_{\infty}=2\epsilon +i\mathbb{R}$ (resp. $\Gamma_{\infty}=\epsilon e^{i\theta}, 0\leq \theta \leq 2\pi,$) oriented from bottom to top (resp. counterclockwise). Set then
\begin{eqnarray}
&&H_{\infty}(\frac{x}{\sigma^2})=\frac{1}{2\pi\sigma^2}\exp{\{-\epsilon \frac{x}{\sigma^2}\}}\int_{\Gamma_{\infty}}\frac{1}{a^k}\exp{\left\{-\frac{a^2}{2\sigma^2}+\frac{x}{\sigma^2}a\right \}}da, 
\cr &&
J_{\infty}(\frac{y}{\sigma^2})=\frac{1}{2\pi \sigma^2}\exp{\{\epsilon \frac{y}{\sigma^2}\}}\int_{\gamma_{\infty}}s^k \exp{\left\{\frac{s^2}{2\sigma^2}-\frac{y}{\sigma^2}s\right \}}ds.\label{JinfiniGUE}
\end{eqnarray}
The aim of the rest of this section is to prove the following result.
\bp\label{Prop: GUE}
Assume $\epsilon>0$ is fixed, and set $Z_N =g(\pi_1) N^{-k/2}\exp{\left\{N F_{C(\pi_1)}(\pi_1)\right\}}.$
For any fixed $y_o\in \mathbb{R}$, there exists constants $C>0$, $c>0$, and an integer  $N_o>0$ such that
\begin{eqnarray} &&\label{JNGGUE2}\Big|\frac{1}{Z_N}J_N(\frac{y}{\sigma^2})-J_{\infty}(\frac{y}{\sigma^2})\Big|\leq \frac{C\exp{\{-c\frac{y}{\sigma^2}\}}}{\sqrt N},\text{ for any }y\geq y_o, \:N\geq N_o.\\
&& \label{HNGGUE2}\Big|Z_NH_N(\frac{x}{\sigma^2})-H_{\infty}(\frac{x}{\sigma^2})\Big|\leq \frac{C\exp{\{-c\frac{x}{\sigma^2}\}}}{\sqrt N}, \text{ for any }x\geq y_o,\: N\geq N_o.\end{eqnarray}
\ep
\brem The fact that Proposition \ref{Prop: GUE} implies Theorem \ref{Th: GUEWig} follows from the equality   \be \label{eganoyau}-\exp{\{\epsilon \frac{(x-y)}{\sigma^2}\}}\int_0^{\infty} H_{\infty}(\frac{x+u}{\sigma^2})J_{\infty}(\frac{y+u}{\sigma^2})du = K(x,y),\ee
where $K(x,y)$ is the correlation kernel of the $k\times k$ GUE with parameter $\sigma^2$.
% Indeed, a straightforward computation yields that
%\begin{eqnarray}
%&&-\exp{\{\epsilon \frac{(x-y)}{\sigma^2}\}}\int_0^{\infty}H_{\infty}(\frac{x+u}{\sigma^2})J_{\infty}(\frac{y+u}{\sigma^2})du\cr
%&&=\int_{\Gamma_{\infty}}\int_{\gamma_{\infty}}\frac{1}{(2i\pi)^2\sigma^2}
%\left (\frac{t}{s}\right )^k\frac{1}{t-s}\exp{\left\{-\frac{s^2-sx}{2\sigma^2}+\frac{t^2-ty}{2\sigma^2}\right \}} dsdt =K(x,y).\end{eqnarray}
Formula (\ref{eganoyau}) follows from (\ref{noyauJohansson}) and a simple change of variables.
Another proof of (\ref{eganoyau})
is given in \cite{BaikGBAPeche}, Section 4.\erem

The proof of Proposition \ref{Prop: GUE} will be obtained in the following subsections.
\subsection{Estimate for $\displaystyle{\frac{1}{Z_N}J_N(\frac{y}{\sigma^2})}$ }
This subsection is devoted to the proof of formula (\ref{JNGGUE2}). The details of the proof will be skipped since the scheme is exactly the same as in the preceding Section. The key points are the following Lemmas. In the first one, we give the descent curve for $F_{C(\pi_1)}$. In the second one, we determine a disk where the second order Taylor expansion of $F_{C(\pi_1)}$ holds.\\

Let $\gamma_1$ be the contour $\gamma_1=\pi_1+\dfrac{2\epsilon}{\sqrt N}+it, t\in \mathbb{R}_+,$ $\displaystyle{\gamma=\gamma_1\cup \overline{\gamma_1}}$ and set $\pi'_1=\pi_1+\dfrac{2\epsilon}{\sqrt N}.$
\bl \label{Lemma21}There exists $C_o>0$ such that $Re \left (F_{C(\pi_1)}(\pi'_1+it)\right )\leq F_{C(\pi_1)}(\pi'_1)-C_ot^2/2, \forall t\in \mathbb{R}. $
\el
\paragraph{Proof of Lemma \ref{Lemma21}: } $\displaystyle{\frac{d}{dt}Re \left (F_{C(\pi_1}(\pi'_1+it)\right )=-t\left(1-\frac{1}{|\pi'_1+it|^2}\right)
\leq -C_o t}$
 since $\pi'_1>\pi_1$ lies in a compact interval of $(1,\infty).\blacksquare$

 \paragraph{}
Let now $\delta $ be such that
 \be
 \label{delta2} \frac{\delta}{(\pi_1-\delta)^3}\leq\frac{1}{4\sigma^2(\pi_1)}\text{ and } \delta\leq \pi_1/2.
  \ee
  \bl \label{TAYlor2JN} In the disk $\{|w-\pi_1|\leq \delta\}$, one has
 \be \left|F_{C(\pi_1)}(w)-F_{C(\pi_1)}(\pi_1)-\frac{F_{C(\pi_1)}''(\pi_1)}{2}(w-\pi_1)^2\right|\leq \frac{F_{C(\pi_1)}''(\pi_1)}{4}|w-\pi_1|^2. \label{formule55}\ee
  \el
\paragraph{Proof of Lemma \ref{TAYlor2JN} :}It is proved as Lemma \ref{Lem: taylor2}.$\blacksquare$\\
%For $|w-\pi_1|\leq \delta$, one has
%$$\left|F_{C(\pi_1)}(w)- F_{C(\pi_1)}(\pi_1)-\frac{F_{C(\pi_1)}''(\pi_1)}{2}(w-\pi_1)^2\right |\leq \max_{|w-\pi_1|\leq \delta} \frac{|F_{C(\pi_1)}^{(3)}(w)||w-\pi_1|^3}{3!}$$
%and $ \left |Re \left (F_{C(\pi_1)}^{(3)}(w)\right )\right|=\dfrac{2}{|w^3|}\leq 16.$
%Using (\ref{deriveesecondeenpi}), this gives (\ref{formule55}). $\blacksquare$\\

As before, we now split the contour into two parts. Let $\gamma'=\gamma\cap \{|w-\pi_1|\leq \delta\}$ and $\gamma''=\gamma\setminus \gamma'.$
Let also $\gamma'_{\infty}$ be the image of $\gamma'$ under the map
$w\mapsto \sqrt N(w-\pi_1)$ and $\gamma''_{\infty}=\gamma_{\infty}\setminus \gamma'_{\infty}.$
Set now $$J_N(\frac{y}{\sigma^2})=J'_N(\frac{y}{\sigma^2})+J''_N(\frac{y}{\sigma^2}),\quad J_{\infty}(\frac{y}{\sigma^2})=J'_{\infty}(\frac{y}{\sigma^2})+J''_{\infty}(\frac{y}{\sigma^2}),$$ where
$\displaystyle{J'_N(\frac{y}{\sigma^2})=\frac{\sqrt N}{2\pi\sigma^2}\int_{\gamma'}(w-\pi_1)^k g(w)\exp{\left\{-\sqrt N \frac{y}{\sigma^2}(w-\tilde \pi_1)\right\}}\exp{\left\{NF_{C(\pi_1)}(w)\right\}}dw }$ and \\$\displaystyle{J'_{\infty}(\frac{y}{\sigma^2})=\frac{1}{2\pi \sigma^2}\exp{\{\epsilon \frac{y}{\sigma^2}\}}\int_{\gamma'_{\infty}}s^k \exp{\left\{\frac{s^2}{2\sigma^2}-\frac{y}{\sigma^2}s\right \}}ds.}$

We only give the main steps of the proof. Let $y_o>0$ be given and assume first that $y$ lies in the interval $[-y_o, y_o].$
Then we show that there exists $ C>0,$ such that for $N$ large enough,
\be \Big |\frac{1}{Z_N}J_N''(y)\Big |\leq C \exp{\{-N\frac{\delta^2}{24}\}}, \quad \Big |J_{\infty}''(y)\Big |\leq C \exp{\{-N\frac{\delta^2}{24}\}}, \quad |\frac{1}{Z_N}J_N'(\frac{y}{\sigma^2})-J_{\infty}'(\frac{y}{\sigma^2})|\leq \frac{C}{\sqrt N} \label{noybla}.\ee

Here we have to take care of the fact that $\gamma$ does not exactly go through the critical point $\pi_1$.
Consider first $\gamma''$ and let $w^*=\gamma\cap \{|w-\pi_1|=\delta\}.$ From Lemma \ref{Lemma21},
\be \label{resplik}Re \bigl (F(w^*+it)-F(w^*)\bigr )\leq -C_o t, \:\forall t>0.\ee
Furthermore, as $N$ goes to infinity, $w^* \rightarrow \pi_1+i\delta,$ so that, for $N$ large enough, by Lemma \ref{Lem: taylor2},
\be Re \left (F_{C(\pi_1)}(w^*)\right )-Re \left (F_{C(\pi_1)}(\pi_1)\right )\leq -\frac{\delta^2}{12}
\label{245}.\ee
Combining (\ref{resplik}), (\ref{245}), and Remark \ref{rem: crucialpoints}, we obtain the first inequality in (\ref{noybla}), for $N$ large enough.
The second inequality is straightforward.
Conversely, Lemma \ref{TAYlor2JN} and Remark \ref{rem: crucialpoints} give the last inequality in (\ref{noybla}), since the perturbative term $|g(w) (w-w_c)^k|$ is uniformly bounded along $\gamma'$. This yields (\ref{JNGGUE2}) in this case.
Finally, we use the fact that $Re (w-\tilde \pi_1)>C, \:\forall w\in \gamma,$ for some constant $C>0$, and the same arguments as in the preceding Section, to obtain (\ref{JNGGUE2}) in the case $y>0$. %The detail is left.

\subsection{Estimate for $Z_N H_N(\frac{x}{\sigma^2})$}
This subsection is devoted to the proof of formula (\ref{HNGGUE2}).
We examine $Z_N H_N(\frac{x}{\sigma^2})$ as a residue integral and show that the sole residue at $z=\pi_1$ gives the leading term in the asymptotic expansion.
We thus split the contour accordingly. Let $\Gamma''$ be a contour that encloses $0$ and $\pi_i,i=2,\ldots, r+1$ but not $\pi_1,$ oriented counterclockwise. Then we readily obtain the following Proposition.
\bp
$Z_NH_N(\frac{x}{\sigma^2})=H_1(\frac{x}{\sigma^2})+H_2(\frac{x}{\sigma^2})$
where
\begin{eqnarray}
&H_{2}(\frac{x}{\sigma^2})=&\frac{g(\pi_1)e^{-\epsilon \frac{x}{\sigma^2}}}{2\pi \sigma^2(\sqrt N)^{k-1}}
\int_{\Gamma''}\frac{\exp{\{\sqrt N \frac{x}{\sigma^2}(z-\pi_1)\}}}{(z-\pi_1)^k g(z)}\exp{\left \{-N(F_{C(\pi_1)}(z)-F_{C(\pi_1)}(\pi_1))\right \}}dz,\cr
&H_{1}(\frac{x}{\sigma^2})=&\frac{e^{-\epsilon \frac{x}{\sigma^2}}}{\sigma^2} \int_{\Gamma_{\infty}} \frac{\exp{\{\frac{x}{\sigma^2}a\}}}{a^k}\frac{g(\pi_1)}{g(\pi_1+\frac{a}{\sqrt N})}\exp{\left\{-N\left (F_{C(\pi_1)}(\pi_1+\frac{a}{\sqrt N})-F_{C(\pi_1)}(\pi_1)\right )\right \}}da.\cr
&& \label{noyauH_1GUE}
\end{eqnarray}
\ep

\paragraph{}The proof of Formula (\ref{HNGGUE2}) is now divided into two facts, in which we examine separately the two kernels $H_1$ and $H_2.$
First we show that $H_1(\frac{x}{\sigma^2})$ behaves as $H_{\infty}(\frac{x}{\sigma^2}).$
\bfa \label{Lem: H_1} Given any fixed $y_o \in \mathbb{R},$ there exists constants $C>0,\: c>0,\: N_o>0$ such that
$$ |H_1(\frac{x}{\sigma^2})-H_{\infty}(\frac{x}{\sigma^2})|\leq \frac{C\exp{\{-c\frac{x}{\sigma^2}}\}}{\sqrt N},\text{ for any }x\geq y_o, \:N\geq N_o.$$
\efa
\paragraph{Proof of Fact \ref{Lem: H_1}: }
We only explain the main changes from \cite{BaikGBAPeche}, since the proof follows the same steps.
For any $l$ , the derivatives
$F_{C(\pi_1)}^{(l)}(\pi_1),$  $g^{(l)}(\pi_1)$ are all bounded, and $|g(\pi_1)|>0$ thanks to Assumption \ref{Hyp1}.
Thus, by a straightforward Taylor expansion, we have that
\begin{eqnarray}
&&\int_{\Gamma_{\infty}}\frac{1}{a^k}\frac{g(\pi_1)}{g(\pi_1+\frac{a}{\sqrt N})}\exp{\left\{-N\left(F_{C(\pi_1)}(\pi_1+\frac{a}{\sqrt N})-F_{C(\pi_1)}(\pi_1)\right )\right \}}\exp{\left\{\frac{x}{\sigma^2}a\right \}}da\cr
&&= \int_{\Gamma_{\infty}}\frac{1}{a^k}\exp{\left\{\frac{x}{\sigma^2}a-\frac{1}{2\sigma^2}a^2\right \}}\left (1+\sum_{j=1}^{k-1}\frac{q_j(a)}{(\sqrt N)^j}\right)da ,\label{62qj}
\end{eqnarray} for some polynomials $q_j, j=1, \ldots, k-1$ independent of $N$. Now (\ref{62qj}) and (\ref{noyauH_1GUE}) give Fact \ref{Lem: H_1}. $\blacksquare$

\paragraph{}
We now turn to the asymptotics of the kernel $H_2$.
\bfa \label{claim H2}For any fixed $y_o \in \mathbb{R},$ there exists  $C>0,\: c>0,\: N_o>0$ such that
$$\left|H_{2}(\frac{x}{\sigma^2})\right |\leq c\exp{\{-\epsilon x -CN\}},\text{ for any }N\geq N_o,\: x\geq y_o.$$
\efa

\paragraph{Proof of Fact \ref{claim H2}:} The proof is obtained by a saddle point analysis of the kernel $H_2$. We define the suitable contour $\Gamma'',$ that depends on some constants $\eta,\: R, \theta_o, x_o^*$ that will be fixed later.
Set $\pi^*=\max (1,\pi_2)$ and define 
\begin{eqnarray}
&\Gamma''_1=\dfrac{\pi_1+\pi^*}{2}+iy,\quad y \leq \eta,
&\Gamma''_2=\dfrac{\pi_1+\pi^*}{2}+i\eta -x, \quad 0\leq x\leq \frac{\pi_1+\pi^*}{2}-x_o^*,\cr %\label{defx*}\\
&\Gamma''_3=\dfrac{C(\pi_1)}{2}e^{i\theta},\quad \theta_o\leq \theta \leq \frac{\pi}{2},
&\Gamma''_4=i\dfrac{C(\pi_1)}{2}-x, \quad 0\leq x\leq R\cr
&\Gamma''_5=i(\dfrac{C(\pi_1)}{2}-t)-R,\quad 0\leq t\leq  \frac{C(\pi_1)}{2}.& \text{    }\nonumber
\end{eqnarray}
Set $\displaystyle{\Gamma''=\cup_{i=1}^5\Gamma''_i \cup \overline{\cup_{i=1}^5\Gamma''_i}.}$
A plot of the contours $\Gamma=\Gamma_{\infty}\cup\Gamma''$ and $\gamma$ is given on Figure \ref{fig:contoursGUE} below.
\begin{figure}[htbp]
 \begin{center}
 \begin{tabular}{c}
 \epsfig{figure=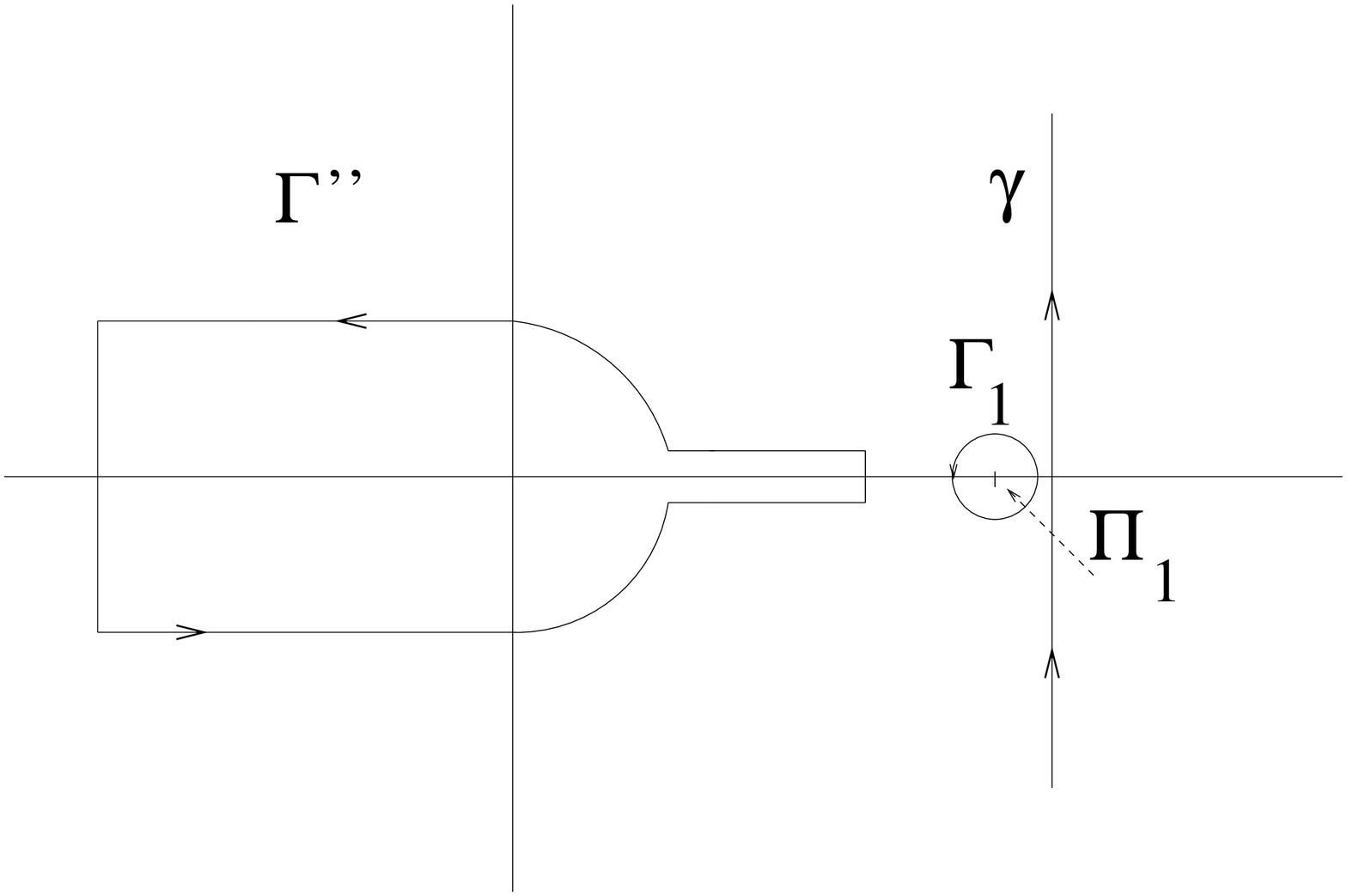,height=3.5cm, angle=0}
 \end{tabular}
 \caption{Contours $\Gamma $ and $\gamma$.
 \label{fig:contoursGUE}}
 \end{center}
\end{figure}

\brem \label{rem: x*}Here, we make some preliminary restrictions on $\eta$ and $R$, that will be fixed in the following Lemma. We assume that $\eta$ is small enough  so that the curve $x+i\eta$, $1\leq x\leq \frac{\pi_1+\pi^*}{2}$ crosses the circle of ray $\frac{C(\pi_1)}{2}.$ As $\frac{C(\pi_1)}{2}>1$,
we will then choose some $\eta \leq \sqrt{(\frac{C(\pi_1)}{2})^2-1}.$ Given such a $\eta,$ we call $x^*=x^*(\eta)=\frac{C(\pi_1)}{2}e^{i\theta_o}=x_o^*+i\eta$ this intersection.
Moreover, $R$ is chosen large enough to enclose all the $\pi_i, i=2, \ldots , r+1$.\erem
The crucial step in the proof of Fact \ref{claim H2} is the following Lemma.
\bl \label{Lem: GUEGamma'}
There exists $0<\eta \leq \sqrt{(\frac{C(\pi_1)}{2})^2-1}$, $R>0$ for which
\begin{itemize}
\item there exists $C=C(\eta)>0$ such that for any $ z\in \Gamma''_1\cup \Gamma''_2,$ $Re \left( F_{C(\pi_1)}(z)-F_{C(\pi_1)}(\pi_1)\right )\geq C >0.$ \item
$Re \left (F_{C(\pi_1)}\right )$ achieves its minimum on $\Gamma''_3\cup\Gamma''_4\cup \Gamma''_5$ at $x^*=x^*(\eta)$ defined in Remark \ref{rem: x*}.
\end{itemize}
\el
\paragraph{Proof of Lemma \ref{Lem: GUEGamma'} }Consider first $\Gamma''_1\cup \Gamma''_2$. The function
$x\mapsto Re \left (F_{C(\pi_1)}(x)-F_{C(\pi_1)}(\pi_1)\right )$ is decreasing on the interval $[\dfrac{1}{\pi_1}, \pi_1].$ Thus, for any $x\in [1, \pi^*],$
which is a compact interval of $(\dfrac{1}{\pi_1}, \pi_1) ,$
$$Re \left (F_{C(\pi_1)}(x)-F_{C(\pi_1)}(\pi_1)\right )\geq Re \left (F_{C(\pi_1)}(\pi^*)-F_{C(\pi_1)}(\pi_1)\right )\geq C_o>0.$$
%Furthermore, for all $x\in [\dfrac{1}{\pi_1}, \pi_1],$
%$\displaystyle{|Re \left (F_{C(\pi_1)}(x+iy)-F_{C(\pi_1)}(x)\right )|\leq C_1 |y|,}$
As $F'_{C(\pi_1)}$ is uniformly bounded in a compact set away from $0$,
we can now choose $\eta$ small enough so that
$\displaystyle{Re \left (F_{C(\pi_1)}(z)\right )\geq  F_{C(\pi_1)}(\pi_1)+ \frac{C_o}{2}, \: \forall z=x+iy, \quad \text{with }x\in[1,\pi^*], \:|y|\leq \eta.}$\\
Now, along $\Gamma''_3,$
$\displaystyle{\frac{d}{d\theta} Re \left (F_{C(\pi_1)}(\frac{C(\pi_1)}{2} e^{i\theta})\right )=\sin \theta C(\pi_1)^2/2(1-\cos \theta)>0,}$
since $\theta\geq \theta_o>0.$
Along $\Gamma''_4,$ for $z=iC(\pi_1)/2-x,$ $x>0$, 
$\displaystyle{\frac{d}{dx} Re \left ( F_{C(\pi_1)}(\frac{iC(\pi_1)}{2}-x)\right ) =C(\pi_1)+x+\frac{x}{|\frac{iC(\pi_1)}{2}-x|^2}>0.}$
Along $\Gamma''_5,$ and for $z=-R +it, t\leq \frac{C(\pi_1)}{2}$
$\displaystyle{Re \left (F_{C(\pi_1)}(z)\right )=\frac{R^2}{2}+C(\pi_1)R -\frac{t^2}{2}-\frac{1}{2}\log (|R+it|^2).}$
We can then choose $R$ large enough so that along $\Gamma''_5,$
$Re \left (F_{C(\pi_1)}(z)\right )\geq Re \left (F_{C(\pi_1)}(x^*)\right ).$ $\blacksquare$

\paragraph{}Now, we fix  $\eta$ and $R$ so that Lemma \ref{Lem: GUEGamma'} holds. Then, one has $$Re \left( F_{C(\pi_1)}(z)-F_{C(\pi_1)}(\pi_1)\right )\geq C >0 \:  \text{ and } Re (z-\tilde \pi_1)>\epsilon, \quad \forall z\in \Gamma''.$$ Using now the fact that $\Gamma''$ is a fixed (independent of $N$) length contour along which $|1/g|$ is uniformly bounded, it is then straightforward to obtain Fact \ref{claim H2} from Lemma \ref{Lem: GUEGamma'}.$\blacksquare$\\
 Combining Fact \ref{claim H2} with Fact \ref{Lem: H_1} gives formula (\ref{HNGGUE2}), which
finally proves Proposition \ref{Prop: GUE}.

\section{ Proof of Theorem \ref{theo: bulkaubord}}
In the whole Section, we assume that $\pi_1$ lies in a compact interval of $(1,\infty)$.
We further make the simplifying assumption
$$W_N=\text{ diag }(\pi_1,\ldots,\pi_1,0,\ldots,0),$$
 with $\pi_1$ of multiplicity $k_N$ such that
$\dfrac{k_N}{N}\rightarrow 0, k_N\rightarrow \infty$
as $N$ goes to infinity. The changes to be made in the case where $W_N$ admits eigenvalues between $0$ and $\pi_1$ (including the case where the number of these eigenvalues is increasing with $N$) will be indicated at the end of this section.

Let $\Gamma$ be a contour
encircling the poles $\pi_1$ and $0$, oriented counter clockwise and $\gamma=A+it, t\in \mathbb{R},$ such that $\Gamma\cap \gamma=\emptyset.$ Then the correlation kernel is now given by
$$K_N(u,v)=\frac{N}{(2i\pi)^2}
\int_{\Gamma}dz \int_{\gamma}dw e^{-N(z^2/2-uz)+N(w^2/2-wv)}\left (\frac{w}{z}\right
)^{N-k_N}\left (\frac{\pi_1-w}{\pi_1-z}\right )^{k_N}\frac{1}{w-z}.$$  Let us briefly indicate the idea of the proof of Theorem \ref{theo: bulkaubord}. Let
$C(\pi_1)$ and $\sigma^2(\pi_1)$ be defined by (\ref{sigma^2})
 and $\alpha_N$ be defined by (\ref{parametrealpha}).
The idea is to make a second order Taylor expansion around $\pi_1$. If $w=\pi_1+\alpha_N s$, and $u=C(\pi_1)+\dfrac{\alpha}{\sigma^2}\alpha_N$, for some $|\alpha|<2\sigma(\pi_1)$, then
the exact exponential term, $F_u$, defined by
\begin{equation}F_u(w):=w^2/2-wu+(1-\alpha_N^2)\log w +\alpha_N^2 \log(w-\pi_1)\label{Fureel}\end{equation}
satisfies $\displaystyle{F_u(\pi_1+\alpha_N s)=Ct(\pi_1)+k_N \left (\frac{s^2}{2\sigma^2}-\frac{\alpha s}{\sigma^2}+\log s+ \alpha_N G(s)\right ),}$
for some constant term $Ct(\pi_1)$ depending on $\pi_1$ and a function $G$ that should not grow much.
The function
$\displaystyle{ H(s)=\frac{s^2}{2\sigma^2}-\frac{\alpha s}{\sigma^2}+\log s}$ is then the exponential term of the correlation kernel (\ref{noyauJohansson}) of the GUE with parameter $\sigma^2=\sigma^2(\pi_1).$
Thus, suitably rescaled, the $k_N$ largest eigenvalues of the deformed Wigner ensemble should exhibit the same fluctuations as the
eigenvalues of a $k_N\times k_N$ GUE with parameter $\sigma^2$. That is what we now show.\\
Let then $\rho$ be the density of the semi-circular law with parameter $\sigma^2=\sigma^2(\pi_1)$
defined in (\ref{semicerclesigma^2}). Let $x_o, y_o$ be fixed and set
\begin{equation} \label{rescainitialbord}
u=C(\pi_1)+\dfrac{\alpha_Nx}{\sigma^2}, \:x=\alpha+\frac{x_o}{k_N \rho(\alpha)};\quad
 v=C(\pi_1)+\dfrac{\alpha_Ny}{\sigma^2},\: y=\alpha+\frac{y_o}{k_N \rho(\alpha)}.
\end{equation}
For $u,v$ satisfying (\ref{rescainitialbord}), we here consider the rescaled correlation kernel
\be \label{KN'bulk}
K_N'(x,y)=\frac{\alpha_N}{k_N \sigma^2 \rho(\alpha)}\exp{\{-N\frac{(x-y)}{\sigma^2}\alpha_N \pi_1\}}K_N(u,v).
\ee
The aim of the rest of this section is to obtain the following result.
\bp \label{prop: sinfinal}
Assume $\alpha=2\sigma \cos \theta_c$ in (\ref{rescainitialbord}), (\ref{KN'bulk}), with $ 0<|\theta_c|<\pi,$ and let $t_{c,\alpha}^{\pm}=\sigma \cos \theta_c$. Then,
\be \label{sinfinal}
\lim_{N \rightarrow \infty} K_N'(x,y)\exp{\left\{(y_o-x_o)Re\left (\frac{
t_{c,\alpha}^{+}}{\sigma^2}\right )\right\}}=\frac{\sin \pi(x_o-y_o)}{\pi(x_o-y_o)}.
\ee
\ep
\brem Theorem \ref{theo: bulkaubord} is then an easy consequence of Proposition \ref{prop: sinfinal} (see  e.g. \cite{DeiftZhou}, Section 6).\erem
Before beginning the proof of Proposition \ref{prop: sinfinal}, it is convenient to make here the following simplifying assumptions.
We assume that $N\geq N_o$, where $N_o$ is such that
\be \label{Nobulk}
\forall N\geq N_o, \forall |t| \leq 2\sigma+1, |\pi_1+\alpha_N t|\geq \frac{\pi_1}{2},|\pi_1+\alpha_N t-1|\geq \frac{\pi_1-1}{2}.
\ee
\subsection{Rewriting the kernel}
In this subsection, we first split the kernel into two subkernels, since
the idea is to prove that the asymptotics of $K_N'(x,y)$ is lead by the integral performed on a neighborhood of $\pi_1.$ Then
%This is the object of the first Proposition. In the second Proposition, 
we obtain an integral representation of these subkernels suitable for the saddle point analysis.\\

 Let then $\Gamma_1$ (resp. $\Gamma_2$) be
the circle of ray $\sigma$ (resp. $1$) centered at $\pi_1$ (resp. the origin). Both contours are oriented counterclockwise. Let $F_u$ be given by (\ref{Fureel}), and define the kernels
 \begin{eqnarray}
 &&K_{N,1}(u,v)=\alpha_N \exp{\{-N\alpha_N \pi_1 \frac{(x-y)}{\sigma^2}\}}
 \int_{\Gamma_1}dz \int_{\gamma}dw
 \exp{\{-NF_u(z)+NF_v(w)\}}\frac{1}{w-z}, \cr
 &&K_{N,2}(u,v)=\alpha_N \exp{\{-N\alpha_N \pi_1 \frac{(x-y)}{\sigma^2}\}}
 \int_{\Gamma_2}dz \int_{\gamma}dw
 \exp{\{-NF_u(z)+NF_v(w)\}}\frac{1}{w-z}.\end{eqnarray}

\bp \label{Prop:expressionoubliee} Let $K_N'(x,y)$ be given by (\ref{KN'bulk}). Then,
$K_N'(x,y)=K_{N,1}'(x,y)+K_{N,2}'(x,y),$ where
\be K'_{N,1}(x,y)=\frac{1}{k_N \sigma^2 \rho(\alpha)(2i\pi)^2}K_{N,1}(u,v)\text{  and  }
K'_{N,2}(x,y)=\frac{1}{k_N \sigma^2 \rho(\alpha)(2i\pi)^2}K_{N,2}(u,v).\ee
\ep
As $x\simeq y \simeq \alpha$ in (\ref{rescainitialbord}), it is not hard to see that the
two integrands $F_u(w)$, $F_v(w)$ have the same critical points lying
around $\pi_1.$ While this should not cause any trouble for the saddle point analysis of $K_{N,2}'$, this prevents that of $K_{N,1}',$ because of the singularity
$\displaystyle{\frac{1}{w-z}}.$ Thus, we have to remove the singularity of the kernel $ K_{N,1}'.$ This is the object of the following Proposition.\\

Set, for $s\in \mathbb{C}$ such that $Re(\pi_1+\alpha_N sx)>0,$ $\forall x\in[0,1]$, \be G(s)=\alpha_N^2s^3\int_0^1 \frac{(1-x)^2}{(\pi_1+\alpha_Nsx)^3}dx -s\int_0^1 \frac{1}{\pi_1+s\alpha_Nx}dx.\label{definitionGpert}\ee
\bp \label{prop: noyaupresqueGUE}Assume $N \geq N_o$, with $N_o$ defined in (\ref{Nobulk}), then with the rescalings (\ref{rescainitialbord}),
\begin{eqnarray}& K_{N,1}'(x,y)&=
\frac{k_N }{(2i\pi)^2(y_o-x_o)}\int_{\Gamma'_1} \int_{\gamma'}\exp{\left\{k_N\left (\frac{t^2-2yt}{2\sigma^2}+\alpha_N G(t)-\frac{s^2-2sx}{2\sigma^2}-\alpha_N G(s)\right )\right \}}   \cr &&\times\left (\frac{t}{s}\right )^{k_N}\left (1-\exp{\{\frac{s(y_o-x_o)}{\sigma^2 \rho (\alpha)}\}}\right )
\frac{1}{s}\left (\frac{s+t-y}{\sigma^2}+\alpha_N\left (\frac{tG'(t)-sG'(s)}{t-s}\right )\right)ds dt,\cr
&&
\end{eqnarray}
where
$\Gamma'_1$ is a circle of ray $\sigma$ around the origin  and $\gamma'=A+i\mathbb{R}, $ with $A \geq -2\sigma-1.$
\ep
\brem
%The above integral representation of $K_{N,1}'$ has no more singularity. Thus 
$\Gamma'_1$ can now cross $\gamma'.$
\erem

\paragraph{Proof of Proposition \ref{prop: noyaupresqueGUE}:}
Assume that $\gamma''=A+i\mathbb{R}, A>0 $ large enough not to cross a circle of ray $\sigma$ around $\pi_1.$ We first show the formula
\begin{eqnarray}&K_{N,1}(u,v)=&\frac{k_N}{(2i\pi)^2}\int_{\Gamma'_1} ds \int_{\gamma''} dt\left (\frac{t}{s}\right )^{k_N}\frac{1}{t-s}\cr
 & &\times
\exp{ \left\{k_N\left (\frac{t^2-2yt}{2\sigma^2}+\alpha_N G(t)-\frac{s^2-2sx}{2\sigma^2}-\alpha_N G(s)\right )\right \}}
 .\label{termeentropGUE}
\end{eqnarray}

Define $\displaystyle{{\tilde F_u}(z):=z^2/2-uz+\log z.}$
Here we choose the principal branch of the logarithm. We now make the change of variables $z=\pi_1+\alpha_Ns.$
Then one has that
${\tilde F_u}(\pi_1+\alpha_Ns)=\tilde F_{C(\pi_1)}(\pi_1+\alpha_Ns)-(u-C(\pi_1))(\pi_1+\alpha_Ns).$
Performing now a Taylor expansion for the real and imaginary part gives
\be \tilde F_{C(\pi_1)}(\pi_1+\alpha_N s)=\tilde F_{C(\pi_1)}(\pi_1)+\frac{{\tilde F}''_{C(\pi_1)}(\pi_1)}{2} \alpha_N^2 s^2+\alpha_N^3 s^3\int_0^1 \frac{(1-x)^2}{(\alpha_N s x +\pi_1)^3}dx.
\ee
Finally, as $\pi_1+\alpha_N s$ does not reach $\mathbb{R}_-$, because of (\ref{Nobulk}), we can write \\
$\displaystyle{\alpha_N^{k_N}\frac{s^{k_N}}{(\pi_1+\alpha_N s)^{k_N}}=\alpha_N^{k_N}s^{k_N}\exp{\{-k_N \left(\log(\pi_1)+\alpha_N s \int_0^1 \frac{1}{\pi_1+\alpha_N sx}dx\right)\}}.}$
Thus we obtain (\ref{termeentropGUE}) for contours $\Gamma'_1$ and $\gamma''$ chosen as above (as neither $\pi_1+\alpha_N t$ nor $\pi_1+\alpha_N s$ reaches $\mathbb{R}_-$)
.\\ 
Finally, we use the same method as in \cite{Johansson} to remove the singularity. In
(\ref{termeentropGUE}), we make the change of variables $s\mapsto \beta s$, $t\mapsto \beta t$ for
$\beta $ real close to one. Thanks to Cauchy's theorem, we can deform back these contours to
$\gamma$ and $\Gamma.$ Taking the derivative at $\beta=1$ gives
\begin{eqnarray}
&K_{N,1}(u,v)=&-\frac{k_N^2}{(2i\pi)^2}\int_{\Gamma'_1}ds\int_{\gamma''}dt  \exp{\left\{k_N \left(-\frac{s^2-2sx}{2\sigma^2}-\alpha_N G(s)+\frac{t^2-2ty}{2\sigma^2}+\alpha_N G(t)\right)\right\}} \cr
&&\times \Bigl (\frac{t^2-s^2}{\sigma^2}+\frac{xs-yt}{\sigma^2}+\alpha_N tG'(t)-\alpha_N sG'(s)\Bigr )\left( \frac{t}{s}\right)^{k_N}\frac{1}{t-s}.
\end{eqnarray}
Now this gives, using (\ref{termeentropGUE}) and for the rescalings (\ref{rescainitialbord}),
\begin{eqnarray}
&\dfrac{d((\frac{x}{\sigma^2}-\frac{y}{\sigma^2})K_{N,1}'(x,y))}{d(\frac{x}{\sigma^2})}=&\!\!-\frac{k_N^2}{(2i\pi)^2}\int_{\Gamma'_1}\int_{\gamma''} \exp{\left\{k_N \left(\frac{t^2-2ty}{2\sigma^2}+\alpha_N G(t)-\frac{s^2-2sx}{2\sigma^2}-\alpha_N G(s)\right)\right\}}\cr
&&\quad \quad \:\:\left (\frac{t}{s}\right)^{k_N}\left (\frac{s+t-y}{\sigma^2}+\alpha_N\frac{tG'(t)-sG'(s)}{t-s}\right )ds dt.\label{lebonaint}
\end{eqnarray}
Solving (\ref{lebonaint}) with an integration by parts, we obtain finally Proposition \ref{prop: noyaupresqueGUE} (we can then move $\gamma''$ to $\gamma'$).  %The two contours $\gamma'$ and $\Gamma_1$ can now cross but must not cross the real axis at some point $t$ such that $\pi_1+\alpha_N t<0.$ 
$\blacksquare$

\subsection{A study of critical points }
 In this part, under Assumption \ref{Hyp2}, we show that the exact critical points of the integrands, in $K_{N,1}'$ and $K_{N,2}'$, lie on a curve that is almost the circle of ray
$\sigma(\pi_1)$ around $\pi_1$, provided $\alpha_N$ tends to $0$. Furthermore, we prove that the relevant critical points for the saddle point analysis are well approximated by those of $H_{\alpha/\sigma^2}$ if
\be H_{u}(t):=\frac{t^2}{2\sigma^2}-ut+\log t.
\label{defdeH} \ee 
Consider the exact exponential term to be analyzed,
$F_u(w):=w^2/2-uw+(1-\alpha_N^2)\log(w)+\alpha_N^2 \log(w-\pi_1).$
The equation
$F_u'(w)=w-u+(1-\alpha_N^2)/w+\alpha_N^2/(w-\pi_1)=0$
admits three solutions. One is real, in the interval $(0,\pi_1),$ and two others $w_N^{\pm}$ that are conjugate. We now study these critical points $w_N^{\pm}$.
It is an easy fact that any critical point $w$ for $F_u$ with non zero imaginary part satisfies the equation
\be\label{eqcourbes}\frac{1-\frac{1}{|w|^2}}{1-\frac{1}{|w-\pi_1|^2}}=\frac{-\alpha_N^2}{1-\alpha_N^2}.
\ee
Then the solution of (\ref{eqcourbes}) define one or two (depending on $\alpha_N^2$) curves encircling $0$ and $\pi_1$.
%Let $0<\epsilon:=\frac{\alpha_N^2}{1-\alpha_N^2}<\infty$, otherwise the result is straightforward.
% We have to find, given an $x$, $y(x)$ such that
%$$1+\epsilon=\frac{1}{x^2+y^2}+\frac{\epsilon}{(x-\pi_1)^2+y^2}\quad (E).$$

%If $\frac{1}{x^2}+\frac{\epsilon}{(x-\pi_1)^2}>1+\epsilon,$ then $(E)$ admits two %solutions $y(x)$ and $-y(x)$. $(E)$ has a single solution $y(x)=0$ if
%$\frac{1}{x^2}+\frac{\epsilon}{(x-\pi_1)^2}=1+\epsilon$ and none otherwise.
%Now the equation $\frac{1}{x^2}+\frac{\epsilon}{(x-\pi_1)^2}=1+\epsilon$
%admits four real roots (note that this equation is $F_u''(x)=0$), provided $\epsilon %\not=0,$ and $\epsilon<\infty$. $\blacksquare$

\paragraph{}Consider now a sequence $\alpha_N$ such that $\displaystyle{\lim_{N \rightarrow \infty}\alpha_N=0}$. We now show that critical points for $F_u$ almost lie on the curve $C_1=\{\pi_1+\alpha_N \sigma e^{i\theta}, 0\leq \theta\leq 2\pi\},$ where $\sigma=\sigma(\pi_1)$ in the following.

\bl \label{Lemm: loccriticalpoints} Let $u$ be given by (\ref{rescainitialbord}) with
$\alpha=2\sigma \cos(\theta_c),0<|\theta_c|<\pi $. Then, the critical
points $w_{N}^{\pm}$ are non real and $\exists\: C(\pi_1)>0$ such that 
$w_N^{\pm}=\pi_1+\alpha_N t_N^{\pm}$ with $|t_N^{\pm}-\sigma e^{i\pm\theta_c}|\leq
C(\pi_1)\alpha_N+ \frac{|x_o|}{k_N \rho(\alpha)}.$

\el
\paragraph{Proof of Lemma \ref{Lemm: loccriticalpoints}:}If $u=C(\pi_1)+\alpha_N\alpha/\sigma^2$, then
\be
F_u'(\pi_1+\alpha_Nt)=\alpha_N\left (H'_{\frac{\alpha}{\sigma^2}}(t)+\alpha_NG'(t)\right)=\alpha_N\left (H'_{\frac{\alpha}{\sigma^2}}(t)+\alpha_N \frac{t^2-\pi_1^2}{\pi_1^2(\pi_1+\alpha_Nt)}\right ),
\ee
with $\displaystyle{H'_{\frac{\alpha}{\sigma^2}}(t)=\frac{t}{\sigma^2}-\frac{\alpha}{\sigma^2}+\frac{1}{t}.}$
Set now $T_o=max\{ 2\sigma(\pi_1), 4\pi_1^2\}$. As $\pi_1$ lies in a compact set of $(1, \infty)$, it is not hard to see that, for $|t|<T_o$, and $N$ large enough so that $\alpha_NT_o<\pi_1/2$, there exists $C(\pi_1)>0$ such that
$\displaystyle{\Big|\bigl (\frac{t^2-\pi_1^2}{\pi_1^2(\pi_1+\alpha_Nt)}\bigr )^{(l)}\Big|\leq C( \pi_1), \quad 0\leq l\leq 4.}$ Thus, if now $u$ is now given as in (\ref{rescainitialbord}),
\be\label{yauneracine}\Big|\frac{F_u'(\pi_1+\alpha_N t)}{\alpha_N}-H'_{\frac{\alpha}{\sigma^2}}(t)\Big|\leq \alpha_NC(\pi_1)+\frac{|x_o|}{k_N \rho(\alpha)} ,\quad
|F_u''(\pi_1+\alpha_N t)-H''_{\frac{\alpha}{\sigma^2}}(t)|\leq \alpha_NC(\pi_1).\ee
Now, if $\alpha=2\sigma\cos \theta_c,$ with $ 0<|\theta_c | <\pi,$  $H_{\frac{\alpha}{\sigma^2}}$ admits  two critical points that are conjugate, and given by
$ t_{c,\alpha}^{\pm}=\sigma e^{i\pm\theta_c}.$
Thus using (\ref{yauneracine}), we obtain Lemma \ref{Lemm: loccriticalpoints}.
$\blacksquare$

\subsection{Estimate for $K_{N,1}'$}
This subsection is devoted to the proof of the following Proposition. Let $K_{N,1}'$ be the kernel defined in Proposition \ref{Prop:expressionoubliee}.
\bp \label{Prop: sinus}Assume $\alpha=2\sigma \cos \theta_c,$ with $0<|\theta_c|<\pi,$ and let $t_{c,\alpha}^{\pm}=\sigma e^{i\pm\theta_c}$. $$\lim_{N \rightarrow \infty} K_{N,1}'(x,y)\exp{\left\{(y_o-x_o)Re\left (\frac{
t_{c,\alpha}^{+}}{\sigma^2}\right )\right\}}=\frac{\sin \pi(x_o-y_o)}{\pi(x_o-y_o)}.$$ \ep

\paragraph{Proof of Proposition \ref{Prop: sinus} } The proof is organized as follows. As the correlation kernel given in Proposition \ref{prop: noyaupresqueGUE} is not of the form (\ref{formeproduit}), we analyze the double integral "simultaneously". First we define ascent and descent contours for $H_{\alpha/\sigma^2}$ and show that the perturbative terms, due to $G$ defined in (\ref{definitionGpert}), do not grow too much.  We then slightly deform these contours to go through the effective critical points of $F_v,$ so that we can then perform the saddle point analysis.
\brem From now on, as $t_{c,\alpha}^{\pm},$ as well as $ t_N^{\pm},$ are conjugate, we may drop the $\pm$ sign (when possible) in the following, if results proved for $t_N^+$ hold for $t_N^-$ up to conjugation.
\erem
\paragraph{}Set $\gamma'=t_{c,\alpha}^++it, t\in \mathbb{R}, $ oriented from bottom to top. Let also $0<\epsilon<<Im(t_{c,\alpha}^+)$ be given. 
\bl \label{Lem:avtdec} One has $\max\bigl\{\Big|e^{k_N H_{\frac{\alpha}{\sigma^2}}(t_{c,\alpha}^++it)}\Big|, -Im(t_{c,\alpha}^+)\leq t\leq -Im(t_{c,\alpha}^+)+\epsilon\bigr\}=\Big|e^{k_N H_{\frac{\alpha}{\sigma^2}}(Re(t_{c,\alpha}^+)+i\epsilon)}\Big|$ and 
there exists $c_o>0$ such that $\Big|e^{k_N H_{\frac{\alpha}{\sigma^2}}( t_{c,\alpha}^++it)}\Big|\leq \Big|e^{k_N H_{\frac{\alpha}{\sigma^2}}(t_{c,\alpha}^+)}\Big|e^{-c_ok_Nt^2},$ $\forall t\in [-Im(t_{c,\alpha}^+)+\epsilon, \infty].$
 \el
\paragraph{Proof :}This follows from the fact that
$\displaystyle{\frac{d}{dt} \log \Big|e^{H_{\frac{\alpha}{\sigma^2}}(\sigma \cos \theta_c+it )}\Big|=-t\left (\frac{1}{\sigma^2}-\frac{1}{\sigma^2 \cos^2{\theta_c}+t^2}\right )}$ ($t>0$ if $\theta_c=\frac{\pi}{2}).$ And $t\mapsto  |e^{H(\sigma \cos \theta_c+it )}|, 0\leq t\leq \epsilon,$ is a decreasing function if $\epsilon$ is small enough.$\blacksquare$

\paragraph{}We now show that $Re (F_{u})$ decreases faster than (resp. almost as) $H_{\frac{\alpha}{\sigma^2}}$ on $\gamma'$, if $t>0$ is large enough and $Re \left(t_{c,\alpha}^{\pm}\right) \geq 0$ (resp $Re \left(t_{c,\alpha}^{\pm}\right) <0$). Let $\epsilon$ be as in Lemma \ref{Lem:avtdec}, $\eta>0$ (small) be given.
\bl \label{Lemma: decroissance} There exist $T_o>0,$ $N_1$ depending
on $\pi_1$ only, $C_o>0$, $C_{T_o}>0$ such that, for $N \geq N_1$, 
\begin{eqnarray}
&&\!\!\!\!\!\!\!\!\big |e^{\{N F_u(\pi_1+\alpha_N(t_{c,\alpha}^{+}+it))\}}\big |\leq e^{\{NRe \left (
Fu(\pi_1+\alpha_Nt_{c,\alpha}^{+})\right )-k_NC_o\epsilon^2/8\}},\text{ } t\in
[-Im(t_{c,\alpha}^{+}),-Im(t_{c,\alpha}^{+})+\epsilon],\label{termgenant}\\
&&\!\!\!\!\!\!\!\!\Big |e^{\{N F_u(\pi_1+\alpha_N(t_{c,\alpha}^{+}+it))\}}\Big |\leq e^{\{NRe \left (
Fu(\pi_1+\alpha_Nt_{c,\alpha}^{+})\right )-k_NC_ot^2/4\}},\text{ } t\in
[-Im(t_{c,\alpha}^{+})+\epsilon,-\eta]\cup [\eta, T_o],\label{petitepsilon1}\\ &&\!\!\!\!\!\!\!\!\Big|e^{N
\{F_u(\pi_1+\alpha_N(t_{c,\alpha}^++it))\}}\Big |\leq \Big |e^{N
\{
Fu(\pi_1+\alpha_Nt_{c,\alpha}^+ )-k_NC_oT_o^2/4-k_NC_{T_o}(t^2-T_o^2)/4\}}\Big |, \text{ } T_o \leq t. \label{petitepsilon}
\end{eqnarray}
\el
\paragraph{Proof of Lemma \ref{Lemma: decroissance}: }
We first examine the case where $Re (t_{c,\alpha}^{\pm})=\sigma cos \theta_c > 0.$ 
Using that for $t>T_o= \max(4\pi_1^2, 2\sigma(\pi_1))$, $Im ( G'(t))> 0,$ we obtain that for $t>T_o$,
\be \label{decplusrap}\frac{d}{dt} Re \left (\frac{F_u(\pi_1+\alpha_N(t_{c,\alpha}^++it))}{\alpha_N^2}\right )<-Im \left (H'_{\alpha}(t_{c,\alpha}^++it)\right )\leq -C_{T_o}t,\ee
where $C_{T_o}=1/\sigma^2-1/|t_{c,\alpha}^++iT_o|^2\geq
1/\sigma^2-1/|\sigma e^{i \theta_c}+iT_o|^2>0$. Now, $G, G'$ are uniformly bounded on a compact set $K$ (independant of $N$) containing $\gamma'\cap\{|Im (w)| \leq T_o+2\sigma
\}$. Thus, using Lemma \ref{Lem:avtdec}, we know that $\exists C_o>0$ such that, for $N$ large enough,
$$\frac{d}{dt}Re \left( H_{\frac{\alpha}{\sigma^2}}(t_{c, \alpha}^++it)+\alpha_NG(t_{c, \alpha}^++it)\right)\leq -C_ot/2, \forall t\in
[-Im(t_{c,\alpha}^{+})+\epsilon,-\eta]\cup [ \eta, T_o] .$$
This gives (\ref{petitepsilon1}). The fact that $G$ is bounded on $K$ also gives that (\ref{termgenant}) holds for $N$ large enough. Combining (\ref{petitepsilon1}) with (\ref{decplusrap}) gives then (\ref{petitepsilon}). 
This proves Lemma \ref{Lemma: decroissance} in this case.\\
If $Re (t_{c,\alpha}^{\pm})=\sigma \cos \theta_c\leq 0, $ (\ref{termgenant}) and (\ref{petitepsilon1}) are proved as above.  
One can then check that 
 $\exists \:C(\pi_1)>0$ such that $\displaystyle{Im \left (G'(Re (t_{c,\alpha}^+) +iT)\right
)\geq -C(\pi_1)T,}$ provided $Re (\pi_1+\alpha_Nt_{c,\alpha}^+) \geq \pi_1/2.$ This holds
for $N$ large enough and we can then find $N_1>0,$ such that
$\displaystyle{C_{T_o}-C(\pi_1)\alpha_N >\frac{C_{T_o}}{2}, \:\forall N \geq N_1.}$
Thus for $N \geq N_1,$ and $t\geq T_o$, one has that $\displaystyle{\frac{d}{dt} Re \left (\frac{F_u(\pi_1+\alpha_N(t_{c,\alpha}^++it))}{\alpha_N^2}\right )<-\frac{C_{T_o}t}{2}.}$ This finishes the proof of Lemma
\ref{Lemma: decroissance}.$\blacksquare$

\paragraph{}We now turn to the second contour.
Define then $\Gamma'_1=\sigma e^{i\theta}, \theta \in [0,2\pi],$ oriented counterclockwise. Note that $\Gamma'_1$ describes the curve of critical points for $H_{x}$ when $x$ describes $[-2\sigma, 2\sigma].$
\bl \label{LemGAMMA'_1}
Assume that $\alpha=2\sigma \cos \theta_c$. Then, there exists $c_o>0$ such that, for any $\theta\in[0,2\pi],$\\
$|e^{-k_NH_{\alpha/\sigma^2}(\sigma e^{i\theta})}|\leq |e^{-k_NH_{\alpha/\sigma^2}(\sigma e^{i\theta_c})}|e^{-k_N c_o(\theta-\theta_c)^2}.$ 

\el
\paragraph{Proof: }If $|\theta|<\pi,$
$\dfrac{d}{d\theta}Re \left (H_{\alpha/\sigma^2}(\sigma e^{i\theta})\right )=2\sin \theta(\cos \theta_c-\cos\theta) .$ The computation for $\theta=\pi$ is here left.
$\blacksquare \hfill$\\
As the contour $\Gamma'_1$ lies in a fixed compact set away from the singularities of $G$, we know that the contribution of $G$ will not perturb the saddle point analysis on $\Gamma'_1.$\\

Before performing the asymptotic expansion of $K'_{N,1}(x,y),$ one should take care of the remaining terms, which should not explode due to the perturbation $G$.
Set
\begin{eqnarray}
&&h(s,t)=\frac{\exp{\{\frac{s(x_o-y_o)}{\sigma^2 \rho(\alpha)}\}}-1}{y_o-x_o}, \: K_y(s)=H_{\frac{y}{\sigma^2}}(s)+\alpha_N G(s),\label{defdehpoursinus}\\
&&g(s,t)=\frac{1}{s}\left (\frac{s+t-y}{\sigma^2}+\alpha_N\frac{tG'(t)-sG'(s)}{t-s}\right)=
\frac{tK_y'(t)-sK_y'(s)}{s(t-s)}.\label{defdegpoursinus}
\end{eqnarray}
Then $\displaystyle{ K_{N,1}'(x,y)=\frac{k_N}{(2i\pi)^2}\int_{\Gamma'}\int_{\gamma}h(s,t)g(s,t)\exp{\{-k_N K_y(s)+k_N K_y(t)\}}. }$\\
We have to check that the function $g(s,t)$ will not perturb the saddle point analysis.
As the contour $\Gamma'_1$ is compact and for $|Im (w)|\leq T_o$, the functions $G(t), G'(t)$ are bounded by some constant depending on $\pi_1$ only. Thus $g(s,t)$ is bounded on $\Gamma'_1\cup \left (\gamma\cap \{|Im (w)|\leq T_o\}\right ).$ Note also that
along $\gamma'$,
$\displaystyle{\frac{1}{|\pi_1+\alpha_N t|}\leq \frac{2}{\pi_1}}$ so that $|G'(t)|\leq \alpha_N t^2.$
%Now it is easy to see that $g(s,t)$ behaves as a polynomial of $|t|$.
Thus, there exists some constant $C>0$ such that, for $t>T_o$, using Lemma \ref{Lemma: decroissance},
$$|g(s,t)|\Big| \frac{\exp{\{N F_u(\pi_1+\alpha_N(t_{c,\alpha}+it))\}}}{\exp{\{NF_u(\pi_1+\alpha_Nt_{c,\alpha})\}}}\Big|
\leq Ct^3\exp{\left\{-Ck_N\frac{t^2}{4}\right \}}\leq \exp{\left\{-Ck_N\frac{t^2}{8}\right \}},$$
for $N$ large enough. This is the needed estimate to perform the saddle point analysis.\\
Now, and this is the core of the argument, we slightly deform the contours $\gamma$ and $\Gamma'_1$ to contours $\gamma_N$ and $\Gamma_N$ that go through the effective critical points $t_N^{\pm}$ of $K_y$. By Lemma \ref{Lemm: loccriticalpoints}, these contours lie within a $C^1$ distance of $\gamma$ (resp. $\Gamma'_1$) smaller than $C\alpha_N$ for some constant $C>0.$ Furthermore,
 $\gamma_N$ and $\Gamma_N$ coincide with $\gamma$ and $\Gamma'_1$ outside the disks $|t-t_{c,\alpha}^{\pm}|<\eta$ ($\eta $ small).
Then , by Proposition \ref{prop: noyaupresqueGUE}, Lemmas \ref{Lem:avtdec}, \ref{Lemma: decroissance}, \ref{LemGAMMA'_1} and (\ref{petitepsilon}), we obtain, by a standard saddle point argument that

\be\lim_{N \rightarrow \infty} K_{N,1}'(x,y)=\sum_{b,d=\pm 1}\frac{ sgn(b)}{(2i\pi)^2}\frac{2\pi \exp{\{k_N (K_y(t_N^b)-K_y(t_N^d))\}}}{i\sqrt{K_y''(t_N^{b})K_y''(t_N^{d})}}g(t_N^{b},t_N^{d})h(t_N^{b},t_N^{d}).\label{eqsinus}\ee

Now, the critical points are conjugate, thus $K_y(t_N^{+})=\overline {K_y(t_N^{-})}.$
Using (\ref{defdegpoursinus}), one can check that $$g(t_N^+,t_N^-)=g(t_N^-,t_N^+)=0, \quad g(t_N^{\pm},t_N^{\pm})=K_y''(t_N^{\pm}),$$ so that only the contribution of equal critical points have to taken into account in (\ref{eqsinus}).
By Lemma \ref{Lemm: loccriticalpoints}, one has $Im(t_N^{\pm})=\pm\pi \sigma^2 \rho(\alpha)+o(1)$, so that for $h$ given by (\ref{defdehpoursinus}),\\
 $\dfrac{h(t_N^-,t_N^-)-h(t_N^+,t_N^+)}{2i}\exp{\{(y_o-x_o) Re (t_N^+/\sigma^2)\}}=\dfrac{\sin \pi(x_o-y_o)}{\pi(x_o-y_o)}.$
This yields Proposition \ref{Prop: sinus}.$\blacksquare$
%$$\lim_{N \rightarrow \infty} K_{N,1}'(x,y)\exp{\{(y_o-x_o)Re (t_N^{+}/\sigma^2)\}}=\frac{\sin \pi(x_o-y_o)}{\pi(x_o-y_o)}.\blacksquare$$
%This finishes the proof of 

\subsection{Estimate for $K'_{N,2}(x,y)$}

This subsection is devoted to the proof of the following Proposition. Let $K_{N,2}'$ be the kernel defined in Proposition \ref{Prop:expressionoubliee}.
\bp \label{Prop: restesinus} There exists $C_o>0,\: N_o>0$ such that  
$\displaystyle{\:\:|K'_{N,2}(x,y)|\leq \exp{\{-C_oN/2\}}, \: \forall N \geq N_o.}$
\ep
\paragraph{Proof of Proposition \ref{Prop: restesinus}:} We first show that the function $\dfrac{1}{w-z}$ is bounded as $z\in \Gamma_2$ and $w=\pi_1+\alpha_N t, t\in \gamma'.$
By (\ref{Nobulk}), we can assume that the image of $\gamma'$ under the map $t\mapsto \pi_1+\alpha_N t$ lies in the half plane $Re (w)> (\pi_1+1)/2.$ Thus, for $z\in \Gamma_2$ and $w=\pi_1+\alpha_N t, t\in \gamma'$, $\dfrac{1}{|w-z|}\leq \dfrac{2}{\pi_1-1}.$
Now, for $N$ large enough, $\displaystyle{\min_{\Gamma_2}Re (F_u(z))=Re(F_u(1))}$ and $1$ lies in a compact set of $(1/\pi_1, \pi_1)$. Then, we have that
(as $t_N^+=\overline{t_N^-}$, we can consider ${t_N}^+ $ only, and drop the $+$ sign from now on)
\begin{eqnarray}
&&\exp{\left \{-N F_u(1)+N F_v(\pi_1+\alpha_N t_N)\right \}}=\cr
&&\frac{\exp{\left \{N \left ((\pi_1+\alpha_N t_N)^2/2-C(\pi_1)(\pi_1+\alpha_N t_N)\right )\right \}}}{\exp{\left \{N(1/2-C(\pi_1))\right \}}}\left (\pi_1+\alpha_N t_N \right )^N\cr
&&\times \exp{\{N (C(\pi_1)-v)(\pi_1+\alpha_N t_N)-N(C(\pi_1)-u)\}}\left (\pi_1+\alpha_N t_N\right )^{-k_N}\left(\frac{\alpha_N t_N}{1-\pi_1}\right )^{k_N}\label{l;o;t;}
\end{eqnarray}
Now, it is easy to see that $|(\ref{l;o;t;})|\leq e^{C\alpha_N N},$ for some constant $C>0$.
Finally, using that
$\displaystyle{\left (\frac{\pi_1+\alpha_N t_N}{\pi_1}\right)^N =\exp{\left\{N\alpha_N \int_{0}^{t_N/\pi_1}\frac{du}{1+\alpha_N u}\right \}}\leq \exp{\{N \alpha_N C'\}},}$
for some constant $C'>0$ and $N$ large enough,
we obtain that there exists a constant $C$ and $N_o$ such that  for $N\geq N_o,$
$$
\Big|\exp{\{-N F_u(1)+N F_v(\pi_1+\alpha_N t_N)\}}\Big|\leq \exp{\{N\left (\pi_1^2/2-C(\pi_1)\pi_1\right )-N\left (1/2-C(\pi_1)\right )+C \alpha_N N\}}\pi_1^N.
$$
Now, $\exists \:C_o>0$ such that $\displaystyle{\exp{\{N\left(\pi_1^2/2-C(\pi_1)\pi_1\right)\}}\pi_1^N \exp{\{-N\left (1/2-C(\pi_1)\right)\}}\leq \exp{\{-C_o N\}}.}$ This follows from the fact that the function $f: x\mapsto x^2/2-C(\pi_1)x
+\log x, x\geq 1$ is strictly decreasing in the interval $(1/\pi_1,\pi_1),$
as $\pi_1 $ lies in a compact interval of $(1,\infty).$
Therefore for $N$ large enough
$\displaystyle{K_{N,2}'(x,y)\leq \exp{\{-C_oN/2\}}. \blacksquare}$\\
Finally, combining Propositions \ref{Prop:expressionoubliee}, \ref{Prop: sinus} and Proposition \ref{Prop: restesinus} yields Theorem \ref{theo:
bulkaubord}.
\subsection{Extensions \label{subsec: extensionsbulk}}
We now explain the changes to be made to prove Theorem \ref{theo: bulkaubord} in the case where $W_N$ has eigenvalues $\pi_i,i=2,\ldots, r_N+1$ distinct of $0$, under Assumption \ref{Hyp2}.
The exponential term to be analyzed is given by
\be \label{tildefuext}\tilde F_u(w)=F_u(w)-\beta_N \log(w)+\frac{1}{N}\sum_{i=1}^{N \beta_N}\log(w-\pi_{i+1}),\ee
where $F_u$ is given as in (\ref{Fureel}).
Let also $u$ be given as in Definition \ref{def: spacing}.
Then, there exist constants, depending on $\pi_1$ only, such that, for all $t$ in a given compact set of $\mathbb{C}^*$,  \be \displaystyle{N\tilde F_u(\pi_1+\alpha_N t)=N Ct(\pi_1)+\beta_N Ct'(\pi_1)+k_N H_{\frac{\alpha}{\sigma^2}}(t)+k_N O(\alpha_N +\beta_N)}.\label{cont1bulk}\ee
Let then define $\displaystyle{G_1(t):=\dfrac{1}{k_N}\left (N\tilde F_u(\pi_1+\alpha_N t)-N Ct(\pi_1)-\beta_N Ct'(\pi_1)-k_N H_{\frac{\alpha}{\sigma^2}}(t)\right),}$
which plays the role of the function $\alpha_N G $ defined in (\ref{definitionGpert}).
Let also $\tilde {t}_N^{\pm}$ be the critical points for $t\mapsto \tilde F_u(\pi_1+\alpha_N t)$, and $\rho$ be the density of the semi-circular law with parameter $\sigma^2(\pi_1)$ as before.
As $G_1$ and its three first derivatives have no singularity in a given compact neighborhood $K_o$ of $0$, we readily have that
\be \label{cont3bulk}Im(\tilde {t}_N^{\pm})= \pi \rho(\alpha)+O(\alpha_N+\beta_N).\ee
Furthermore, defining $\displaystyle{ u_i=C(\pi_1)+\dfrac{\alpha}{\sigma^2}\alpha_N,}$
and given any compact set $K$ of $\mathbb{C}\setminus\{0, \pi_2, \ldots, \pi_{r+1}\}$, it is easy to check that there exists a constant $C(K)$ such that
\be \label{cont2bulk}\displaystyle{|\tilde F_{u}^{(l)} (w) -F_{u_i}^{(l)}(w) | \leq C(K) \beta_N,\:\forall w\in K},\: l=0,\ldots,3.\ee
Now formulas (\ref{cont1bulk}), (\ref{cont3bulk}) and (\ref{cont2bulk}) readily give that the asymptotics of $K_{N,1}'$ is unchanged. One simply replaces the function $\alpha_N G(\cdot)$ with the function $G_1(\cdot)$ in the proof of Proposition \ref{Prop: sinus}.
For the proof of Proposition \ref{Prop: restesinus}, we choose $\Gamma'_2$
to be the circle of ray $\pi^*=\max\{\pi_2+(\pi_1-\pi_2)/2,1\}$ completed by some   contour encircling the $\pi_i<0$. The latter contour lies in a fixed
compact set $K$ of $\mathbb{C}\setminus\{0, \pi_2, \ldots, \pi_{r+1}\}$, by Assumption \ref{Hyp2}. Then $\displaystyle{Re (\tilde F_u(w))> Re (\tilde F_u(\pi^*)) -C(K)\beta_N,}$  $\forall w\in \Gamma'_2.$ The fact that $Re (\tilde F_u(\pi^*))>Re (\tilde F_u(\pi_1+\alpha_N \tilde{t}_N^{\pm}))$ now follows from the same arguments as in the proof of Proposition \ref{Prop: restesinus}. This finishes the proof of Theorem \ref{theo: bulkaubord} in this case.

\section{Proofs of Theorem \ref{theo: rhoN<<3/7} and Theorem \ref{theopi_1<1}\label{sec: preuvebordgrrang}}
In this Section we first prove Theorem \ref{theo: rhoN<<3/7}
%and assume that $\pi_1$ lies in a compact interval of $(1,\infty)$. Furthermore,
%we also make the
under the following simplifying assumptions. We assume that $\pi_1>1$ is given independently of $N$ and that
$W_N=\text{ diag }(\pi_1,\ldots,\pi_1,0,\ldots,0),$
with $\pi_1$ of multiplicity $k_N$, for some sequence $k_N$ satisfying (\ref{condH2suites}).
Changes to be made in the case where $W_N$ has eigenvalues distinct of $0$ and $\pi_1,$ or to prove Theorem \ref{theopi_1<1},
 will be indicated in subsection \ref{subsec: extensionedge} below.
With the above assumption, $F_u$, defined by (\ref{defindefuengen}), becomes \be
F_u(w)=\frac{w^2}{2}-uw+(1-\alpha_N^2)\log (w)+\alpha_N^2 \log (w-\pi_1).\label{Fuaubord}\ee The
basic idea for the study of the correlation kernel at the edge is to perform a third order Taylor
expansion of $F_u$ close to the degenerate critical point $w_o$ defined by $F_u'(w_o)=F_u''(w_o).$
This point is close to $\pi_1+\alpha_N \sigma,$ which is the degenerate critical point of
$H_{2/\sigma}$. The ascent or descent curves for $F_u(\pi_1+\alpha_N t)$ should then be those for
$H_{2/\sigma}$ slightly modified in a neighborhood of width $k_N^{-1/3}$ of $\pi_1+\alpha_N
\sigma$, to go through the exact degenerate critical point. This simple analysis can be achieved
as long as $k_N<<N^{3/7}$. This is the regime where the bulk of $N-k_N $ eigenvalues does
not interfere with the $k_N$ largest eigenvalues. For the other regimes, one will have to define
new contours, that are descent or ascent paths for $F_u$, and show that the Taylor expansion can
still be made in a neighborhood of $w_o$. 
%This the reason why we have introduced the point $u_o$ defined in (\ref{defdeu_o}), and the associated degenerate critical point $w_o$.
We will however see that the asymptotic expansion is still lead in  some way by $H_{2/\sigma}$. \\

We set as in (\ref{defdeu_o}), $w_o=\pi_1+\alpha_N t_r$ and consider the
rescalings \be \label{rescalphan} u=u_o+x\left
(\frac{\nu_N}{2}\right)^{1/3}\frac{\alpha_N}{k_N^{2/3}} ,\quad v=u_o+y\left (\frac{\nu_N}{2}\right
)^{1/3}\frac{\alpha_N}{k_N^{2/3}}, \ee where \be \nu_N=\alpha_N F_{u_o}^{(3)}(\pi_1+\alpha_N
t_r)=
\frac{2}{t_r^3}+\alpha_N\frac{1-\alpha_N^2}{(\pi_1+\alpha_Nt_r)^3}.\label{nun}\ee

Let $\epsilon>0$ be given. From now on, we consider the rescaled correlation kernel
 \be \label{rescaledkernel}
K_N'(x,y)=\frac{\alpha_N}{k_N^{2/3}}\left (\frac{\nu_N}{2}\right )^{1/3}
K_N(u,v)\exp{\left\{-N(u-v)\left (\pi_1+\alpha_N(t_r+ \frac{\epsilon}{k_N^{2/3}})\right)\right \}}. \ee

The end of this section is now devoted to the proof of Theorem \ref{theo: rhoN<<3/7}. This proof is here indirect, since we will first split the correlation kernel into two subkernels. These subkernels are then analyzed separately, using the same scheme as in Section \ref{Sec: TW}.
\\
Before beginning the proof of Theorem \ref{theo: rhoN<<3/7}, it is convenient to make the
following assumption on $N$. Let then $t_c=\sigma$ be the degenerate critical point for $H_{2/\sigma}$ and
define sequences $\mu_N, \mu'_N$ by
$$u_o=C(\pi_1)+\alpha_N\frac{2}{\sigma} (1+\mu_N), \:t_r=t_c(1+\mu'_N).$$ Then it is easy to check
that there exists some constant $C$, depending on $\pi_1$ only, such that $ |\mu_N|,|\mu'_N|\leq C
\alpha_N$. Let also $R_o$ and $\nu_N$ be defined as in (\ref{defdeR_o}) and (\ref{nun}).
From now on, we assume that $N\geq N_o,$ where $N_o$ is such that
\begin{eqnarray}&\forall N\geq N_o, &\forall |t| \leq 2\sigma R_o+1, \:|\pi_1+\alpha_N t|\geq \frac{\pi_1}{2},\text{ and }|\pi_1+\alpha_N t-1|\geq \frac{\pi_1-1}{2}\cr
&& t_r \in [\frac{\sigma}{2},\frac{3\sigma}{2}],\: \: |\mu'_N|\leq \frac{1}{2},\:|\mu_N|\leq \frac{1}{2},\quad
\frac{3}{t_c^3}\geq \nu_N\geq \frac{1}{t_c^3}. \label{noedgefin}
\end{eqnarray}

\subsection{Rewriting the kernel}
 In this subsection, we split the kernel $ K_N'(x,y)$, defined in (\ref{rescaledkernel}), into two subkernels, to get rid of the integrals performed away from a small neighborhood of $\pi_1$. Then we bring these subkernels to the form (\ref{formeproduit}). Set $\tilde t_r=t_r+\dfrac{\epsilon}{k_N^{2/3}},$ and let $\exp{\{-N F_u(w)\}}$ stand for $\exp{\{-N
w^2/2+wu\}}\dfrac{w^{k_N-N}}{(w-\pi_1)^{k_N}}.$
Define the kernels
\begin{eqnarray} \label{newJn} &J_N(y)&=k_N^{1/3}(\frac{\nu_N}{2})^{1/3}\int_{\gamma'}\exp{\left\{N
F_{u_o}(\pi_1+\alpha _N t)\right\}}\exp{\left\{-k_N^{1/3}y(\frac{\nu_N}{2})^{1/3}(t-\tilde t_r)\right\}}dt, \\
&\label{newHn} H_N(x)&=k_N^{1/3}(\frac{\nu_N}{2})^{1/3}\int_{\Gamma_1'}\exp{\left\{-N
F_{u_o}(\pi_1+\alpha_Ns)\right \}}\exp{\left\{k_N^{1/3}x(\frac{\nu_N}{2})^{1/3}(s-\tilde t_r)\right\}}ds,\\
&\label{newHn''}H_N''(x)&=k_N^{1/3}(\frac{\nu_N}{2})^{1/3}\int_{\Gamma''}\exp{\left\{-N
F_{u_o}(\pi_1+\alpha_Ns)\right\}}\exp{\left\{k_N^{1/3}x(\frac{\nu_N}{2})^{1/3}(s-\tilde t_r)\right\}}ds,
 \end{eqnarray} where $\Gamma_1'$ is a contour encircling $0$ not crossing $\gamma':=a+i\mathbb{R}, a>0$ and $\Gamma''$ is such that its image under the map $t\mapsto\pi_1+\alpha_N t,$ is the circle of ray
one centered at the origin. Both $\Gamma_1'$ and $\Gamma''$ are oriented counterclockwise and $\gamma'$ is oriented from bottom to top.

\bp \label{Prop:bordnoyaucomplet} $K_N'(x,y)=K_N^1(x,y)+K_N^2(x,y),$ with
 $$K_N^1(x,y)=-\int_{0}^{\infty}H_N(x+u)J_N({y+u})du \text{ and }
 K_N^2(x,y)=-\int_{0}^{\infty}H_N''(x+u)J_N({y+u})du.$$
\ep

\paragraph{Proof of Proposition \ref{Prop:bordnoyaucomplet}: }
We first split the contour $\Gamma$ into the contours
$\displaystyle{\Gamma=\Gamma_1\cup\Gamma_2}$, where $\Gamma_1$ is encircling $\pi_1$ and crosses
the real axis at $\pi_1\pm \sigma \alpha_N.$ $\Gamma_2$ is a contour encircling $0$. Then,
let $\gamma=A+i\mathbb{R}$ with $A>0$ large enough so that $\gamma \cap \Gamma_1=\emptyset$ .
We call $K_N^1$ the part of the integral formula defining (\ref{rescaledkernel}) integrated on $\Gamma_1,$ and $ \gamma$. Then we obtain
\begin{eqnarray}
&K_N^1(x,y)&=\frac{N\alpha_N}{(2i\pi)^2k_N^{2/3}}(\frac{\nu_N}{2})^{1/3}\int_{\Gamma_1}dz\int_{\gamma}dw\cr
&&\frac{w^{N-k_N}(w-\pi_1)^{k_N}}{z^{N-k_N}(z-\pi_1)^{k_N}}\frac{1}{w-z}\frac{\exp{\left\{Nw^2/2-Nu_ow-N(v-u_o)(w-\tilde
\pi_1)\right\}}}{\exp{\left\{Nz^2/2-Nu_oz-N(u-u_o)(z-\tilde
\pi_1)\right\}}}\cr
&&=\frac{k_N^{2/3}}{(2i\pi)^2}(\frac{\nu_N}{2})^{2/3}\int_{\Gamma_1'}ds \int_{\gamma'}dt \int_{0}^{\infty}du
\exp{\left \{N F_{u_o}(\pi_1+\alpha_Nt)-N F_{u_o}(\pi_1+\alpha_Ns)\right \}}\cr
&&\exp{\left \{-k_N^{1/3}(y+u)(\frac{\nu_N}{2})^{1/3}(t-\tilde
t_r)+k_N^{1/3}(x+u)(\frac{\nu_N}{2})^{1/3}(s-\tilde
t_r)\right \}}.
\end{eqnarray}
The last equality follows from a change of variables.$\blacksquare$

We now set \be \label{ZN98}Z_N =\exp{\{N F_{u_o}(\pi_1+\alpha_N t_r) \}}.\ee The end of this section is aimed at obtaining the asymptotics of the rescaled kernels
$Z_N H_N'',$ $Z_NH_N,$ and $1/Z_N J_N$. It is then straightforward to deduce the asymptotics for the correlation kernel (\ref{rescaledkernel}).

\subsection{Estimate for$Z_NH_N''$}
The aim of this subsection is to prove the following Proposition. Let $H_N''$ be the kernel defined in (\ref{newHn''}), $Z_N$ as in (\ref{ZN98}).
\bp \label{Prop:loindepi_1}For any fixed $y_o \in \mathbb{R},$ $\exists$ $ C>0,c>0, C'>0$, an integer $N_o>0$ such that \be |Z_N
H_N''(x)|\leq \frac{C \exp{\{-cx\}}}{k_N^{1/3}}\exp{\{-C' N\}}, \text{ for any }x\geq y_o,\: N \geq
N_o.\ee  \ep
\paragraph{Proof of Proposition \ref{Prop:loindepi_1}: }Let $\Gamma''$ be such that its image under the map $\pi_1+\alpha_N t$ is the
circle of ray one, oriented counterclockwise. Then, it is easy to see that $\min_{\Gamma''}Re F_{u_o}(\cdot)=F_{u_o}(1).$ 
Now, one can check that 
$\displaystyle{F_{u_o}'(x)=-\frac{(x-\alpha_o)(x-(\pi_1+\alpha_Nt_r))^2}{x(\pi_1-x)},}$
where $\alpha_o<1$ is the second critical point, of mutliplicity one, of $F_{u_o}.$ 
Thus for $N$ large enough, as $\pi_1$ lies in a compact interval of $(1, \infty),$ one has that $Re F_{u_o}'(x)<0$ $\forall$ $x \in
(1,\pi_1)$.
Let then $0<\eta_1<\eta_2 <(\pi_1-1)/2$ be given and set $I=[1+\eta_1,1+\eta_2].$ Then, there exist $N_o$ and $\eta>0$, depending on $\pi_1$ only, such that 
 $|F_u'(x)|>2\eta,\: \forall x \in I \text{ and } \eta_2<\pi_1-\alpha_Nt_r, \: \forall N\geq N_o.$ From this, we deduce that there exists $\eta'>0$ such that $\displaystyle{\Big|\exp{\{-N
F_{u_o}(1)\}}\Big|\leq \Big|\exp{\{-N (F_{u_o}(\pi_1-\alpha_N t_r)+2\eta')\}}\Big|.}$ Now there exists $C>0$ such that $| F_{u_o}(\pi_1\!+\alpha_N t_r)-F_{u_o}(\pi_1\!-\alpha_N t_r)|\leq C\alpha_N,$ so that, for $N$ large enough,
$$\Big|\exp{\{-N
F_{u_o}(1)\}}\Big|\leq \!|\exp{\{-N F_{u_o}(\pi_1+\alpha_N t_r)-N\eta'\}}\Big |.$$
Let then $y_o>0$ be given and assume first that
$x \in[-y_o,y_o]$. Using (\ref{newHn''}) and the fact that the contour $\Gamma''$ is of length
$\frac{2\pi}{\pi_1 \alpha_N},$ we can see that for $N$ large enough,\\
$\displaystyle{|Z_N H_N''(x)| \leq \frac{k_N^{1/3}}{\alpha_N}\exp{\left \{\frac{k_N ^{1/3}}{\alpha_N}(\pi_1+2\sigma)y_o-N \eta'\right \}},}$
which goes to zero as $N$ goes to infinity (since $k_N^{1/3}/\alpha_N<<\sqrt N.$) Thus,
for $N$ large enough,
$\displaystyle{|Z_N H_N''(x)|\leq \exp{\{-N\frac{\eta'}{4}\}}.}$ This yields Proposition \ref{Prop:loindepi_1}
in this case. The case where $x$ is positive is handled as in the preceding sections. Indeed, $Re
(s-\tilde t_r)\leq -\frac{(\pi_1-1)}{4\alpha_N}-\epsilon$ along $\Gamma'',$ for $N$ large enough.
Thus, we readily obtain from the above proof that, for $x>0$ and $N$ large enough,
$\displaystyle{|Z_N H_N''(x)|\leq \exp{\{-N\frac{\eta'}{4}-\epsilon x\}}.\blacksquare}$
%This finishes the proof of Proposition \ref{Prop:loindepi_1}.$\blacksquare$

\subsection{\label{subsec: sousectiondernier}Estimate for $Z_NH_N$, $1/Z_N J_N$}The aim of this subsection is to obtain the following estimates for the kernels $H_N$ and $J_N$ defined in (\ref{newHn}) and (\ref{newJn}).

\bp \label{prop: bordGUEWig} Assume $\epsilon>0$ is fixed and let $\nu_N$ be given by (\ref{nun}),
$Z_N$ by (\ref{ZN98}).
 For any fixed $y_o\in \mathbb{R}$, $\exists \: C>0,c>0,N_o>0$ such that for any $y\geq y_o$, $x\geq y_o$ and $N\geq N_o,$
$$\Big|\frac{J_N(y)}{Z_N}-ie^{\epsilon y(\frac{\nu_N}{2})^{1/3}}Ai(y)\Big|\leq \frac{C\exp{\{-cy\}}}{k_N^{1/3}}\text{ and }
\Big|Z_NH_N(x)-ie^{-\epsilon x(\frac{\nu_N}{2})^{1/3}}Ai(x)\Big|\leq \frac{C\exp{\{-cx\}}}{k_N^{1/3}}.$$
\ep
The proof of Proposition \ref{prop: bordGUEWig} is divided into three parts. First, we establish three basic lemmas that enable us to get rid of some negligible parts of the contours and to perform the third order Taylor expansion. In the second part, we give the contours needed to perform the saddle point analysis and obtain, in the last part, the asymptotic expansion of the kernels $H_N$ and $J_N$.

\subsubsection{Preliminary lemmas \label{subsec: lemmas}}
In this part we prove that there exists a disk, $D=D(t_r, \delta'),$ such that the exponential term is driven by $H_{2/\sigma}$ outside $D$, and by its third order Taylor expansion inside $D.$
First, we fix the left frontier of $\Gamma_1$ and show that on this frontier, the exponential term behaves as $\exp{\{H_{2/\sigma}\}}$ despite the artificial singularity we have introduced (due to the $\log$). This is the object of the following Lemma.
\bl \label{bordenR} Let $R_o$ be defined in (\ref{defdeR_o}) and assume $t=\sigma (-R_o+ix), |x|\leq \sqrt 3.$
Then, $\exists C_o(\pi_1)>0$ depending on $\pi_1$ only such that
$$|\exp{\{NF_{u_o}(\pi_1+\alpha_N t_r)-N F_{u_o}(\pi_1+\alpha_N t)\}}|\leq
|\exp{\{k_N(H_{2/\sigma}(t_r)-H_{2/\sigma}(t))\}}|\exp{\{ C_o(\pi_1)\alpha_N k_N\}},$$
where $\exp{\{-k_NH_{2/\sigma}(-R_o)\}}$ stands for $\exp{\{-k_N \frac{R_o^2+4\sigma R_o}{2\sigma^2}\}}(-R_o)^{-k_N}.$\el
\paragraph{Proof of Lemma \ref{bordenR}: }
We set $t=-R+ix$ where $R=\sigma R_o$ and $x \in [0,\sigma \sqrt 3]$. The case where $x\in[-\sigma \sqrt 3,0]$ is obtained by using that $F_u(w)=\overline F_u(\bar w).$
As $N\geq N_o$, where $N_o$ has been defined in (\ref{noedgefin}), $\pi_1+\alpha_N t $ does not lie on the negative real axis, thus by a straightforward Taylor expansion
\begin{eqnarray}
&&\exp{\left \{N F_{u_o}(\pi_1+\alpha_N t_r)-N F_{u_o}(\pi_1+\alpha_N (-R+ix))\right \}}\cr
%&&=\exp{\left\{N\alpha_N (\pi_1-u_o)(t_r+R-ix)\right\}}
%\exp{\left\{k_N \left(\frac{t_r^2-(-R+ix)^2}{2}\right)\right\}}
%\left(\frac{t_r}{-R+ix}\right)^{k_N}\cr
%&&\times \Bigl (\frac{1+\alpha_N t_r/\pi_1}{1+\alpha_N (-R+ix)/\pi_1}\Bigr )^{N}\Bigl (\frac{1+\alpha_N t_r/\pi_1}{1+\alpha_N (-R+ix)/\pi_1}\Bigr )^{-k_N} \cr
&&=\exp{\left\{k_N \left(\frac{t_r^2-(-R+ix)^2}{2}\right)\right\}}
\left(\frac{t_r}{-R+ix}\right)^{k_N}\exp{\left\{-\frac{2k_N}{\sigma}(1+\mu_N)(t_r-(-R+ix))\right\}}\cr
&&\times\exp{\left\{N\alpha_N\left ( \int_{(-R+ix)/\pi_1}^{t_r/\pi_1}\frac{du}{1+\alpha_N u} -(t_r+R-ix)/\pi_1\right)\right\}}
\left (\frac{1+\alpha_N t_r/\pi_1}{1+\alpha_N (-R+ix)/\pi_1}\right )^{-k_N}\cr
&&\label{estbord}
\end{eqnarray}
where we have used that $u_o-C(\pi_1)=\alpha_N\dfrac{2}{\sigma}(1+\mu_N).$
Now
\begin{equation}
\int_{(-R+ix)/\pi_1}^{t_r/\pi_1}\dfrac{du}{1+\alpha_N u}-\dfrac{t_r+R-ix}{\pi_1}
%=-\alpha_N\int_{(-R+ix)/\pi_1}^{t_r/\pi_1}\frac{u}{1+\alpha_N u}du\cr
%&&
=-\frac{\alpha_N}{2}\frac{1}{\pi_1^2}(t_r^2-(-R+ix)^2)-\alpha_N^2\int_{(-R+ix)/\pi_1}^{t_r/\pi_1}\frac{u^2}{1+\alpha_N u}. 
\label{restebord}
\end{equation}
Inserting (\ref{restebord}) in (\ref{estbord}) yields
\begin{eqnarray}
&&\exp{\left\{N F_{u_o}(\pi_1+\alpha_N t_r)-N F_{u_o}(\pi_1+\alpha_N (-R+ix))\right\}}\cr
&&=\exp{\left\{k_N H_{2/\sigma}(t_r)-k_N H_{2/\sigma}(-R+ix)\right\}}\times \left (\frac{1+\alpha_N (-R+ix)/\pi_1}{1+\alpha_N t_r/\pi_1}\right )^{k_N}\label{resteparenthses}
\\
&& \times \exp{\left\{-\frac{2}{\sigma} \mu_N k_N(t_r+R-ix)-\alpha_N k_N \int_{(-R+ix)/\pi_1}^{t_r/\pi_1}\frac{u^2}{1+\alpha_N u}\right\}}\label{restepetit}
\end{eqnarray}
Now, as $N\geq N_o$, (\ref{restepetit}) is $O((\alpha_N+\mu_N) k_N)$, and this $O$ is uniform, since $\pi_1$ lies in a compact interval of $(1,\infty).$ Indeed, we can choose a segment $S$ for the $u$-path from $-R+ix $ to $t_r$, of length  smaller than $R^2+3\sigma^2+t_r^2\leq \sigma^2(R_o^2+3)+t_r^2$, which is uniformly bounded. Thus, as $t_r\in [\dfrac{\sigma}{2},\dfrac{3\sigma}{2} ]$, there exists $C_1(\pi_1, R_o)>0$ such that
$\displaystyle{\int_S\dfrac{|u|^2}{|1+\alpha_N u/\pi_1|}|du|\leq C_1(\pi_1, R_o).}$
The remaining bracket in (\ref{resteparenthses}) is obviously bounded. This finishes the proof of Lemma \ref{bordenR}.$\blacksquare$

\paragraph{}In the following lemma, we prove that, in a suitably chosen compact set of $\mathbb{C}$,
$NF_{u_o}(\pi_1+\alpha_N t)$ behaves, up to constants or lower order terms, as $k_N H_{2/\sigma}(t).$ Let $\delta'>0$ be given and
define
\begin{eqnarray}&& t_r^*(\Gamma'_1)=t_r(1+\delta')e^{2i\pi/3},\: t_c^*(\Gamma'_1)=t_c(1+\delta')e^{2i\pi/3},\label{lest*Gamma}\\
&& t_r^*(\gamma')=t_r(1+\delta')e^{i\pi/3},\: t_c^*(\gamma')=t_c(1+\delta')e^{i\pi/3}.\label{lest*gamma}
 \end{eqnarray}
Define also $D(\Gamma'_1)$ (resp. $D(\gamma')$) to be the segment joining $t_r^*(\Gamma'_1)$ to $t_c^*(\Gamma'_1)$ ( $t_r^*(\gamma')$ to $t_c^*(\gamma')$).
Let finally $R_o$ be chosen as in Lemma \ref{bordenR} and $\eta>0$ be given.

\bl \label{Lem: controle} There exists constants $Ct(\pi_1)$ depending on $\pi_1$ only, and $C>0$ (depending on $\eta $ and $\pi_1$) such that
\begin{eqnarray}&&|N F_{u_o}(\pi_1+\alpha_N t)-NCt(\pi_1)-k_N H_{2/\sigma}(t)|\leq C\alpha_Nk_N ,\quad \forall \: \eta<|t|\leq 2\sigma R_o,\cr
&&|N F_{u_o}(\pi_1+\alpha_N t_r )-NF_{u_o}(\pi_1+\alpha_N t_c)|\leq C \alpha_N k_N,\cr
&& |H_{2/\sigma}(t_r^*(\Gamma'_1))-H_{2/\sigma}(t)|\leq C\alpha_N , \: \forall t\in D_{\Gamma'_1},\text{ and } |H_{2/\sigma}(t_r^*(\gamma'))-H_{2/\sigma}(t)|\leq C\alpha_N , \forall t\in D_{\gamma'}.\nonumber \end{eqnarray}
\el
\paragraph{Proof of Lemma \ref{Lem: controle}:}One has 
$\displaystyle{\frac{d}{dt}NRe \left (F_{u_o}(\pi_1+\alpha_N t)\right)=k_N \left (Re (H_{2/\sigma }'(t)+\frac{2}{\sigma}\mu_N +\alpha_N G'(t))\right ).}$ The first estimate follows from the fact that $G$ and $H_{2/\sigma}$ are uniformly bounded in the annulus considered. 
Combining the first estimate and the inequality $|H_{2/\sigma}(t_r)-H_{2/\sigma}(t_c)|\leq \mu_N^3$ (which follows from the facts that $H_{2/\sigma}'(t)=\frac{(t-t_c)^2}{t\sigma^2}$ and $t_r, t_c$ are greater than
$\sigma/2$), yields the second estimate.
The last ones follow from the fact that both  $|t_c^*(\Gamma'_1)-t_r^*(\Gamma'_1)|\leq C'\alpha_N $ and $|t_c^*(\gamma')-t_r^*(\gamma')|\leq C'\alpha_N $ for some constant $C',$ and that $|H_{2/\sigma}'|$ is bounded on the two segments considered.
$\blacksquare$

\paragraph{} In the third lemma, we then determine a disk where the third order Taylor expansion for the exact exponential term $F_u(.)= F_{u,N}(.)$, depending on $N$, can still be made.
Let $\delta$ be given by (\ref{delta}).

\bl \label{Lem: Taylordelta'} There exist $0<\delta'<\delta/2 <1,$  $N_1$
independent of $\delta',$ a constant $C_o=C_o(\pi_1)>0,$ such that, for any $N \geq N_1$, for any $t\in D(t_r,
\delta'):=\{|t-t_r|\leq t_r \delta'\}$
\begin{eqnarray}
&&\!\!\!\!\!\!\!\!\!|F_{u_o}^{(4)}(\pi_1+\alpha_N t)\alpha_N^2|\leq C_o,\cr
&&\!\!\!\!\!\!\!\!\!\left|F_{u_o}(\pi_1+\alpha_N t)-F_{u_o}(\pi_1+\alpha_N t_r)-\frac{\alpha_N^3(t-t_r)^3}{3!}F_{u_o}^{(3)}(\pi_1+\alpha_N t_r)\right|\leq \frac{|\alpha_N(t_r-t)|^3}{24}|F_{u_o}^{(3)}(\pi_1+\alpha_N t_r)|\nonumber
\end{eqnarray}
\el

\brem \label{rem: taylordelta'}The above Lemma implies in particular, for $N$ large enough (to ensure that $\nu_N =F_{u_o}^{(3)}(\pi_1+\alpha_N t_r)\alpha_N\geq 1/t_r^3$), that
$\displaystyle{Re \left (NF_{u_o}(\pi_1+\alpha_N t_r^*(\gamma'))-NF_{u_o}(\pi_1+\alpha_N t_r)\right )\leq -k_N{\delta'}^3/8}$, and 
$\displaystyle{Re \left (NF_{u_o}(\pi_1+\alpha_N t_r^*(\Gamma_1'))-NF_{u_o}(\pi_1+\alpha_N t_r)\right )
\geq k_N{\delta'}^3/8.}$
\erem
\paragraph{Proof of Lemma \ref{Lem: Taylordelta'}:}
We prove the second inequality of Lemma \ref{Lem: Taylordelta'} (the first one will be established within this proof). This inequality will be established if we %prove
%For any $t \in D(t_r,\delta')$,
%$$\left|F_{u_o}(\pi_1+\alpha_N t)- F_{u_o}(\pi_1+\alpha_N t_r)-\frac{{F_{u_o}}^{(3)}(\pi_1+\alpha_N t_r)}{6}(\alpha_N (t-t_r))^3\right|\leq
%\frac{1}{4!}\max_{D(t_r, \delta')}|F_{u_o}^{(4)}(t)||\alpha_N (t_r-t)|^4.$$
%We then want to 
find $\delta'>0$ such that
$$\displaystyle{\frac{1}{4!}\max_{D(t_r, \delta')}\left|F_{u_o}^{(4)}(t)\right|\left|\alpha_N (t_r-t)\right|^4
\leq \frac{{F_{u_o}}^{(3)}(\pi_1+\alpha_N t_r)}{24}\left|\alpha_N (t-t_r)\right|^3.}$$
Assume $\delta'<1/2$, then, as $ t_r\in [\frac{\sigma}{2},\frac{3\sigma}{2}]$,
 $D(t_r,\delta')\subset D(t_c,\frac{\delta'+1}{2}).$
Define then $v_o=\dfrac{2}{\sigma}(1+\mu_N)$, so that $u_o=C(\pi_1)+\alpha_N v_o$, and
let $H_{v_o}$ be given by (\ref{defdeH}). Then,
$F_{u_o}'(\pi_1+\alpha_N t)=\alpha_N( H_{v_o}'(t)+\alpha_N G'(t)),$
where $\displaystyle{G'(t)=\frac{t^2-\pi_1^2}{\pi_1^2(\pi_1+\alpha_N t)}}.$ 
%and for $l=0,1$, $\displaystyle{
%F_{u_o}^{(3+l)}(\pi_1+\alpha_Nt)=\frac{1}{\alpha_N^{1+l}}\left(H_{v_o}^{(3+l)}(t)+\alpha_N G^{(3+l)}(t)\right).}$
Now, as $N\geq N_o$, for $t \in  D(t_r, \delta')\subset D(t_c,\frac{\delta'+1}{2})$, as $\pi_1+\alpha_N t\geq \pi_1/2,$
there exists constants $C_3(\pi_1)>0$, $C_4(\pi_1)>0$, depending on $\pi_1$ only, such that
$$\max_{t \in D(t_c,\frac{\delta'+1}{2})}|G^{(4)}(t)|\leq C_4(\pi_1),\: \: |G^{(3)}(t_r)|\leq C_3(\pi_1).$$  Note that this gives the first inequality in Lemma \ref{Lem: Taylordelta'} with $C_o=C_4(\pi_1).$
Furthermore, one has $\displaystyle{\max_{t\in D(t_r,\delta')} |H_{v_o}^{(4)}(t)|=\frac{6}{t_r^4(1-\delta')^4}.}$
Thus to prove Lemma \ref{Lem: Taylordelta'}, it is enough to determine $\delta'$ such that 
\be \label{condasat}\forall t\in D(t_r,\delta'),\:\:\frac{1}{4!} |t-t_r|^4\left (\frac{6}{t_r^4(1-\delta')^4}+\alpha_N C_4(\pi_1)\right)+\frac{\alpha_N}{24} C_3(\pi_1)|t-t_r|^3\leq \frac{|t-t_r|^3}{24} H_{v_o}^{(3)}(t_r).\ee
%Indeed $u_o$ is real so the third derivatives $H_{v_o}^{(3)}$ and $G^{(3)}$ are also real, and we can assume the worst case, that is, $G^{(3)}<0.$  \\
Let now $0<\delta'<1$ be such that
$\dfrac{\delta'}{(1-\delta')^4}< \dfrac{1}{32}.$ As 
$H_{v_o}^{(3)}(t_r)=2/t_r^3$, we then have that 
\be \frac{6}{24 t_r^4(1-\delta')^4}|t-t_r|^4<\frac{6}{32}|t-t_r|^3 \frac{H_{v_o}^{(3)}(t_r)}{24}, \quad \forall t\in D(t_r,\delta').
\label{frsteq}\ee
And there exists $N_2=N_2$, depending on $\pi_1$ only, such that, as $t_r \in [\frac{\sigma}{2}, \frac{3\sigma}{2}],$ and $\delta'<1,$
\be\label{lasteq}\frac{3!}{4!}\delta' t_r C_4(\pi_1)\alpha_N +\frac{\alpha_N}{24} C_3(\pi_1)\leq
\frac{\alpha_N}{24}\left(18\sigma C_4(\pi_1)+C_3(\pi_1)\right )\leq \frac{2}{96\sigma^3}\leq \frac{2}{96 t_r^3}=\frac{1}{4}\frac{H_{v_o}^{(3)}(t_r)}{24}. \ee
Formulas (\ref{frsteq}) and (\ref{lasteq}) now imply (\ref{condasat}).
This finishes the proof of Lemma \ref{Lem: Taylordelta'}.$\blacksquare$

\subsubsection{Contours}
\begin{figure}[htbp]
 \begin{center}
 \begin{tabular}{c}
 \epsfig{figure=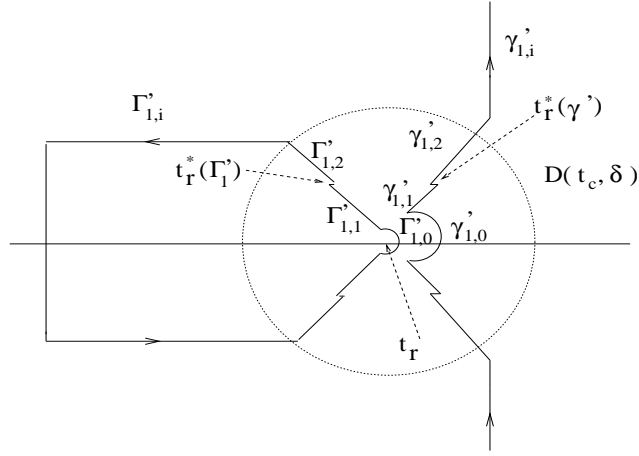, height=6cm,width=8.5cm, angle=0}
 \end{tabular}
 \caption{Contours $\Gamma'_1$ and $\gamma'$.
 \label{fig:contourGrhon}}
 \end{center}
\end{figure}

We now define the contours $\Gamma'_1$ and $\gamma',$ suitable for the saddle point analysis of $H_N$ and $J_N$.\\ %$\Gamma'_1$ and $\gamma '$ coincide with ascent and descent curves for $H_{2/\sigma}$ outside some disk $D(t_c, \delta)$, as we now explain.
Let $\delta$ be given by (\ref{delta}) and 
$\delta'\leq \delta/2$ be chosen so that Lemma \ref{Lem: Taylordelta'} holds. From now on, we assume that $N$ is large enough to ensure that $D(t_r,\delta')\subset D(t_c,\delta).$
Let then $\Gamma_{\sigma}$ and $\gamma_{\sigma }$ be the image of the contours defined in
Figure \ref{fig:contourGamma} under the map $t\mapsto \sigma t$.
Then
$\Gamma_{\sigma}$ (resp. $\gamma_{\sigma}$) is an ascent (resp. descent ) curve for
$H_{2/\sigma},$ as $H_{2/\sigma}(\sigma t)=F(t)+\log \sigma $ where $F$ has been defined in (\ref{defdeF}). \\
We now define the contour $\Gamma'_1$, which coincides with $\Gamma_{\sigma}$ outside $D(t_c,\delta)$.
Let then
$$\Gamma'_{1,i}=\Gamma_{\sigma}\cap D(t_c, \delta)^c;\: \:
\Gamma'_{1,0} = t_r +\frac{\epsilon}{2k_N ^{1/3}}e^{i\theta}, 0\leq \theta\leq 2\pi/3;\: \:
\Gamma'_{1,1}= t_r +te^{2i\pi/3}, \frac{\epsilon}{2k_N ^{1/3}}\leq t\leq \delta't_r.$$
Let then $t_r^*(\Gamma'_1)$ and $t_c^*(\Gamma'_1)$ be given as in (\ref{lest*Gamma}) and note that they are the respective endpoints of $\Gamma'_{11}$ and $\Gamma_{\sigma}.$
We then join $t_r^*(\Gamma'_1)$ to $t_c^*(\Gamma'_1)$ by a 
%contour within a $C^1$ distance of $\Gamma_{\sigma}$ 
segment (of length smaller than $ C \alpha_N$), and finally join $t_c^*(\Gamma'_1)$ to $t_c(1+\delta)e^{2i\pi/3}$ along $\Gamma_{\sigma}.$ We call $\Gamma'_{1,2}$ this last contour.
Finally we set $\Gamma'_{1}=\Gamma'_{1,i}\cup \Gamma'_{1,0}\cup \Gamma'_{1,1}\cup\Gamma'_{1,2}\cup\overline{\Gamma'_{1,0}\cup \Gamma'_{1,1}\cup\Gamma'_{1,2}},$ and this contour is oriented counterclockwise.
Similarly, $\gamma'$ is the contour $\gamma_{\sigma}$ modified in the disk $D(t_c, \delta)$, in the following way.
$$\gamma'_{1,i}=\gamma_{\sigma}\cap D(t_c, \delta)^c;\:\:
\gamma'_{1,0} = t_r +\frac{3\epsilon}{2k_N ^{1/3}}e^{i\theta}, 0\leq \theta\leq \pi/3;\: \:
\gamma'_{1,1}= t_r +te^{i\pi/3}, \frac{\epsilon}{2k_N ^{1/3}}\leq t\leq \delta't_r;
$$
Let also $t_r^*(\gamma')$ and $t_c^*(\gamma')$ be given by (\ref{lest*gamma}).
We then join $t_r^*(\gamma')$ to $t_c^*(\gamma')$ by a 
%contour within a $C^1$ distance of $\gamma_{\sigma }$ 
segment (of length smaller than $ C \alpha_N$), and finally join $t_c^*(\gamma')$ to $t_c(1+\delta)e^{i\pi/3}$ along $\gamma_{\sigma}$.
We call $\gamma'_{1,2}$ this last contour and define
$\gamma'=\gamma'_{1,i}\cup \gamma'_{1,0}\cup \gamma'_{1,1}\cup\gamma'_{1,2}\cup\overline{ \gamma'_{1,0}\cup \gamma'_{1,1}\cup\gamma'_{1,2}},$ oriented from bottom to top.
A plot of the contours $\Gamma'_1$ and $\gamma'$ is given on Figure \ref{fig:contourGrhon}.
\brem \label{rempoureta}There exists $\eta>0$ such that $\gamma'\cap D(0,\eta)=\emptyset$ and $\Gamma_1'\cap D(0,\eta)=\emptyset.$
\erem
The contours defined above coincide with the steepest ascent and descent curves for $F_{u_o}$ in a small disk $D(t_r, \delta')$, where the third order Taylor expansion is known to hold. Thus we now
introduce the expected limiting kernels. Let $\Gamma_{\infty,N}$ (resp. $\gamma_{\infty,N}$) be a contour such that it coincides with the  image of $\Gamma'_1$ (resp. $\gamma'$) under the map $t\mapsto
k_N^{1/3}(t-t_r)$, in the disk $D(t_r,\delta'),$ and then follows the
curve $te^{\pm i 2\pi/3},\:|t|\geq \delta',$ (resp. $te^{\pm i \pi/3},\: |t|\geq \delta'$). Set then
 \begin{eqnarray}
&&\label{newHinfini}H_{\infty,N }(x):=(\frac{\nu_N}{2})^{1/3}\exp{\{-\epsilon
x(\frac{\nu_N}{2})^{1/3}\}}\int_{\Gamma_{\infty,N}} \exp{\{x(\frac{\nu_N}{2})^{1/3}a-\frac{a^3}{3!}
\nu_N\}}da, \\
&&  \label{newJinfini}J_{\infty,N }(y):=(\frac{\nu_N}{2})^{1/3}\exp{\{\epsilon
y(\frac{\nu_N}{2})^{1/3}\}}\int_{\gamma_{\infty,N}} \exp{\{-y(\frac{\nu_N}{2})^{1/3}b-\frac{b^3}{3!}
\nu_N\}}db.              \end{eqnarray}
Then, 
$H_{\infty,N }(x)=i\exp{\{-\epsilon
x(\frac{\nu_N}{2})^{1/3}\}}Ai(x) \text{ and }J_{\infty,N }(y)=i\exp{\{\epsilon
y(\frac{\nu_N}{2})^{1/3}\}}Ai(y).$
We now split the contours. 
Set $H_N'(x)=H_N(x)-H_{N,2}''(x),$ $J_N'(y)=J_N(y)-J_{N,2}''(y),$ where
\begin{eqnarray}&&H_{N,2}''(x)=k_N^{1/3}(\frac{\nu_N}{2})^{1/3}\int_{\Gamma_1'\cap D(t_r,\delta')^c}\exp{\left\{-N
F_{u_o}(\pi_1+\alpha_Ns)\right \}}\exp{\left\{k_N^{1/3}x(\frac{\nu_N}{2})^{1/3}(s-\tilde t_r)\right\}}ds,\cr
&&J_{N,2}''(x)=k_N^{1/3}(\frac{\nu_N}{2})^{1/3}\int_{\gamma'\cap D(t_r,\delta')^c}\exp{\left\{N
F_{u_o}(\pi_1+\alpha _N t)\right\}}\exp{\left\{-k_N^{1/3}y(\frac{\nu_N}{2})^{1/3}(t-\tilde t_r)\right\}}dt.\nonumber\end{eqnarray} 
Similarly, $H''_{\infty,N }(x)$ (resp. $J''_{\infty,N }(y)$) is the part of (\ref{newHinfini}) (resp (\ref{newJinfini}) corresponding to the integral performed on the curve $te^{\pm i 2\pi/3},\:|t|\geq \delta'$ (resp. $te^{\pm i \pi/3},\:|t|\geq \delta'$).

\subsubsection{Saddle point estimates}
We now prove Proposition \ref{prop: bordGUEWig} in the case
$x$ and $y$ lie in a fixed compact interval; the case where they are positive follows from arguments similar to those of the preceding sections.\\
 We first show that the contribution of the contour outside $D(t_r, \delta')$ is negligible, because the exponential term behaves as $k_N H_{2/\sigma}(t)$ outside this disk.
\bfa \label{factZ_N"}
 Let $y_o>0$ be fixed and assume that $x,y$ $\in [-y_o,y_o]$.
There exists $N_1>0$ such that,
\begin{eqnarray}
 && |Z_N H_{N,2}''(x)| \leq \exp{\{-\frac{k_N }{16} \delta'^3\}}, \:|H''_{\infty,N }(x)|\leq \exp{\{-k_N \frac{\delta'^3}{12}\}}, \forall N \geq N_1,\label{HN2''}\\
&&
  |\frac{1}{Z_N} J_{N,2}''(y)| \leq \exp{\{-\frac{k_N }{16} \delta'^3\}},\:|J''_{\infty,N }(y)|\leq \exp{\{-k_N \frac{\delta'^3}{12}\}}, \forall N \geq N_1.\label{JN2''}
\end{eqnarray}
\efa
\paragraph{Proof of Fact \ref{factZ_N"}: }
We first prove (\ref{HN2''}) and consider $\Gamma'_1\cap D(t_r,\delta')^c$.
Lemma \ref{Lem: Taylordelta'} and Remark \ref{rem: taylordelta'} first ensure that
$\displaystyle{Re \left (NF_{u_o}(\pi_1+\alpha_N t_r^*(\Gamma'_1))\right )-NF_{u_o}(\pi_1+\alpha_N t_r) \geq k_N\frac{{\delta'}^3}{8}.}$
Let $\eta>0$ be chosen as in Remark \ref{rempoureta}. Then, from Lemma \ref{Lem: controle}, we obtain that $\forall t\in \Gamma'_{12}\cap \Gamma_{\sigma}^c,,$
\begin{eqnarray}
&  NRe \left (F_{u_o}(\pi_1+\alpha_N t)-F_{u_o}(\pi_1+\alpha_N t_r)\right )
&>\: NRe \left (F_{u_o}(\pi_1+\alpha_N t)-F_{u_o}(\pi_1+\alpha_N t_r^*(\Gamma'_1))\right ) +k_N \frac{{\delta'}^3}{8}\cr
&&\geq k_N \frac{{\delta'}^3}{8}-C\alpha_N k_N \geq k_N \frac{{\delta'}^3}{16}, \label{premdom}
\end{eqnarray}
for $N$ large enough.
Similarly for $t\in \Gamma'_{12}\cap \Gamma_{\sigma}$, using Lemma \ref{bordenR}, Lemma \ref{Lem: controle}, and the fact that $\Gamma_{\sigma}$ is an ascent curve for $H_{2/\sigma},$ we obtain that
$\displaystyle{ Re \left (NF_{u_o}(\pi_1+\alpha_N t)-NF_{u_o}(\pi_1+\alpha_N t_r)\right )\geq k_N \frac{{\delta'}^3}{16}}$
%k_N \frac{{\delta'}^3}{8}-C\alpha_N k_N \geq , \label{premdom2}\ee
for $N$ large enough. Combining the latter inequality and (\ref{premdom}), we obtain that\\
 $\displaystyle{|Z_N H_{N,2}''(x)| \leq \exp{\left \{-k_N  \frac{\delta'^3}{8} +k_N^{1/3} y_o+C\alpha_N k_N\right\}},}$
 for some constant $C$ uniformly bounded. Thus, for $N$ large enough, we obtain the first part of (\ref{HN2''}). The second part is straightforward using that $\nu_N \in [1/\sigma, 3/\sigma].$\\
 We now turn to the proof of (\ref{JN2''}). Remark \ref{rem: taylordelta'} also ensures that
 $$Re \left (NF_{u_o}(\pi_1+\alpha_N t_r^*(\gamma'))-NF_{u_o}(\pi_1+\alpha_N t_r)\right ) \leq -k_N\frac{{\delta'}^3}{8}.$$
Let then $t_o$ be chosen as in (\ref{gamma2}) and
large enough so that $(1+t_o/2)^2<3t_o^2/4.$ Using again Lemma \ref{Lem: controle}, we have that
$\forall t \in \gamma'\cap \{|Im (t)|<t_o\sqrt 3\sigma/2\}$
 \be Re \left (NF_{u_o}(\pi_1+\alpha_N t)-NF_{u_o}(\pi_1+\alpha_N t_r)\right )\leq
-k_N \frac{{\delta'}^3}{8}+C\alpha_N k_N \leq -k_N \frac{{\delta'}^3}{16}, \label{premdom3}\ee
for $N$ large enough. And for $t\in \gamma',$ with $ t=t(s)=t_c+t_oe^{i\pi/3}\sigma+is,$ $s\geq 0$, it is easy to check that, as
$(1+t_o/2)^2<3t_o^2/4,$ there exists $C>0$ depending on $\pi_1$ only, such that 
%if$N$ is large enough to ensure that $\alpha_N<2/Re (t_c+t_oe^{i\pi/3}\sigma), $ then,
\be
\label{premdom4}Re \frac{d}{ds}N F_{u_o}(\pi_1+\alpha_N t(s))<-k_N Im\left (H_{2/\sigma}(t(s))\right
)\leq -k_N CIm (t(s)).\ee Now, (\ref{premdom3}) and (\ref{premdom4})
give that
$\displaystyle{\Big |\frac{1}{Z_N} J_{N,2}''(y)\Big | \leq \exp{\left \{-k_N  \frac{\delta'^3}{8} +k_N^{1/3} y_o+C\alpha_N k_N\right\}},}$
which proves  the first part of (\ref{JN2''}). The second part of (\ref{JN2''}) is easy to check. $\blacksquare$

\paragraph{} We now show that the contribution from the contours in the disk $D(t_r,\delta')$ gives the leading term of the asymptotic expansion for both kernels $Z_N H_N$ and $1/Z_N J_N$.

\bfa \label{factZ8N'}Let $y_o>0$ be fixed and assume $x,y \: \in [-y_o, y_o].$ Then,  $\exists\: C=C(y_o)>0,\:  N_o ,\text{ such that } \forall N \geq N_o$, one has
\begin{eqnarray}&& \label{HN'}|Z_N H_N'(x)-H'_{\infty,N}(x)|\leq \frac{C }{k_N^{1/3}}, \quad
\frac{1}{Z_N} |J_N'(y)-J'_{\infty,N}(y)|\leq \frac{C }{k_N^{1/3}}.
\end{eqnarray}
\efa
\paragraph{Proof of Fact \ref{factZ8N'}: } We will only prove the first inequality of (\ref{HN'}), since the second follows from similar arguments. Then,
\begin{eqnarray}
&|Z_N H_N'(x)-H_{\infty,N }'(x)|\leq&\frac{k_N^{1/3}}{2\pi}\left (\frac{\nu_N}{2}\right )^{1/3}\int_{\Gamma'_{1,0}\cup\Gamma'_{1,1}}e^{k_N^{1/3}y_o(\frac{\nu_N}{2})^{1/3}
Re (t-\tilde t_r)}\cr
&&\times |e^{-NF_{u_o}(\pi_1+\alpha_N t)+NF_{u_o}(\pi_1+\alpha_N t_r)}-e^{-k_N \nu_N(t-
t_r)^3/3!}|dt|.\cr &&\label{127}
\end{eqnarray}

We first consider the $\Gamma'_{1,0}$ integral in (\ref{127}). Thus $t=t_r+\dfrac{\epsilon}{2k_N^{1/3}}e^{i\theta}$,
and using Lemma \ref{Lem: Taylordelta'} to mimick the proof of (\ref{158}), we obtain that
$$\Big|\exp{\{N F_{u_o}(\pi_1+\alpha_N t)-N F_{u_o}(\pi_1+\alpha_N t_r)\}}-\exp{\{-k_N \frac{\nu_N(t-t_r)^3}{3!}\}}\Big|\leq \tilde C_o \exp{\{\nu_N \epsilon^3\}}\frac{1}{k_N^{1/3}}.$$
Now, we use the fact that $\nu_N \leq \frac{3}{\sigma^3(\pi_1)}$, by (\ref{noedgefin}), to obtain that there exists $C>0$ so that in (\ref{127})
$$\frac{k_N^{1/3}}{2\pi}\left (\frac{\nu_N}{2}\right )^{1/3}\!\int_{\Gamma'_{1,0}}e^{k_N^{1/3}y_o(\frac{\nu_N}{2})^{1/3}
Re (t-\tilde t_r)}\Big|e^{-NF_{u_o}(\pi_1+\alpha_N t)+NF_{u_o}(\pi_1+\alpha_N t_r)}-e^{-k_N \nu_N(t-
t_r)^3/3!}\Big||dt|\leq \frac{C}{k_N^{1/3}}.$$

And for $t=t_r+pe^{i2\pi/3}\in \Gamma'_{1,1}$, there exists $C_o$, depending on $\pi_1$ only, such that, by Lemma \ref{Lem: Taylordelta'},
$$|\exp{\{NF_{u_o}(\pi_1+\alpha_N t_r)-NF_{u_o}(\pi_1+\alpha_N t)\}}-\exp{\{-k_N(t-t_r)^3\frac{\nu_N}{3!}\}}|
\leq \exp{\{-k_Np^3\frac{\nu_N}{4!}\}} C_o(k_Np^4+p).$$
%Furthermore, $\forall t\in \Gamma'_{1,1}$, $\Big |\exp{\{
%x\left(\frac{\nu_N}{2}\right)^{1/3}(t-\tilde t_r)\}}\Big |\leq \exp{\{\left(\frac{\nu_N}{2}\right)^{1/3}k_N^{1/3}y_op+\epsilon y_o
%\left(\frac{\nu_N}{2}\right )^{1/3}\} }$.
Now, following the same scheme as in Section \ref{Sec: TW}, we obtain that in (\ref{127})
\begin{eqnarray}
&&\frac{k_N^{1/3}}{2\pi}\left (\frac{\nu_N}{2}\right )^{1/3}\!\int_{\Gamma'_{1,1}}e^{k_N^{1/3}y_o(\frac{\nu_N}{2})^{1/3}
Re (t-\tilde t_r)}\Big|e^{-NF_{u_o}(\pi_1+\alpha_N t)+NF_{u_o}(\pi_1+\alpha_N t_r)}-e^{-k_N \nu_N(t-
t_r)^3/3!}\Big||dt|\cr
&&\leq C_ok_N^{1/3}\int_{\frac{\epsilon}{2k_N^{1/3}}}^{\delta' t_r}
(k_Np^4+p)\exp{\left \{\epsilon y_o\left(\frac{\nu_N}{2}\right )^{1/3}+\left(\frac{\nu_N}{2}\right)^{1/3}k_N^{1/3}y_o \frac{p}{2}-k_N\frac{\nu_N p^3}{4!}\right \}}\leq \frac{C}{k_N^{1/3}}.\nonumber\end{eqnarray}
Here, we have used that both $\nu_N$ and $t_r$ are uniformly bounded.
This finally gives from (\ref{127}) that
$\displaystyle{|Z_N H_N'(x)-H'_{\infty,N}(x)|\leq \frac{C}{k_N^{1/3}}.}$
This proves (\ref{HN'}). 
$\blacksquare$\\
Combining formulas (\ref{HN'}),(\ref{HN2''}) and  (\ref{JN2''}) yield then
Proposition \ref{prop: bordGUEWig} in the case $x$ or $y$ lie in a fixed compact interval. The
case where $x>0$ (resp. $y>0$), is analyzed in a similar way than in the preceding sections, using the fact that the
whole contour $\Gamma_1$ (resp. $\gamma'$) lies on the left (resp. right) handside of $\tilde t_r.$ The detail is left. This finishes the proof of Theorem \ref{theo: rhoN<<3/7}.

\subsection{Extensions\label{subsec: extensionedge}}
In this part, we explain how the proof has to be modified to consider more general diagonal perturbations $W_N$. It is easy to see that the core of the proof of Theorem \ref{theo: rhoN<<3/7} are the three Lemmas obtained in Subsection \ref{subsec: lemmas}.\\ %Once these lemmas proved, the end of the proof should follow, as long as the eigenvalues of $W_N$ remain in a fixed compact interval.\\
We now indicate the main changes to prove Theorem \ref{theo: rhoN<<3/7} under Assumption \ref{Hyp2}, when some eigenvalues of $W_N$ differ from $0$ or $\pi_1$.
Let $\tilde F_u$ be given by (\ref{tildefuext}) and $\tilde w_o$, $\tilde u_o$ be defined as in (\ref{defwophrase}) and (\ref{defdeu_o}).
Let also $F_u$ be as in (\ref{Fuaubord}) and set
set $\tilde G=\tilde F_u -F_u.$ Then, under assumption \ref{Hyp2}, there exist sequences $\mu'_N$, $\eta_N$, $\mu''_N$, a constant $C>0,$ such that
\begin{eqnarray}&\tilde w_o=\pi_1+\alpha_N \sigma(\pi_1)(1+\mu'_N +\eta_N),
&\text{ where }|\mu'_N|\leq C\alpha_N\text{ and }|\eta_N|\leq C\beta_N,\cr
&&\cr
&\tilde u_o =C(\pi_1)+\tilde G'(\pi_1) +\alpha_N \frac{2}{\sigma(\pi_1)} +C\alpha_N \mu''_N, &\text{ with }|\tilde G'(\pi_1)|\leq C\beta_N \text{ and }\lim_{N\rightarrow \infty}\mu''_N=0.\nonumber
\end{eqnarray}
This implies that $\tilde t_r =\dfrac{w_o-\pi_1}{\alpha_N}$ still lies in an arbitrarily small neighborhood of $t_c=\sigma(\pi_1)$ and also gives that
 $\displaystyle{0<\lim_{N \rightarrow \infty}\alpha_N{\tilde F_{\tilde u_o}}^{(3)}(\tilde w_o) <\infty.}$
And, given a compact set $K$ of $\mathbb{C}^*$, there exist positive constants $C_o,$ $C_1$, $C$, depending on $\pi_1$ and $K$, and a sequence $\mu_N$ with $\displaystyle{\lim_{N\rightarrow \infty}\mu_N=0},$ such that,
\begin{eqnarray}&&N\tilde F_{\tilde u_o}(\pi_1+\alpha_N t)=N C_o+r_N C_1+k_N H_{2/\sigma}(t) +O(\mu_N k_N),\quad \forall t \in K,\label{premlem}\\
&&|\tilde F_{\tilde u_o}^{(l)}(\pi_1+\alpha_N t)-F_{\tilde u_o}^{(l)}(\pi_1+\alpha_N t)|\leq C \beta_N,\quad \text{ for }l=3,4, \quad \forall t \in K.\label{deuxlem}
\end{eqnarray}
In this case, Lemma \ref{Lem: controle} (resp. Lemma \ref{Lem: Taylordelta'}) follows from (\ref{premlem}) (resp. (\ref{deuxlem})). We also choose $\Gamma''$ as in the proof of Fact \ref{claim H2}. The end of the proof is a simple rewriting of the arguments used in the preceding subsections.
This gives Theorem \ref{theo: rhoN<<3/7} in this case.
\paragraph{}We now indicate the idea of the proof of Theorem \ref{theopi_1<1}, when $0<\pi_1\leq 1$. For ease of explanatory, we here assume that $r_N=0$.
Assume first that $\pi_1<1.$ Then, the exponential to be considered is given by
$\displaystyle{F_{u_o}(w)=w^2/2-u_o w+(1-\alpha_N^2)\log w +\alpha_N^2 \log (w-\pi_1).}$
Let then $w_o$ and $u_o$ be defined as in (\ref{defwophrase}) and (\ref{defdeu_o}).
Then there exists some sequences $C_N$, $C'_N$, $\nu_N$ such that
\begin{eqnarray}
&&w_o=1+\alpha_N^2C_N,\quad u_o=2+\alpha_N ^2C'_N , \quad F_{u_o}^{(3)}(w_o)=\nu_N, \text{ with}\cr
&&\lim_{N\rightarrow \infty }C_N = C_o:=\dfrac{1}{2}\left (\dfrac{1}{(1-\pi_1)^2}-1\right ), \quad
\lim_{N\rightarrow \infty}C'_N=\frac{1}{1-\pi_1}-1,\quad
\lim_{N\rightarrow \infty}\nu_N=2.\nonumber
\end{eqnarray}
The function that now leads the exponential term is
\be\label{Fdernier}F(w)=\frac{w^2}{2}-2w+\log w,\ee
and given any compact set $K$ of $\mathbb{C}\setminus \{0, \pi_1\},$ we have that
\be \Big|F_{u_o}^{(l)}(w)-F^{(l)}(w)\Big |\leq C(K)\alpha_N^2, \quad \forall l=0,\ldots, 4,\label{controle<1}
\ee
where $C(K)$ depends on the compact set $K$ only. Formula (\ref{controle<1}) ensures that  Lemma \ref{Lem: Taylordelta'} can be established, as $F^{(3)}(2)=2>0.$ It also readily gives Lemma \ref{Lem: controle}. We choose the contours $\Gamma$ and $\gamma$ as in Section \ref{Sec: TW}, slightly modified in a small disk around $w_o.$
Then replacing, in the whole Section \ref{sec: preuvebordgrrang}, the function $H_{2/\sigma}$ with $F$ defined above, it is not hard to deduce Theorem \ref{theopi_1<1}.\\
If $\pi_1=1$ then there exists sequences $\mu_N, \mu'_N$ such that
\begin{eqnarray}
&& w_o=1+2^{-1/3}\alpha_N^{2/3}(1+\mu'_N), \quad u_o=2+3\alpha_N^{4/3}2^{-2/3}(1+\mu_N),\text{ with }\lim_{N \rightarrow \infty}\mu^{(')}_N=0,\cr
&& F_{u_o}^{(3)}(w_o)=\nu_N,\text{ with }\lim_{N \rightarrow \infty}\nu_N
=6.\end{eqnarray}
Let then $K$ be a given compact set of $\mathbb{C^*}.$
By as straightforward Taylor expansion, one has that
$F_{u_o}(1+x\alpha_N^{2/3})=Ct(N) +\alpha_N^2 H(x) -\alpha_N^2\log(1+\alpha_N^{2/3}x) +O(\alpha_N^{2/3}|x|)^4, \: \forall\: x \in K,$ where $Ct(N)$ depends on $N$ only and
\be H(x)=x^3/3-3x 2^{-2/3}+\log x\label{Hdernier}.\ee
The function $H$ admits the degenerate critical point $x_c=2^{-1/3},$ and  $H^{(3)}(x_c)=6,$ $H^{(4)}(x_c)=2^{1/3}\times 12.$
Set then $G(x)=\left(F_{u_o}(1+x\alpha_N^{2/3})-Ct(N) +\alpha_N^2 H(x)\right)/\alpha_N^2.$ 
Then there exists $C>0$ such that $\displaystyle{|G^{(l)}(x)|\leq C\alpha_N^{2/3}, \forall x \in K, \forall l=0, \ldots, 4.}$ 
This ensures that Lemma \ref{Lem: controle} can be established in a suitably chosen neighborhood of width $\alpha_N^{2/3}$ of $w_o.$ Lemma \ref{Lem: Taylordelta'} also holds in some disk centered at $1+\alpha_N^{2/3}x_c$ of ray $\delta'\alpha_N^{2/3},$ for some $\delta'>0.$ 
Now, the steepest descent and ascent curves for $H$ can be computed. Indeed, one can check that
 $\displaystyle{\frac{d}{dt}Re \left (H(x_c+te^{2i\pi/3})\right )=\frac{t^2(t^2-2x_ct+3x_c^2)}{t^2-x_ct+t^2}>0,}$ $ \forall t\not=0.$
Then, the contours for the saddle point analysis are chosen as follows. Here, for short, we do not make the change of variables $w\rightarrow 1+\alpha_N^{2/3}x$ to define the contours as in Subsection \ref{subsec: sousectiondernier}. Let $t_o>\delta'$ be given and define
\begin{eqnarray}
&&\Gamma_{1, +}=\{1+\alpha_N^{2/3}(x_c+te^{2i\pi/3}), 0\leq t\leq 2x_c\}\cup \{1+\alpha_N^{2/3}\sqrt 3 x_c e^{i\theta}, \pi/2\leq \theta \leq \pi\},\cr
&&\Gamma''=1/2 e^{i\theta}, 0\leq \theta\leq 2\pi,\cr
&&\gamma_+=\{1+\alpha_N^{2/3}(x_c+te^{i\pi/3}), 0 \leq t\leq 2t_o\}\cup\{1+\alpha_N^{2/3}(x_c+2t_oe^{i\pi/3})+it, t\geq 0\},\nonumber
\end{eqnarray}
and set $\Gamma_1=\Gamma_{1, +}\cup\overline{\Gamma_{1, +}},$ $\gamma=\gamma_+\cup\overline{\gamma_+}.$ We then slightly modify the contours  $\Gamma_1$ and $\gamma$ in a small neighborhood of width $\alpha_N^{2/3}$ of $w_o$, as in Subsection \ref{subsec: sousectiondernier}.
Then, considering the rescalings $u=u_o+\alpha_N^{4/3}k_N^{-2/3}y=u_o+N^{-2/3}y$, it is enough to replace $H_{2/\sigma}$ with $H$ defined in (\ref{Hdernier}) and $\alpha_N$ with $\alpha_N^{2/3}$ in the whole Section \ref{sec: preuvebordgrrang}. The fact that the contribution of $\Gamma''$ is negligible is also clear. This is because, far from $w=1$, the exponential term $F_u(\cdot)$ behaves as $F$ defined in (\ref{Fdernier}). The proof of Theorem 
\ref{theopi_1<1} is then straightforward.

\end{document}